\documentclass{lms}

\usepackage{amssymb} 
\usepackage{amsmath} 
\usepackage{graphicx} 
\usepackage{psfrag} 

\newcommand{\qed}{\null\hspace*{\fill}\endmark}

\newcommand{\secno}{\null\hspace*{\fill}}

\newcommand{\pa}{\partial}
\newcommand{\half}{{\textstyle\frac{1}{2}}}
\newcommand{\tsum}{\textstyle\sum}

\newcommand{\norm}[1]{\lvert#1\rvert} 

\newcommand{\abs}[1]{\lvert#1\rvert} 

\newcommand{\inn}[1]{\langle#1\rangle} 
\newcommand{\biginn}[1]{\big\langle#1\big\rangle} 

\renewcommand{\vec}[1]{\boldsymbol{#1}} 

\newcommand{\dd}{{\mathrm d}}  
\newcommand{\UU}{{\mathrm U}}  
\newcommand{\ww}{{\mathrm w}}  
\newcommand{\WW}{{\mathrm W}}  
\newcommand{\JJ}{{\mathrm J}}  

\newcommand{\RR}{{\mathbb R}}  
\newcommand{\CC}{{\mathbb C}}  
\newcommand{\HH}{{\mathbb H}}  

\newcommand{\MM}{{\mathbb M}}  
\newcommand{\FF}{{\mathbb F}}  

\newcommand{\Cc}{{\mathcal C}}  
\newcommand{\Ff}{{\mathcal F}}  
\newcommand{\Ll}{{\mathcal L}}  
\newcommand{\Jj}{{\mathcal J}}  
\newcommand{\Hh}{{\mathcal H}}  

\newcommand{\Qq}{{\mathcal Q}}  

\newcommand{\Rt}{\widetilde{\mathbb R}^3}  
\newcommand{\Rf}{\widetilde{\mathbb R}^4}  
\newcommand{\Cf}{\widetilde{\mathbb C}^4}  
\newcommand{\Mf}{\widetilde{\mathbb M}^4}  

\newcommand{\Ci}{{\mathbb C} \cup \{ \infty\} }  

\newcommand{\RP}{{\mathbb R}P} 
\newcommand{\CP}{{\mathbb C}P} 

\newcommand{\RH}{H} 

\newcommand{\cC}{{\mathrm C}}   
\newcommand{\tT}{{\mathrm t}}  

\newcommand{\ii}{{\rm i}}  
\newcommand{\jj}{{\rm j}}  
\newcommand{\kk}{{\rm k}}  

\newcommand{\ra}{\rightarrow}
\newcommand{\lra}{\longrightarrow}
\newcommand{\llra}{\longleftrightarrow}

\renewcommand{\phi}{\varphi}
\newcommand{\la}{\lambda}
\newcommand{\La}{\Lambda}
\newcommand{\na}{\nabla}
\newcommand{\al}{\alpha}
\newcommand{\be}{\beta}

\newcommand{\Ga}{\Gamma}

\newcommand{\si}{\sigma}
\newcommand{\Si}{\Sigma}

\newcommand{\wt}{\widetilde} 
\newcommand{\ov}{\overline} 
\newcommand{\wh}{\widehat} 

\newcommand{\Orthog}{\mbox{\rm O}} 
\newcommand{\SO }{\mbox{\rm SO}} 
\newcommand{\SL}{\mbox{\rm SL}} 
\newcommand{\SU}{\mbox{\rm SU}} 
\newcommand{\GL}{\mbox{\rm GL}} 
\newcommand{\C}{\mbox{\rm C}}   
\newcommand{\Spin}{\mbox{\rm Spin}}  
\newcommand{\Symp}{\mbox{\rm Sp}} 

\DeclareMathOperator{\grad}{grad} 

\newcommand{\const}{\text{\rm const.}} 

\renewcommand{\Re}{{\rm Re}\,} 
\renewcommand{\Im}{{\rm Im}\,} 

\newcommand{\nablacirc}{\overset{\circ}{\na}}

\psfrag{x_1}{$x_1$}
\psfrag{x_2}{$x_2$}
\psfrag{smallw}{$\vec{w}$}
\psfrag{you}{$\UU$}
\psfrag{r3subx}{$\RR^3_{\vec{x}}$}
\psfrag{dbyst}{$\pa/\pa t$}

\begin{document}

\setcounter{chapter}{19}
\setcounter{section}{0}
\setcounter{subsection}{0}
\setcounter{page}{1}

\def\thechapter{\Alph{chapter}}

\chapter*{Harmonic morphisms and shear-free ray congruences}


\begin{center} {\Large\bf P.\ Baird and J.C.\ Wood}\end{center}
\begin{center}
{\small Universit\'e de Bretagne Occidentale,
Facult\'e des Sciences, D\'epartement de Math\'ematiques\\
6 Avenue Le Gorgeu,  29285 Brest, France, and,} \\ 
{\small Department of Pure Mathematics, University of Leeds, 
Leeds LS2 9JT, G.B.} \\
{\small e-mail addresses: Paul.Baird@univ-brest.fr, j.c.wood@leeds.ac.uk}
\end{center}

\medskip

\begin{center}We describe the relationship between \emph{complex-valued harmonic morphisms from
Minkowski\\ $4$-space} and the \emph{shear-free ray congruences} of mathematical physics.
Then we show how a horizontally conformal
submersion on a domain of $\RR^3$ gives the boundary values at infinity of a complex-valued harmonic morphism
on hyperbolic $4$-space.\end{center}
  
\medskip

\section*{Introduction}

\noindent
A \emph{harmonic morphism} is a map that preserves Laplace's equation. 
More explicitly, a smooth map
$\phi:M \to N$ between Riemannian or semi-Riemannian manifolds
is called a harmonic morphism if its composition $f \circ \phi$
with any locally defined harmonic function on the codomain $N$
is a harmonic function on the domain $M$; it thus `pulls back' germs of harmonic functions to germs
of harmonic functions.   Harmonic morphisms are characterized as harmonic \emph{maps} which
are `horizontally weakly conformal', see below.

Our purpose here is to describe the relationship between \emph{complex-valued harmonic morphisms from
Minkowski $4$-space}
and the \emph{shear-free ray congruences} of mathematical physics.  This paper
is a revised and expanded version of Baird and Wood (1998); there are no plans to publish it.
It can be regarded as a supplement to the book (Baird and Wood 2003),
see\\ \centerline{{\tt http://www.amsta.leeds.ac.uk/Pure/staff/wood/BWBook/BWBook.html},}
which the reader may consult for background on harmonic morphisms and related topics;
references in boldface such as `Definition {\bf \ref{pre:def:conffoln}}\,', `({\bf \ref{pre:Bott}})'
refer to that book.  In particular, the reader may find it helpful to read
Chapters {\bf\ref{ch:pre}}--{\bf\ref{ch:fund}} of the book which give a self-contained introduction to harmonic maps and
morphisms between Riemannian manifolds, Chapter {\bf\ref{ch:twistor}} which gives the background in twistor theory and its relationship to harmonic
morphisms in the Riemannian setting, and Chapter {\bf 14} which discusses harmonic maps
and morphisms between semi-Riemannian manifolds.

In Sections \ref{cxha:sec:conffoln}--\ref{cxha:sec:unification} of the present paper
we describe one-to-one correspondences between \emph{conformal
foliations by curves of\/ $\RR^3$}, \emph{Hermitian structures on\/ $\RR^4$},
\emph{shear-free ray congruences 
on Minkowski $4$-space} and \emph{complex-analytic foliations by null planes on $\CC^4$}; then, in Sections
\ref{cxha:sec:twistor}--\ref{cxha:sec:Kerr} we 
explain how to find these quantities by twistor theory, on the way discussing group actions on the twistor
spaces.   A CR interpretation is given in Section \ref{cxha:sec:CR}.  Although this material is essentially known
in mathematical physics, we give a
detailed self-contained treatment as it is not widely known outside that domain.

In Sections \ref{cxha:sec:loccoords}--\ref{cxha:sec:ha-morph-SFR}, we show how a complex-valued
harmonic morphisms from Minkowski $4$-space determines a shear-free ray congruence and conversely, and give
a version relating \emph{complex-harmonic morphisms} and null planes.

In Sections \ref{cxha:sec:hypha}--\ref{cxha:sec:examples} we show how a horizontally conformal
submersion on a domain of $\RR^3$ gives the boundary values at infinity of a complex-valued harmonic morphism on an open subset of
hyperbolic $4$-space; then we give explicit constructions of such submersions, relating these to the 
quantities described above.

We follow the format of the book (Baird and Wood 2003) by giving some `Notes and comments' on each section
at the end of the paper.  Conventions follow that book, as does notation,
as much as possible.  For a regularly updated list of papers on harmonic morphisms,
the reader should consult the \emph{Bibliography of Harmonic Morphisms} at\\
\centerline{\tt http://www.maths.lth.se/matematiklu/personal/sigma/harmonic/bibliography.html.}

\bigskip

\section{Conformal foliations by curves on\/ $\RR^3$}
\label{cxha:sec:conffoln}
Let $\RR^3$ denote three-dimensional Euclidean space, i.e.,
$\RR^3 = \{(x_1,x_2,x_3): x_i \in \RR\}$ equipped with its standard metric
$g = \dd {x_1}^{\! 2} + \dd {x_2}^{\! 2} + \dd {x_3}^{\! 2}$ and orientation.
Let $\Cc$ be a smooth foliation by curves of an open subset $A^3$ of $\RR^3$.
Assume that $\Cc$ is oriented, i.e., we can find a
smooth unit vector field $\UU$ tangent to its leaves.
The orthogonal complement $\UU^{\perp}$ is then oriented.

For such a foliation we can reformulate the definition of conformal
foliation (Definition {\bf \ref{pre:def:conffoln}}), as follows.
Let $\JJ^{\perp}$ denote rotation through
$+\pi/2$ on $\UU^{\perp}$, and let $\{ \ \}^{\perp}$ denote orthogonal projection onto $\UU^{\perp}$.
Let $\Ll_Y$ denote the Lie derivative with respect to a vector field $Y$. 

\begin{proposition} \label{cxha:prop:conffoln}
The foliation\/ $\Cc$ is conformal if and only if, at all points of\/ $A^3$,
\begin{equation} \label{cxha:conffoln1}
\bigl\{({\Ll}_{\UU} \JJ^{\perp})(X) \bigr\}^{\perp} = 0  \qquad  (X\in \UU^{\perp} ),
\end{equation}
and this holds if and only if
\begin{equation} \label{cxha:conffoln2}
\na_{\JJ^{\perp} X}\UU = \JJ^{\perp} \na_X \UU \qquad (X \in \UU^{\perp} ) \,.
\end{equation}
\end{proposition}

\begin{proof}
That \eqref{cxha:conffoln1} is equivalent to conformality follows from Proposition {\bf\ref{pre:prop:conf-codim2}}(ii).
To show the equivalence of \eqref{cxha:conffoln1} and \eqref{cxha:conffoln2}, we have,
for $X \in \Ga(\UU^{\perp})$,
\begin{eqnarray}
({\Ll}_{\UU} \JJ^{\perp})(X)  & = & \{ ( {\Ll}_{\UU}(\JJ^{\perp}X) \}^{\perp} -
\JJ^{\perp} \{ {\Ll}_{\UU}(X) \}^{\perp} \nonumber \\
& = & \{\na_{\UU}(\JJ^{\perp}X) \}^{\perp} - \na_{\JJ^{\perp}X}\UU -
\JJ^{\perp} \{ \na_{\UU} X \}^{\perp} + \JJ^{\perp}\na_X \UU
\hspace{1cm} \label{cxha:LieJ} \end{eqnarray}
noting that $\na_X \UU \in \Ga(\UU^{\perp})$
since $g(\na_X \UU, \UU ) = \half X \bigl( g(\UU,\UU)\bigr)  = 0$.  Further, since
$\UU^{\perp}$ has rank $2$, from Proposition {\bf\ref{pre:prop:conf-codim2}}(i), we have
$\na_{\UU}^{\mbox{\scriptsize{End}}\,\UU^{\perp}}\JJ^{\perp} = 0$ so that
$$
\{ \na_{\UU}(\JJ^{\perp}X) \}^{\perp} -\JJ^{\perp}\{\na_{\UU} X
\}^{\perp} = (\na_{\UU}^{\mbox{\scriptsize{End}}\,\UU^{\perp}}
\JJ^{\perp})(X) = 0 \;,
$$
hence (\ref{cxha:conffoln1}) holds if and only if (\ref{cxha:conffoln2})
holds.
\end{proof}

\begin{remark}  \rm  
In terms of the \emph{Bott partial connection} $\nablacirc _{\UU}$ on ${\rm End}(\UU^{\perp})$
(see ({\bf \ref{pre:Bott}}) and ({\bf\ref{pre:prop:conf-codim2}})),
\eqref{cxha:conffoln1} can be written as $\nablacirc_{\UU}\JJ^{\perp} = 0$.
\end{remark}

Recall that a $C^1$ map $M^m \to N^n$ between Riemannian or semi-Riemannian manifolds
is called \emph{horizontally (weakly) conformal} if, at each point of $M^m$,
the \emph{adjoint} of its differential is zero or conformal,
see Section {\bf \ref{pre:sec:HWC}} for equivalent and more geometrical ways of saying this.
In particular, a smooth map $f:M^m \to \RR^n$, $f =(f_1,\ldots, f_n)$, is horizontally weakly conformal if
the vector fields $\grad f_i$ are all orthogonal and of the same square norm:
$$
\inn{\grad f_i\,,\,\grad f_i} = \inn{\grad f_j\,,\,\grad f_j} \quad \mbox{and} \quad
\inn{\grad f_i\,,\, \grad f_j} = 0 \qquad \bigl(i,j \in \{1,\ldots, n\}, \ i \neq j \bigr) \,.
$$
(Here, $\inn{\ , \ }$ denotes the inner product on $TM$ given by the metric on $M$.)
For a smooth function $f:A^3 \to \CC$ this condition can be written nicely as
\begin{equation} \label{cxha:conffoln-f}
\inn{\grad f \,,\, \grad f} \equiv \sum_{i=1}^3 \left( \pa f/\pa x_i\right)^2 = 0 \qquad \bigl( (x_1,x_2,x_3) \in A^3 \bigr) \,.
\end{equation}

By Corollary {\bf \ref{pre:cor:folHC}}, conformal foliations
$\Cc$ are given locally as the fibres of horizontally conformal
submersions; in the case of interest this can be said explicitly as follows.

\begin{lemma} \label{cxha:lem:conffolnfn}
Let\/ $f:A^3 \to \CC$ be a smooth (respectively, real-analytic) function on an
open subset of\/ $\RR^3$ which satisfies \eqref{cxha:conffoln-f}
Then its fibres form a smooth (respectively, real-analytic) conformal foliation
$\Cc$ by curves.  All such foliations are given this way locally.
\qed \end{lemma}

\begin{remark} \label{cxha:rem:conffolnfn}  \rm
The foliation $\Cc$ is oriented with unit positive tangent
\begin{equation} \label{cxha:f-to-U}
\UU = \grad f_1 \times \grad f_2 \big/ \abs{\grad f_1 \times \grad f_2} \,.
\end{equation}
\end{remark}

\begin{example} \label{cxha:ex:linearfn}
Set $f(x_1,x_2,x_3) = a_1 x_1 + a_2 x_2 + a_3 x_3$ where the $a_i$ are
complex constants satisfying
${a_1}^{\! 2} + {a_2}^{\! 2} + {a_3}^{\! 2} = 0$;  this gives a foliation by parallel lines.
\end{example}

\begin{example} \label{cxha:ex:bunch}
For any constant $c \in \RR$, define $f:\RR^3 \setminus \{x_1\mbox{-axis}\}
\to \CC$ by
$$
f(x_1,x_2,x_3) = \big\{ (x_1+c)^2 + {x_2}^{\! 2} + {x_3}^{\! 2} \big\} \big/
	\{ x_2 + \ii x_3 \} \,.
$$
Then $f$ is horizontally conformal and submersive.  On identifying the extended
complex plane $\Ci$ with the $2$-sphere by stereographic projection ({\bf \ref{Euclid3:stereo0}}),
$f$ can be extended to
a horizontally conformal submersive map from $\RR^3 \setminus \{(-c,0,0)\}$ to
$\Ci \cong S^2$.
It can easily be checked that its level sets are the circles through
$(-c,0,0)$ tangent to the $x_1$-axis, together with the $x_1$-axis; these
are thus the leaves of a conformal foliation $\Cc_c$ on $\RR^3 \setminus \{(-c,0,0)\}$.
For later reference, the unit tangent to these leaves is calculated to be
\begin{equation} \label{cxha:bunch}
\UU = \frac{1}{\wt{x}_1^{\: 2}+{x_2}^{\! 2}+{x_3}^{\! 2}} \big({\wt{x}_1}^{\; 2}-{x_2}^{\! 2}-{x_3}^{\! 2}, \,2\wt{x}_1 x_2,
\,2\wt{x}_1 x_3 \big)\,;
\end{equation}
where $\wt{x}_1 = x_1 + c$.
\end{example}

For more examples, see Section \ref{cxha:sec:SFR}.

\section{Hermitian structures on Euclidean\/ $4$-space}
\label{cxha:sec:Herm}
Let $\RR^4$ denote four-dimensional Euclidean space, i.e.,
$\RR^4 = \{(x_0,x_1,x_2,x_3):x_i \in \RR \}$ equipped with the standard
Euclidean metric $g = \dd {x_0}^{\! 2} + \dd {x_1}^{\! 2} + \dd {x_2}^{\! 2} + \dd {x_3}^{\! 2}$
and standard orientation.
Let $\vec{x} \in \RR^4$.  Recall (see, e.g., Section {\bf \ref{twistor:sec:twistorspace}})
that an \emph{almost Hermitian structure at $\vec{x}$} is an isometric linear
transformation $J_{\vec{x}}: T_{\vec{x}} M \to T_{\vec{x}} M$ which satisfies ${J_{\vec{x}}}^{\! 2} = -I$ (where $I$ denotes the identity map).
We can always find an orthonormal basis $\{e_0, e_1, e_2, e_3\}$ of
$T_{\vec{x}}\RR^4$ such that $J_{\vec{x}} e_0 = e_1, J_{\vec{x}} e_2 = e_3$.
Then $J_{\vec{x}}$ is called
\emph{positive} (respectively, \emph{negative}) according as
$\{e_0, \ldots,e_3\}$ is positively (respectively, negatively) oriented.
As explained in Section {\bf \ref{twistor:sec:twistorspace}}, the set of all positive
Hermitian structures at a point can be identified with $S^2 \cong \CP^1$.

Let $A^4$ be an open subset of $\RR^4$. Recall (see, e.g., Section
{\bf\ref{twistor:sec:twistorspace}} again) that by an \emph{almost Hermitian
structure $J$ on $A^4$} we mean a smooth choice of almost Hermitian structure
at each point of $A^4$, i.e., a smooth map (which we shall denote by the
same letter)
$J:A^4 \to \CP^1$.  Recall, also, that $J$ is integrable if and only if its
Nijenhuis tensor $N$ satisfies
\begin{equation} \label{cxha:Jint1}
N \equiv 0 \,;
\end{equation}
by Proposition {\bf\ref{twistor:prop:JsigmaJ}}(ii) this is equivalent to
$(\na_{JX}J)(JY) = (\na_X J)(Y)$ \ $(\vec{x} \in A^4, \ X, Y \in T_{\vec{x}}\RR^4)$.
A simple
calculation shows that $(\na_X J)(Y) = J(\na_X J)(JY)$, so that \eqref{cxha:Jint1} is equivalent to
\begin{equation} \label{cxha:Jint2}
\na_{JX}J = J\na_X J \qquad (\vec{x} \in A^4, \ X \in T_{\vec{x}}\RR^4),
\end{equation}
explicitly,
$$
(\na_{JX}J)(Y) = J\big((\na_{X}J)(Y)\big) \qquad
(\vec{x} \in M, \ X,Y \in T_{\vec{x}} M) \,.
$$
If $A^4$ is given the almost Hermitian structure $J$, and
$\CP^1$ is given its standard K\"ahler structure,
then (\ref{cxha:Jint2}) is the condition that the map $J:A^4 \to \CP^1$
be holomorphic.   An integrable almost Hermitian structure is called
a \emph{Hermitian structure}.

\section{Shear-free ray congruences} \label{cxha:sec:SFR}
Let $\MM^4$ denote four-dimensional Minkowski space $\MM^4
= \{(t,x_1,x_2,x_3):t,x_i \in \RR \}$ equipped with the standard
Minkowski metric
$g^M = -\dd t^2 + \dd {x_1}^{\! 2} + \dd {x_2}^{\! 2} + \dd {x_3}^{\! 2}$
and time- and space-orientations.
 By an \emph{$\RR^3$-slice} we
mean a hyperplane given by $t = \mbox{const}$.  Let ${\vec{x}} \in \MM^4$.  We denote the
$\RR^3$-slice through ${\vec{x}}$ by $\RR^3_{\vec{x}}$.  Recall (e.g., from Definition {\bf \ref{semi:def:tns}})
that a non-zero vector $w \in T_{\vec{x}}\MM^4$ is called \emph{null} if $g^M(w,w) = 0$.
It spans a one-dimensional null subspace, we shall call such a subspace a \emph{null direction}.
Any non-zero null vector $w \in T_{\vec{x}}\MM^4$ can be normalized so that it is of
the form $\ww = \pa/\pa t +\UU$ where $\UU$ is a unit vector in $T_{\vec{x}}\RR^3_{\vec{x}}$.
Let $\ell$ be a smooth foliation by null lines (often called
a \emph{ray congruence}) of an open subset $A^M$ of
Minkowski space, let $\ww = \pa/\pa t + \UU$ be its tangent vector field and
write $\WW = \mbox{span}\{\ww\}$.   Then $\WW$ is a field of null directions.
The distribution $\WW^{\perp}$ orthogonal to $\WW$
(with respect to the Minkowski metric $g^M$) is
three-dimensional and contains $\WW$.  Choose any complement
$\Si$ of $\WW$ in $\WW^{\perp}$, such a complement is called a
\emph{screen space}, then the restriction of the Minkowski metric $g^M$
to $\Si$ is positive definite.   Note that any screen space is naturally isomorphic to the factor space
$\WW^{\perp}\big/\WW$, however, it is usually more convenient to work with a specific screen space.
A special choice of screen space at $\vec{x}$ is $\UU_{\vec{x}}^{\perp} \cap \RR^3_{\vec{x}}$;
this has a canonical orientation from which $\WW^{\perp}\big/\WW$, and so all screen spaces,
acquire an orientation.
  
Note that, for all  $X \in \Ga(\Si)$, we have
$\inn{\na_X \ww, \ww} = \half X \inn{\ww, \ww} = 0$; also
since the integral curves of $\ww$ are geodesics, we have
\begin{equation} \label{cxha:na_w-w}
\na_{\ww} \ww = 0 \,.
\end{equation}
It follows that
$\inn{\Ll_{\ww}X, \ww}
= \inn{\na_{\ww} X, \ww} - \inn{\na_X \ww, \ww} = -\inn{X, \na_{\ww}\ww} - 0 = 0$,
so that
\begin{equation} \label{cxha:welldef}
{\rm (i)} \ \na_X \ww \in \Ga(\WW^{\perp}) \quad \text{and}
	\quad {\rm (ii)} \ \Ll_{\ww}X \in \Ga(\WW^{\perp})
	\qquad \bigl( X \in \Ga(\Si) \bigr). 
\end{equation}

Let $J^{\perp} \in \Ga(\mbox{End}\,\Si)$ denote rotation through $+\pi/2$.
Then the Lie derivative of $J^{\perp}$ along $W$ is measured by
\begin{equation} \label{cxha:LJ}
\bigl\{ (\Ll_{\ww} J^{\perp}) (X)\bigr\}^{\Si}
	= \bigl\{ \Ll_{\ww} (J^{\perp}X)\bigr\}^{\Si} - \bigl\{ J^{\perp}(\Ll_{\ww}X)\bigr\}^{\Si}
	\qquad \bigl( X \in \Ga(\Si) \bigr).
\end{equation}
where $\{ \ \}^{\Si}$ indicates projection from $\WW^{\perp}$ onto $\Si$
along $\WW$;  this is well-defined by \eqref{cxha:welldef}(ii).

Then $\ell$ (or $\WW$) is said to be a 
\emph{shear-free ray (SFR) congruence} if Lie transport along $\WW$ of
vectors in $\Si$ is conformal, i.e.,
\begin{equation} \label{cxha:SFR1}
\bigl\{ (\Ll_{\ww} J^{\perp}) (X)\bigr\}^{\Si} = 0 \qquad \bigl( X\in \Ga (\Si )\bigr)\,,
\end{equation}
{}From \eqref{cxha:LJ},  we see that this is equivalent to
\begin{equation} \label{cxha:SFR2} 
\{\na_{J^{\perp}X}\ww\}^{\Si} = J^{\perp}\{\na_X \ww\}^{\Si}
	\qquad \bigl( X \in \Ga(\Si) \bigr) \,;
\end{equation}
it is easily checked that this condition is independent of the
choice of screen space. 
Comparison of equations (\ref{cxha:SFR2}) and (\ref{cxha:conffoln2})
shows that the
restriction of $\UU$ to any $\RR^3$-slice is a vector field
whose integral curves form an oriented conformal foliation by curves
of an open subset $A^3$ of the slice.  We shall call this process
\emph{projection onto the slice}.
Conversely, given an oriented conformal foliation of an open subset
$A^3$ of an $\RR^3$-slice, set $\UU$ equal to its unit positive tangent vector
field, then the null lines of $\MM^4$ tangent to
$\pa/\pa t + \UU$ at points of $A^3$ define a ray congruence
$\ell$ on an open neighbourhood of $A^3$ in $\MM^4$ which is shear-free
at points of $A^3$.  The next result shows that $\ell$ is a SFR
congruence; we shall call it the \emph{extension} of $\Cc$.

\begin{figure}[htb]
\centerline{\includegraphics*[bb=70 300 540 530,scale=0.5]{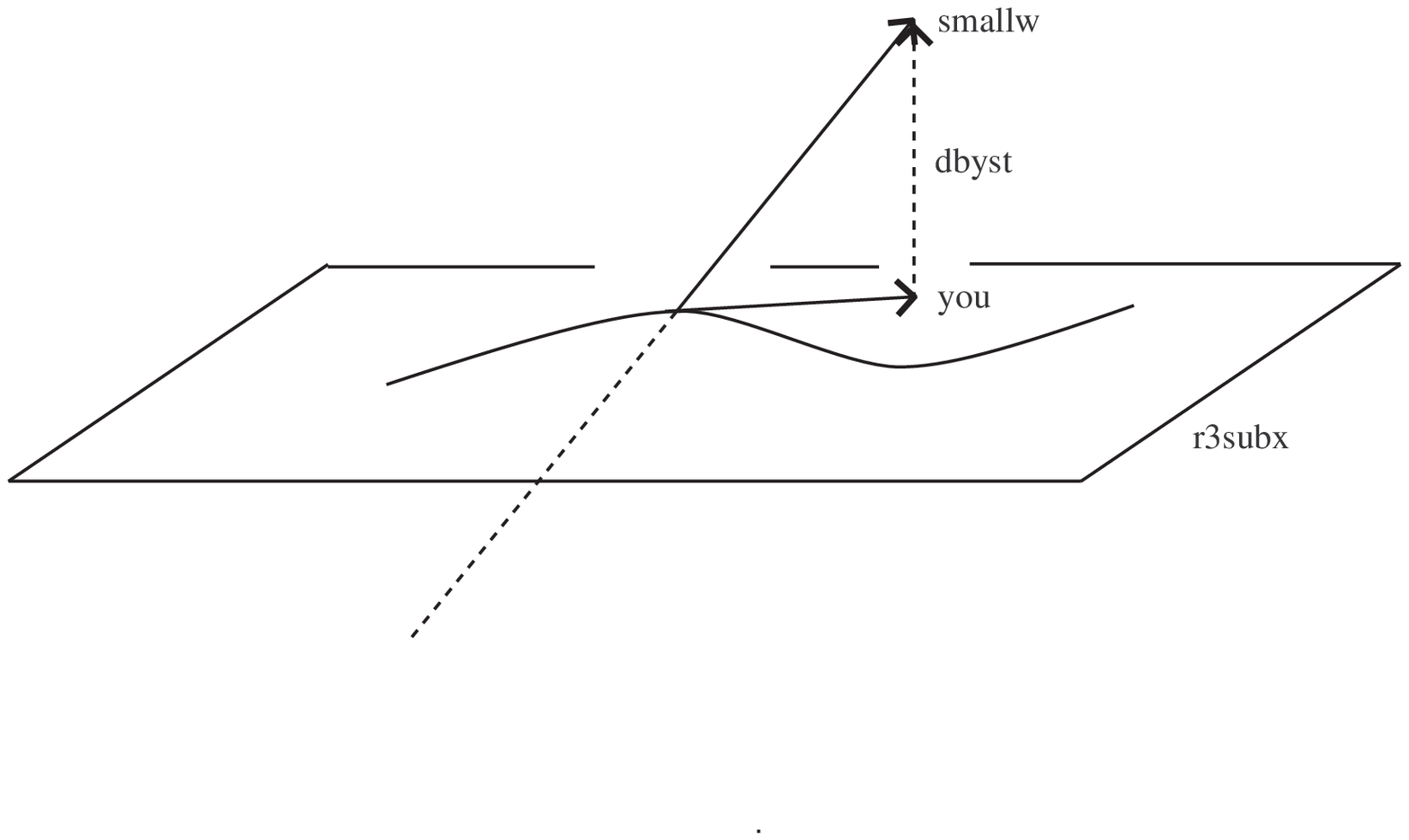}}
\caption{Projection onto a slice.}
\label{cxha:fig:projection}
\end{figure}

\begin{lemma} \label{cxha:lem:Sachs}
If a ray congruence\/ $\ell$ is shear-free at one point of a ray, it is
shear-free at all points of that ray.
\end{lemma}

\begin{proof}
Let $X$ be a section of $\Si$ and set $Z = X + \ii J^{\perp}X$.  Then it is
easily seen that, at any point, condition (\ref{cxha:SFR2}) is
equivalent to
\begin{equation} \label{cxha:SFR3}
\inn{\na_Z \ww, Z} = 0 \,.
\end{equation}
As in \eqref{cxha:welldef},
$\Ll_{\ww} Z \equiv [\ww,Z] \in \Ga(\WW^{\perp})$.  Extending the
terminology of Section {\bf \ref{pre:sec:confoln}}, 
call $Z$ \emph{basic} if the component of $\Ll_{\ww} Z$ in $\Si$ is zero, i.e.
\begin{equation} \label{cxha:Lie-vZ}
[\ww,Z] \in \Ga(\WW) \,.
\end{equation}
Then, using (\ref{cxha:na_w-w}), (\ref{cxha:Lie-vZ}) and the zero
curvature of $\MM^4$ we obtain
\begin{eqnarray*}
\ww\big(\biginn{\na_Z \ww, Z} \big) & =
	& \biginn{\na_{\ww} \na_Z \ww, Z}
	+\biginn{\na_Z \ww, \na_{\ww} Z} \\
	& = & \biginn{\na_Z \na_{\ww} \ww , Z} +
	\biginn{\na_Z \ww, \na_Z \ww} \,.
\end{eqnarray*}
Using (\ref{cxha:na_w-w}) and (\ref{cxha:SFR3}) we see that
this is zero; hence if (\ref{cxha:SFR3}) holds at one point of a
ray, it holds at all points of that ray.
\end{proof}

With projection
and extension defined as above, the precise statement of our construction is the following.

\begin{theorem} \label{cxha:th:SFRconf}
Let\/ $\vec{p} \in \MM^4$, and let\/ $A^3$ be an open subset of\/ $\RR^3_{\vec{p}}$. 
Then projection onto the slice\/ $\RR^3_{\vec{p}}$ defines a bijective
correspondence between germs at\/ $A^3$ of\/ $C^{\infty}$  shear-free
ray congruences\/ $\ell$ defined on an open neighbourhood of\/ $A^3$ in\/
$\MM^4$ and\/ $C^{\infty}$ conformal foliations\/ $\Cc$ by curves
on\/ $A^3$. The inverse of projection is extension.
\qed \end{theorem}

\begin{remark} \label{cxha:rem:twist}  \rm
The \emph{twist}, \emph{rotation} or \emph{vorticity} tensor of a ray congruence is
defined at a point ${\vec{x}}$ by
\begin{align*}
T(X,Y) &= \half g(\nabla_X Y - \nabla_Y X, \ww) = \half g([X,Y], \ww) \\
	&= \half \big\{g(\nabla_Y \ww, X) - g(\nabla_X \ww, Y)\big\}  \qquad\qquad  (X,Y \in \Si_{\vec{x}}).
\end{align*}
It measures
how much (infinitesimally) nearby null lines passing through the
screen space $\Si_{\vec{x}}$ twist around the null line through ${\vec{x}}$.
It is easily seen that \emph{a congruence of null lines $\ell$ is twist-free
if and only if its projection $\Cc$ onto {\rm one} $\RR^3$-slice
has integrable horizontal distribution}.
\end{remark}

\begin{example} \label{cxha:ex:radial}
Let $f(x_1,x_2,x_3) = (x_2 + \ii x_3) \big/ (x_1 \pm \abs{\vec{x}})$. It is
easily checked that $f$ is a horizontally conformal map from
$\RR^3 \setminus \{ (0,0,0) \}$ to $\Ci \cong S^2$.  Its level curves are
radii from the origin and are the leaves of a conformal
foliation $\Cc$ of $\RR^3 \setminus \{ (0,0,0) \}$ whose unit tangent
vector field is given by
\begin{equation}
\UU(x_1, x_2, x_3) = \pm \frac{1}{\sqrt{{x_1}^{\! 2} + {x_2}^{\! 2} + {x_3}^{\! 2}}}
\,(x_1, x_2, x_3) \:. \label{cxha:Urad}
\end{equation}
We remark that, when $\RR^3$ is conformally compactified to $S^3$
(see Section \ref{cxha:sec:twistor}), this example
becomes projection from the poles: $S^3 \setminus \{(\pm
1,0,0,0)\} \ra S^2$ given by the formula $(x_0,x_1,x_2,x_3)
\mapsto \pm (x_1,x_2,x_3)/\sqrt{{x_1}^{\! 2} + {x_2}^{\! 2} + {x_3}^{\! 2}}\;.$

It is easily checked that the shear-free ray congruence
$\ell$ extending $\UU$ is defined on $\MM^4 \setminus
\{(t,0,0,0):t\in \RR\}$ and has tangent vector field
$\ww = \pa/\pa t + \UU$ where
\begin{equation}
\UU(t, x_1, x_2, x_3) = \pm
\frac{1}{\sqrt{{x_1}^{\! 2} + {x_2}^{\! 2} +{x_3}^{\! 2}}}\,(x_1, x_2, x_3) \;.
\label{cxha:radialSFR} \end{equation}
For each $t$, \ $\ell$ projects to a conformal foliation on
$\RR^3_t$ with tangent vector field given by (\ref{cxha:radialSFR});
note that this conformal foliation is independent of $t$.   See also Example
\ref{cxha:ex:radial-twistor}.
\end{example}

\begin{example} \label{cxha:ex:circles}
Let $f(x_1,x_2,x_3) = -\ii x_1 \pm \sqrt{{x_2}^{\! 2} + {x_3}^{\! 2}}$. This
is a horizontally conformal submersion from $A^3 = \RR^3
\setminus \{ (x_1,x_2, x_3) : x_2 = x_3 = 0 \}$ to $\CC$. Its
level curves are circles in planes parallel to the
$(x_2,x_3)$-plane and centred on points of the $x_1$-axis; these
give a conformal foliation $\Cc$ of $A^3$ with a unit tangent vector
field given by
\begin{equation}
\UU(x_1, x_2, x_3) = \pm \frac{1}{\sqrt{{x_2}^{\! 2} + {x_3}^{\! 2}}}
\,(0, -x_3, x_2) \;.
\label{cxha:Ucircles} \end{equation}

We compute the tangent vector field to the SFR congruence $\ell$ which
extends $\UU$. The affine null geodesic of $\ell$ in $\MM^4$
through $(x_1, x_2, x_3)$ with direction $\pa / \pa t + \UU$ is
given parametrically by
\begin{eqnarray}
T & \mapsto & \left(T, x_1, x_2 
\mp x_3T\big/\sqrt{{x_2}^{\! 2} + {x_3}^{\! 2}}\, ,\ 
x_3 \pm x_2T\big/\sqrt{{x_2}^{\! 2} + {x_3}^{\! 2}}  \right)
\label{cxha:NG} \\
& = & (T,X_1,X_2,X_3) \ , \mbox{ say.} \nonumber
\end{eqnarray}
Conversely, given ${\vec{x}} = (T,X_1,X_2,X_3) \in \MM^4$ the null geodesic
of $\ell$ through ${\vec{x}}$ hits the $\RR^3$-slice $t=0$ at
$(x_1, x_2, x_3)$ where
$$
x_1 = X_1, \quad x_2 \mp x_3T\big/\sqrt{{x_2}^{\! 2} + {x_3}^{\! 2}}  =X_2,
\quad x_3 \pm x_2T\big/\sqrt{{x_2}^{\! 2} + {x_3}^{\! 2}}  =X_3 \;.
$$
This has solution
\begin{equation}
(x_1, x_2, x_3) = \left( X_1,\frac{R}{{X_2}^{\! 2} + {X_3}^{\! 2}}(RX_2\pm TX_3), 
\frac{R}{{X_2}^{\! 2}+{X_3}^{\! 2}}(RX_3\mp TX_2) \right)
\label{cxha:circlesSFRsoln} \end{equation}
where $R = \sqrt{{X_2}^{\! 2}+{X_3}^{\! 2}-T^2}$\,.
Hence the tangent vector field to the null geodesic of $\ell$
through $(t, x_1, x_2, x_3) \in \MM^4$ is given by $\ww = \pa / \pa t + \UU$ with
\begin{equation}
\UU = \UU_t(x_1,x_2,x_3) = \UU(t,x_1,x_2,x_3) =
\frac{r}{\sqrt{{x_2}^{\! 2}+{x_3}^{\! 2}}}
	\left(0\,,\,-x_3 \pm \frac{t}{r}x_2,\,\, x_2 \pm\frac{t}{r}x_3 \right)
\label{cxha:E2} \end{equation}
where $r = \sqrt{{x_2}^{\! 2}+{x_3}^{\! 2}-t^2}$\,.

For each $t$, the integral curves of $\UU_t$ give a conformal (in fact,
Riemannian) foliation $\Cc_t$ of the $\RR^3$-slice $t = \const$;
it is easy to see that the leaves of $\Cc_t$ are the
involutes of circles (see Figure \ref{cxha:fig:evolutes}).
See Examples \ref{cxha:ex:circles-hyp} and \ref{cxha:ex:circles-twistor}
for further developments.
\end{example}

\section{Complex-analytic foliations by null planes on\/ $\CC^4$}
\label{cxha:sec:nullplanes}
\begin{definition} {\rm (LeBrun 1983)} \label{cxha:def:holmetric}
A\/ {\rm holomorphic metric} on a complex manifold\/ $M$ is a
holomorphic section\/ $g^{\cC}$ of\/ $\odot^2 T^*_{1,0}M$ which defines a
(complex-sym\-metric bilinear) non-degenerate inner product on\/ $T^{1,0}_xM$ for
each\/ $x \in M$.
\end{definition} 

Let $\CC^4 = \{(x_0,x_1,x_2,x_3) : x_i \in \CC^4 \}$ equipped with the
standard holomorphic metric $g^{\cC} = \dd {x_0}^{\! 2} + \dd {x_1}^{\! 2} + \dd {x_2}^{\! 2} + \dd {x_3}^{\! 2}$.
A vector $w$ is called \emph{(complex-)null} or \emph{isotropic} if
$g^{\cC}(w,w) = 0$. Let ${\vec{p}} \in \CC^4$.  A subspace in $T_{\vec{p}}\CC^4$ is called
\emph{null} if it consists of null vectors.
For any $\vec{p} = (p_0,p_1,p_2,p_3) \in \CC^4$, by the $\RR^4$-\emph{slice through
${\vec{p}}$} we mean the real four-dimensional
affine subspace 
$$
\RR^4_{\vec{p}} = \bigl\{ (x_0,x_1,x_2,x_3) \in \CC^4 : \Im x_i = \Im
p_i \ (i=0,1,2,3) \bigr\}
$$ 
parametrized by
\begin{equation} \label{cxha:R4p}
\RR^4 \ni (x_0,x_1,x_2,x_3) \mapsto (p_0 + x_0, p_1 + x_1, p_2 + x_2, p_3 +
x_3) \in \RR^4_{\vec{p}} \subset \CC^4
\end{equation}
and with the orientation induced from the standard orientation of $\RR^4$.
If $\vec{p}$ is the zero vector ${\vec{0}} = (0,0,0,0)$ we write $\RR^4$ for $\RR^4_{\vec{0}}$.   Note
that $T_{\vec{p}} \CC^4$ can be identified canonically with $T_{\vec{p}} \RR^4_{\vec{p}} \otimes
\CC$; we shall frequently make this identification.

In a similar way, extending the notation of Section \ref{cxha:sec:SFR}, we define the \emph{$\RR^3$-slice
through $\vec{p}$} by
$$
\RR^3_{\vec{p}} = \bigl\{ (x_0,x_1,x_2,x_3) \in \CC^4 : x_0 = p_0,\, \Im x_i = \Im
p_i \ (i=1,2,3) \bigr\}
$$ 
parametrized by
\begin{equation} \label{cxha:R3p}
\RR^3 \ni (x_1,x_2,x_3) \mapsto (p_0, p_1 + x_1, p_2 + x_2, p_3 +
x_3) \in \RR^3_{\vec{p}} \subset \CC^4
\end{equation}
and with the orientation induced from the standard orientation of $\RR^3$;
we write $\RR^3$ for $\RR^3_{\vec{0}}$.

Given any orthonormal basis $\{ e_0,e_1,e_2,e_3 \}$ of the tangent space $T_{\vec{p}}\RR^4_{\vec{p}}$,
the plane $\Pi_{\vec{p}} = \mbox{span}\{e_0+\ii e_1, e_2 + \ii
e_3\}$ is null; we call $\Pi_{\vec{p}}$ an \emph{$\alpha$-plane} (respectively,
\emph{$\beta$-plane}) according as the basis $\{ e_0,e_1,e_2,e_3\}$
is positively (respectively, negatively) oriented. 

By a \emph{complex-analytic distribution of $\alpha$-planes} on an
open set $A^C$ of $\CC^4$ we mean a map $\Pi$ which assigns to
each point ${\vec{p}}$ of $A^C$ an $\alpha$-plane $\Pi_{\vec{p}} \subset
T_{\vec{p}}\CC^4$ in a complex-analytic fashion, i.e., $\Pi_{\vec{p}} = \mbox{span}
\{w_1({\vec{p}}), w_2({\vec{p}})\}$ where the maps $w_i:A^C \ra T\CC^4$ are complex analytic.
 We can identify $\Pi$ with its image, a complex-analytic subbundle of
$T\CC^4\vert_{A^C}$.  

By Frobenius' Theorem, a complex-analytic distribution $\Pi$ of
$\alpha$-planes on $A^C$ is integrable if and only if
\begin{equation} \label{cxha:alphaint}
[w_1,w_2] \in \Ga(\Pi) \qquad \big( w_1, w_2 \in \Ga(\Pi) \big) \,.
\end{equation}

\begin{lemma} \label{cxha:lem:intauto}
A complex-analytic distribution\/ $\Pi$ of\/ $\alpha$-planes on\/
$A^C$ is integrable if and only if it is\/ {\rm autoparallel}, i.e.,
\begin{equation} \label{cxha:autopar}
\na_{w_1}w_2 \in \Ga(\Pi) \qquad \big( w_1, w_2 \in
\Ga(\Pi) \big) \,.
\end{equation}
\end{lemma}

\begin{proof}
Condition (\ref{cxha:alphaint}) holds if and only if
$$
g^{\cC} \big([w_1,w_2],w_1 \big) = 0 \qquad \big( w_1, w_2 \in
\Ga(\Pi) \big) \,.
$$
Since $g^{\cC}(\na_{w_2}w_1, w_1) = \half w_2 \big( g^{\cC}(w_1, w_1) \big) = 0$,
this holds if and only if
\begin{equation} \label{cxha:autopar2}
g^{\cC}(\na_{w_1}w_2, w_1) = 0 \qquad \big( w_1, w_2 \in
\Ga(\Pi) \big) \,.
\end{equation}
Now, $g^{\cC}(\na_{w_1}w_2, w_2) = 0$, hence (\ref{cxha:autopar2})
holds if and only if (\ref{cxha:autopar}) holds.
\end{proof}

Thus an integrable complex-analytic distribution
of $\alpha$-planes has integral submanifolds which are (affine)
$\alpha$-planes, i.e., planes of $\CC^4$ whose tangent spaces are
$\alpha$-planes; thus, \emph{an integrable complex-analytic
distribution of $\alpha$-planes on $A^C$ is equivalent to a
complex-analytic foliation by affine $\alpha$-planes of $A^C$.} 
Similar considerations apply if we replace `$\alpha$-plane' by
`$\beta$-plane'.

\section{Unification} \label{cxha:sec:unification}
The formulae (\ref{cxha:conffoln1}), (\ref{cxha:Jint1}),
(\ref{cxha:SFR1}) and (\ref{cxha:alphaint}) show that the four `distributions' (or `fields')
studied in Sections \ref{cxha:sec:conffoln}--\ref{cxha:sec:nullplanes}
are invariant if the standard metrics are replaced by conformally
equivalent metrics.  This suggests that there may be some relations
between the distributions; we now explain those relations.

For any ${\vec{p}} = (p_0,p_1,p_2,p_3) \in \CC^4$, by analogy with the $\RR^4$-slice through ${\vec{p}}$
defined in Section \ref{cxha:sec:nullplanes}, we define the \emph{Minkowski slice through ${\vec{p}}$} 
to be the real four-dimensional
affine subspace 
$$
\MM^4_{\vec{p}} = \bigl\{ (x_0, x_1, x_2, x_3 ) : \Re x_0 = \Re p_0, \ \Im x_i
= \Im p_i \ (i=1,2,3) \bigr\}
$$ 
parametrized by 
$$
\MM^4 \ni
(t,x_0,x_2,x_2) \mapsto (p_0 - \ii t, p_1 + x_1, p_2 + x_2, p_3 +
x_3) \in \MM^4_{\vec{p}} \subset \CC^4.
$$
If ${\vec{p}}= {\vec{0}}$, we write
$\MM^4 \equiv \MM^4_{\vec{0}}$.  (Note that the minus sign before the term $\ii t$
is unimportant---it is simply to avoid minus signs elsewhere.)  

Now we use the slice notation to show how all the four quantities
discussed above are equivalent at a point.

\begin{proposition} \label{cxha:prop:unifpt}
Let\/ ${\vec{p}} \in \CC^4$.
There are one-to-one correspondences between the quantities\/{\rm :}

{\rm (QP1)} $\alpha$-planes\/ $\Pi_{\vec{p}} \subset T_{\vec{p}}\CC^4$\,,

{\rm (QP2)} positive almost Hermitian structures\/ $\JJ_{\vec{p}}:T_{\vec{p}}\RR^4_{\vec{p}} \to 
T_{\vec{p}}\RR^4_{\vec{p}}$\,,

{\rm (QP3)} unit vectors\/ $\UU_{\vec{p}} \in T_{\vec{p}}\RR^3_{\vec{p}}$\,,

{\rm (QP4)} null directions\/ $\WW_{\vec{p}} \subset T_{\vec{p}}\MM^4_{\vec{p}}$\,,

\noindent given by

\begin{equation} \left. \begin{array}{rclcl}
\JJ_{\vec{p}}  & = & \JJ_{\vec{p}}(\Pi_{\vec{p}}) & =  &-\ii \mbox{ on } \Pi_{\vec{p}}\,, \quad
 +\ii \mbox{ on } \ov{\Pi}_{\vec{p}}\,;
\\
\Pi_{\vec{p}} & = & \Pi_{\vec{p}}(\JJ_{\vec{p}}) & = & (0,1)\mbox{-tangent space of\/ } \JJ_{\vec{p}},
\text{ i.e., } \Pi_{\vec{p}} = \{ X+\ii \JJ_{\vec{p}} X : X \in T_{\vec{p}} \RR^4_{\vec{p}} \}\,; \\             
\UU_{\vec{p}}  & = & \UU_{\vec{p}}(\JJ_{\vec{p}}) & = & \JJ_{\vec{p}}(\pa / \pa x_0) \,; \\
\JJ_{\vec{p}} & = & \JJ_{\vec{p}}(\UU_{\vec{p}}) & = & \mbox{the unique positive almost Hermitian
structure at\/ }{\vec{p}} \\
& & & &\mbox{ with\/ } \JJ_{\vec{p}}(\pa / \pa x_0) = \UU_{\vec{p}}\,;  \\ 
\WW_{\vec{p}}  & = & \WW_{\vec{p}}(\UU_{\vec{p}}) & = & \mbox{\rm span}\{\ww_p\} \mbox{ where }
	\ww_p = \pa / \pa t + \UU_{\vec{p}}\,; \\
\UU_{\vec{p}} & = & \UU_{\vec{p}}(\WW_{\vec{p}}) & = & \mbox{normalized projection of\/ } \WW_{\vec{p}}
\mbox{ onto } T_{\vec{p}}\RR^3_{\vec{p}}\,, \mbox{i.e. the unique } \\
& & & & \ \UU_{\vec{p}} \in T_{\vec{p}}\RR^3 \mbox{such that\/ }\WW_{\vec{p}} =
\mbox{\rm span}(\pa / \pa t + \UU_{\vec{p}})\,; \\
\WW_{\vec{p}} & = & \WW_{\vec{p}}(\Pi_{\vec{p}}) & = & \Pi_{\vec{p}} \cap T_{\vec{p}}\MM^4_{\vec{p}}\,; \\
\Pi_{\vec{p}} & = & \Pi_{\vec{p}}(\WW_{\vec{p}}) & = & \mbox{the unique $\alpha$-plane containing
$\WW_{\vec{p}}$} \;.   
\end{array} \right\} \label{cxha:unif} \end{equation}
\qed \end{proposition}

Note that we shall use an upright font:  $\JJ$, $\UU$, $\ww$, etc.\ to denote quantities which are related
to each other as in this proposition.

Next we discuss correspondences between the various \emph{distributions}.

Let ${\vec{p}} \in \CC^4$.  Suppose that $J$ is a real-analytic almost Hermitian
structure on a neighbourhood $A^4$ of ${\vec{p}}$ in $\RR^4_{\vec{p}}$.  Then we can
extend $J$ to a neighbourhood $A^C$ of ${\vec{p}}$ in $\CC^4$ by asking that it be
complex analytic.  Note that this complex analyticity may be expressed by
\begin{equation} \label{cxha:cxan}
\na_{\ii X}J = J\na_X J \qquad (X \in T_{\vec{q}}\CC^4, \ {\vec{q}} \in A^C) \,.
\end{equation}
Via Proposition \ref{cxha:prop:unifpt}, this defines a
complex-analytic distribution of $\alpha$-planes.   
 
Conversely, a complex-analytic distribution of $\alpha$-planes on a
neighbourhood $A^C$ of ${\vec{p}}$ defines a
real-analytic almost Hermitian structure on $A^4 = A^C \cap \RR^4_{\vec{p}}$\,.

Similarly, a real-analytic distribution $\WW$ of null directions on
a neighbourhood $A^M$ of ${\vec{p}}$ in $\MM^4_{\vec{p}}$
gives rise to a complex-analytic distribution of $\alpha$-planes on a
neighbourhood $A^C$ of ${\vec{p}}$ in $\CC^4$, and conversely.  Then we have the following.

\begin{proposition} \label{cxha:prop:unifdist}
Let\/ $\Pi$ be a complex-analytic distribution of\/ $\alpha$-planes on an open
subset\/ $A^C$ of\/ $\CC^4$, let\/ ${\vec{p}} \in A^C$, and let\/ $\JJ$ {\rm (}respectively,
$\WW${\rm )} be the corresponding almost Hermitian structure on\/ $A^4 = A^C \cap
\RR^4_{\vec{p}}$ {\rm (}respectively, distribution of null directions on
$A^M = A^C \cap \MM^4_{\vec{p}}${\rm )}.
Then the following conditions are equivalent\/{\rm :}
\begin{itemize}
	\item[{\rm (i)}] $\Pi$ is integrable on\/ $A^C$, 
	\item[{\rm (ii)}] $\JJ$ is integrable on\/ $A^4$,
	\item[{\rm (iii)}] $\WW$ defines a shear-free ray congruence on\/ $A^M$.
\end{itemize}
\end{proposition}

\begin{proof}
Firstly we show that \emph{$\Pi$ is integrable at points of\/ $A^4$
if and only if\/ $\JJ$ is integrable on\/ $A^4$}.  To do this, let ${\vec{q}} \in A^4$.
In view of Lemma \ref{cxha:lem:intauto}, the distribution $\Pi$ is
integrable at ${\vec{q}}$ if and only if
\begin{equation} \label{cxha:auto}
\na_{X+\ii \JJ X}\JJ = 0 \qquad (X \in T_{\vec{q}}\RR^4_{\vec{q}}) \,.
\end{equation}
By (\ref{cxha:cxan}), this is equivalent to
$$
\na_{\JJ X}\JJ = \JJ\na_X \JJ \qquad (X \in T_{\vec{q}}\RR^4_{\vec{q}})
$$
which is the integrability condition (\ref{cxha:Jint2}) for $\JJ$ at ${\vec{q}}$.

Next we show that \emph{$\Pi$ is integrable at points of\/ $A^M$ if
and only if\/ $\WW$ is a shear-free ray congruence on\/ $A^M$}.  To do
this, take ${\vec{q}} \in A^M$.  Then with $\UU = \JJ(\pa/\pa x_0)$, since
$\pa/\pa x_0$ is parallel, (\ref{cxha:auto}) is equivalent to
\begin{equation} \label{cxha:autoU}
\na_{\JJ X}\UU = \JJ\na_X \UU \qquad (X \in T_{\vec{q}}\RR^4_{\vec{q}})\,,
\end{equation}
and so, with $\ww = \pa/\pa t + \UU$,
\begin{equation}
\na_{X+\ii\JJ X}\ww = 0 \qquad (X \in T_{\vec{q}}\RR^4_{\vec{q}}) \;.
\label{w} \end{equation}
Choose $X = \pa/\pa x_0 = \ii\pa/\pa t$; then this reads
\begin{equation} \label{cxha:geod}
\na_{\ww} \ww = 0\,;
\end{equation}
this is the condition that the integral curves of $\ww$ be geodesic at ${\vec{q}}$.
With
$X$ chosen instead in the screen space $\WW_{\vec{q}}^{\perp} \cap \RR^3_{\vec{q}} = \UU_{\vec{q}}^{\perp}
\cap \RR^3_{\vec{q}}$ we obtain the shear-free condition (\ref{cxha:SFR2}).
Conversely, since the two choices of $X$ give vectors $X+\ii \JJ X$
spanning $\Pi_{\vec{q}}\,$, if (\ref{cxha:geod})  and (\ref{cxha:SFR2}) hold
so does (\ref{cxha:auto}).

To finish the proof, if any of the above conditions
holds at all points of $A^4$ or $A^M$, by analytic continuation
it holds throughout $A^C$.
\end{proof}

Combining this result with Theorem \ref{cxha:th:SFRconf} we obtain
relations between the following four sorts of distributions
(for any ${\vec{p}} \in \CC^4$):

\begin{equation} \label{cxha:Q}
\left. \begin{array}{rl}
\mbox{(Q1)} & \mbox{holomorphic foliations } \Ff \mbox{ by }
\alpha\mbox{-planes } \Pi \mbox{ of an
open subset } A^C \mbox{ of } \CC^4 \,,
\\
\mbox{(Q2)} & \mbox{positive Hermitian structures } \JJ \mbox{ on an
open subset } A^4 \mbox{ of } \RR^4_{\vec{p}} \,,
\\
\mbox{(Q3)} & \mbox{real-analytic shear-free ray congruences } \ell
     \mbox{ on an open subset } A^M \mbox{ of } \MM^4_{\vec{p}}\,,
\\
\mbox{(Q4)} & \mbox{real-analytic conformal foliations } \Cc
\mbox{ by curves of an open subset } A^3 \mbox{ of } \RR^3_{\vec{p}}\,. 
\end{array} \space \space \space \right\}
\end{equation}

\begin{theorem} \label{cxha:th:unifgerms}
Let\/ $A^C$ be an open subset of\/ $\CC^4$, and let $\vec{p} \in A^C$. Suppose that we are given
\par {\rm (i)} a complex-analytic foliation by null planes on\/ $A^C$.

\noindent Then restriction to slices through\/ ${\vec{p}}$ defines
\par {\rm (ii)} a positive Hermitian structure on\/ $A^C \cap \RR^4_{\vec{p}}$\,,
\par {\rm (iii)} a real-analytic shear-free ray congruence on\/ $A^C \cap
\MM^4_{\vec{p}}$\,,
\par {\rm (iv)} a real-analytic conformal foliation by curves on
$A^C \cap \RR^3_{\vec{p}}$\,.

Further, for a fixed open set\/ $A^3$ of\/ $\RR^3_{\vec{p}}$\,, these maps
define bijections between germs at\/ $A^3$ of the distributions {\rm (i)},
{\rm(ii)}, {\rm(iii)}, {\rm(iv)}.
\qed \end{theorem}

\begin{corollary} \label{cxha:cor:JtoU}
Let\/ ${\vec{p}} \in \RR^4$, and let\/ $A^3$ be an open subset of\/ $\RR^3_{\vec{p}}$. 
Then projection \/ $\JJ \mapsto \UU = \JJ(\pa/ \pa x_0)$ onto the slice\/
$\RR^3_{\vec{p}}$ defines a bijective
correspondence between germs at\/ $A^3$ of positive Hermitian structures\/ $\JJ$
defined on an open neighbourhood of\/ $A^3$ in\/ $\RR^4$ and
real-analytic conformal foliations\/ $\Cc$ of\/ $A^3$ by curves.
\qed \end{corollary}

\begin{remark} \label{cxha:rem:projn}  \rm
Given a \emph{real-analytic} conformal foliation $\Ff$ by curves
of an open subset $A^3$ of $\RR^3$, there corresponds a
$5$-parameter family of `associated' real-analytic conformal foliations by
curves of open subsets of $\RR^3$ obtained by extending $\Ff$ to a
holomorphic foliation by $\alpha$-planes and then projecting this
to all possible $\RR^3$-slices.
\end{remark}

\section{Local coordinates and harmonic morphisms}
\label{cxha:sec:loccoords}
A smooth map between Riemannian manifolds
is a harmonic morphism if and only if it is both (i) harmonic
and (ii) horizontally weakly conformal (Fuglede 1978, Ishihara 1979); this remains true for
semi-Riemannian manifolds (Fuglede 1996).  For  a complex-valued map
$\phi:M \to \CC$ from a (semi-)\linebreak
Riemannian manifold these conditions read, respectively,
\begin{equation} \label{cxha:hamo}
{\rm (i)} \ \Delta^M\phi = 0\,, \quad \text{and} \quad {\rm (ii)} \  \inn{\grad\phi\,,\, \grad\phi} = 0\,,
\end{equation}
where $\Delta^M$ is the Laplace-Beltrami operator on $M$, see Section {\bf \ref{pre:sec:ha-fn}}
for more details.  Note that a complex-valued harmonic morphism from a real-analytic Riemannian manifolds
is real analytic (see Proposition {\bf \ref{fund:prop:cts-hamorph}}); this is false, in general, for a
semi-Riemannian domain (see Section {\bf\ref{semi:sec:hamorph}}).  The
pair of equations \eqref{cxha:hamo} is \emph{conformally invariant} in the sense that
the composition $f \circ \phi$ of a harmonic morphism $\phi:M \to \CC$ and
a weakly conformal (i.e., holomorphic or antiholomorphic) map $f$ to $\CC$
is a harmonic morphism. Hence, we can replace $\CC$ by a \emph{Riemann surface} $N$
and then \eqref{cxha:hamo} are the equations for a harmonic morphism $\phi:M \to N$
in any complex (or isothermal) coordinates (cf.\ Sections
{\bf\ref{Euclid3:sec:defchar}} and {\bf\ref{fund:sec:def}}). 

It is convenient to introduce \emph{standard null coordinates} 
$(q_1, \wt{q}_1, q_2, \wt{q}_2)$ on $\CC^4$ by setting 
\begin{equation}
q_1 = x_0 + \ii x_1, \ \wt{q}_1 = x_0 - \ii x_1, \ q_2 = x_2 +
\ii x_3, \ \wt{q}_2 = x_2 - \ii x_3.
\label{cxha:nullcoords} \end{equation}
Then $\RR^4$ is given by $\wt{q}_1 = \ov{q}_1, \wt{q}_2 =
\ov{q}_2$ and the holomorphic metric on $\CC^4$ by $g^{\cC} = \dd q_1
\dd\wt{q}_1 + \dd q_2 \dd\wt{q}_2$; this restricts to the
Euclidean metric $g = \dd q_1 \dd \ov{q}_1 + \dd q_2 \dd \ov{q}_2$ on real slices $\RR^4_{\vec{p}}$.
(Here we write $\ov{q}_i$ for the complex conjugate $\ov{q_i}$ of $q_i$.)

For any ${\vec{p}} \in \CC^4$ and $[w_0, w_1] \in \CP^1$, set
\begin{equation}
\Pi_{\vec{p}} = \mbox{span} \left\{w_0 \frac{\pa}{\pa \wt{q}_1}
- w_1 \frac{\pa}{\pa q_2} \,, \ 
w_0 \frac{\pa}{\pa \wt{q}_2} + w_1 \frac{\pa}{\pa q_1} \right\}
\,. \label{cxha:NT0} \end{equation}
It is easily checked that $\Pi_{\vec{p}}$ is an $\alpha$-plane.
Write $\mu = w_1/w_0 \in \Ci$; then if $\mu \neq \infty$, \  $\Pi_{\vec{p}}$
has a basis
\begin{equation}
\left\{ \frac{\pa}{\pa \wt{q}_1} - \mu \frac{\pa}{\pa q_2} \,,
\ \frac{\pa}{\pa \wt{q}_2} + \mu \frac{\pa}{\pa q_1} \right\}
\,; \label{cxha:NT}
\end{equation}
if $\mu = \infty$, a basis for $\Pi_{\vec{p}}$ is given by
$\left\{ \pa/\pa q_2 \,, \ \pa/\pa q_1 \right\}$.

The assignment (\ref{cxha:NT0}) gives an explicit bijection
$$
\CP^1 \llra \{ \alpha\mbox{-planes at } {\vec{p}} \}.
$$
On applying Proposition \ref{cxha:prop:unifpt}(ii), we obtain also an almost
Hermitian structure $\JJ_{\vec{p}} = \JJ_{\vec{p}}(\Pi_{\vec{p}})$ on $T_{\vec{p}}\RR^4_{\vec{p}}$.  This has
$(0,1)$-tangent space with basis
\begin{equation}
\left\{ w_0 \frac{\pa}{\pa \ov{q}_1} - w_1 \frac{\pa}{\pa q_2} \,, \
w_0 \frac{\pa}{\pa \ov{q}_2} + w_1 \frac{\pa}{\pa q_1}  \right\} \,,
\label{cxha:cotanghom} \end{equation} 
or, if $\mu = w_1/w_0 \neq \infty$,
\begin{equation}
\left\{ \frac{\pa}{\pa \ov{q}_1} - \mu \frac{\pa}{\pa q_2} \,, \
\frac{\pa}{\pa \ov{q}_2} + \mu \frac{\pa}{\pa q_1} \right\} \,,
\label{cxha:tang} \end{equation}
(cf.\ Section {\bf \ref{twistor:sec:twistorR4}}).

On applying Proposition \ref{cxha:prop:unifpt}(iii) and (iv), we also obtain a
unit vector $\UU_{\vec{p}}$ and a null direction $\WW_{\vec{p}} = \mbox{span}\{\ww_{\vec{p}}\}$
where $\ww_{\vec{p}} = \pa/\pa t + \UU_{\vec{p}}$.  We
shall say that $\mu$ \emph{represents} $\Pi_{\vec{p}},\ \JJ_{\vec{p}}, \UU_{\vec{p}}$ and $\WW_{\vec{p}}$.
To get more explicit formulae for these quantities, we apply the unitary matrix
$\dfrac{1}{\abs{w_0}^2 + \abs{w_1}^2}\left(
\begin{array}{rr} \ov{w}_0 & \ov{w}_1 \\
                           -w_1  &           w_0  
\end{array} \right)$;
this converts the basis of $\Pi_{\vec{p}}$ to
another basis $\{e_0+\ii e_1, e_2 + \ii e_3 \}$ where $\{ e_0, e_1, e_2, e_3 \}$
is a positive orthonormal basis for $\RR^4$. Explicitly, write $u
= \ii w_1/w_0 = \ii\mu$, and let $\si:S^2 \ra \Ci$ denote
stereographic projection from $(-1,0,0)$, then the
basis $\{ e_0, e_1, e_2, e_3 \}$ is given by

\begin{equation}
\left. \begin{array}{rcl}
e_0 & = & \pa / \pa x_0 \\  \vspace{1ex}
\UU_{\vec{p}} = e_1 &
	= & \displaystyle{\frac{1}{1+\abs{u}^2}} \big(1-\abs{u}^2,\, 2\Re u,\,2\Im u \big) = \si^{-1}(u) \ \\
\vspace{1ex}
e_2 + \ii e_3 &
	= & \displaystyle{\frac{1}{1+\abs{u}^2}} \big(-2u,\,1-u^2,\, \ii (1+u^2) \big) \;.
\end{array} \right\}
\label{cxha:stdbasis} \end{equation}
Thus $\JJ_{\vec{p}}$ is the unique positive almost Hermitian structure at ${\vec{p}}$ which satisfies 
$\JJ_{\vec{p}}(\pa/\pa x_0) = \si^{-1}(\ii\mu)$.

Now let there be given a smooth distribution of any of the quantities $\Pi$, $\JJ$, $\WW$, $\UU$
of (\ref{cxha:Q}) related as in Proposition \ref{cxha:prop:unifpt}.
Write $\UU = \si^{-1}(\ii \mu)$;  then, after applying a Euclidean transformation, if
necessary, to ensure that $\mu$ is finite, we see that (\ref{cxha:auto}) is
equivalent to the equation
\begin{equation}  
Z (\mu )= 0 \qquad (Z \in \Pi_{\vec{q}}) \;.
\label{cxha:nablamu} \end{equation}
Recall that, in standard null coordinates $(q_1, \wt{q}_1, q_2,
\wt{q}_2)$, a basis for $\Pi$
is given by (\ref{cxha:NT0}) or (\ref{cxha:NT}).  Thus
(\ref{cxha:nablamu}) reads
\begin{equation}
\left( \frac{\pa}{\pa \wt{q}_1} - \mu \frac{\pa}{\pa q_2}
\right)\! \mu = 0 \,,
\quad
\left( \frac{\pa}{\pa \wt{q}_2} + \mu \frac{\pa}{\pa q_1}
\right)\! \mu = 0 \,;
\label{cxha:EC} \end{equation}
this is the condition that the distribution $\Pi$ of $\alpha$-planes
represented by $\mu$ be
integrable.  These equations restrict on real slices to the equations
\begin{equation}
\left( \frac{\pa}{\pa \ov{q}_1} - \mu \frac{\pa}{\pa q_2} \right)\!\mu = 0 \,,
\quad
\left( \frac{\pa}{\pa \ov{q}_2} + \mu  \frac{\pa}{\pa q_1}
\right)\!\mu = 0 \,;
\label{cxha:ER} \end{equation}
this is the condition that the almost Hermitian structure $\JJ$
represented by $\mu$ be integrable. 
On Minkowski slices, writing $v = x_1 + t$, $w = x_1 - t$, these equations restrict
to
\begin{equation}
\left( \frac{\pa}{\pa v} + \ii\mu \frac{\pa}{\pa q_2} \right)\!\mu = 0
\ , \quad
\left( \frac{\pa}{\pa \ov{q}_2} - \ii\mu \frac{\pa}{\pa w}
\right)\! \mu = 0 \;;
\label{cxha:EM} \end{equation}
this is the condition that the integral curves of the smooth
distribution $\WW$ of null directions represented by $\mu$
form a shear-free ray congruence.

\begin{example} \label{cxha:ex:Hopf}
Define $\mu: \CC^4 \setminus \{\vec{0}\} \to \Ci$ by
\begin{equation} \label{cxha:Hopfcx}
\mu = - q_2/\wt{q}_1 \,.
\end{equation}
This satisfies equations (\ref{cxha:EC}), and so
defines a holomorphic foliation $\Pi$ by $\al$-planes.  In fact, the
leaves of this foliation form the holomorphic $2$-parameter family
of $\al$-planes:
$$
q_1 - \mu \wt{q}_2 = \nu \,, \ q_2 + \mu\wt{q}_1 = 0 \qquad (\mu \in \Ci, \ \nu
\in \CC) \,.
$$

The map $\mu$ restricts to the map $\RR^4 \setminus \{\vec{0}\} \to \Ci$ given by
$\mu = q_2/\ov{q}_1$; this satisfies (\ref{cxha:ER}) and so defines a
positive Hermitian structure $\JJ$.

On $\MM^4$ it restricts to
\begin{equation} \label{cxha:HopfMin}
\mu = -\ii q_2/(x_1+t)\,;
\end{equation}
this defines a shear-free ray congruence $\ell$.  A short calculation
shows that the vector $\UU = \si^{-1}(\ii\mu)$ is given by (\ref{cxha:bunch}),
so that the projection of $\ell$ onto
a slice $t=c$ is the foliation $\Cc_c$
described in Example \ref{cxha:ex:bunch}. 
\end{example}

We now show how the quantities $\JJ$ and $\ell$ define harmonic morphisms.

\begin{proposition} \label{cxha:prop:R4muhm}
Let\/ $\mu:A^4 \to \CC$ be a smooth function from an open subset of\/
$\RR^4$ which satisfies\/ {\rm(\ref{cxha:ER})}, i.e., 
represents a Hermitian structure.  Then\/
$\mu$ is a harmonic morphism.
\end{proposition}

\begin{proof}
Simply note that, for any smooth function $\mu:A^4 \to \CC$,
\begin{alignat*}{3}
\tfrac{1}{4}\Delta\mu &= \frac{\pa^2\mu}{\pa q_1\pa \ov{q}_1} + \frac{\pa^2\mu}{\pa q_2\pa \ov{q}_2}
	&&=\frac{\pa}{\pa q_1}\!\left( \frac{\pa\mu}{\pa \ov{q}_1}
		- \mu \frac{\pa\mu}{\pa q_2} \right) + \frac{\pa}{\pa q_2}\!
	\left( \frac{\pa\mu}{\pa \ov{q}_2} + \mu \frac{\pa\mu}{\pa q_1} \right) \,,\\[0.5ex]
\tfrac{1}{4}\inn{\grad\mu\,,\, \grad\mu}
	&= \frac{\pa\mu}{\pa q_1}\frac{\pa\mu}{\pa \ov{q}_1} + \frac{\pa\mu}{\pa q_2}\frac{\pa\mu}{\pa \ov{q}_2}
	&&=\frac{\pa\mu}{\pa q_1} \!\left( \frac{\pa\mu}{\pa \ov{q}_1}
	- \mu \frac{\pa\mu}{\pa q_2} \right) + \frac{\pa\mu}{\pa q_2}\!
	\left( \frac{\pa\mu}{\pa \ov{q}_2} + \mu \frac{\pa\mu}{\pa q_1} \right) \,.
\end{alignat*}
By \eqref{cxha:ER}, the right-hand sides vanish, 
so that $\mu$ is a harmonic morphism by the `if' part of the Fuglede--Ishihara characterization
\eqref{cxha:hamo} (see Lemma {\bf\ref{fund:lem:char1}}).

Since $\JJ = \si^{-1}(\ii\mu)$, by conformal invariance (see above or Section {\bf 4.1}),
it follows that $\JJ:A^4 \to S^2$ is also a harmonic morphism.
\end{proof}

\begin{remark}  \rm 
That $\JJ$ is a harmonic morphism also follows from Proposition
{\bf\ref{twistor:prop:superint}}, as in Example {\bf\ref{twistor:ex:R4-converse}}.
\end{remark}

We have an analogue for Minkowski space, as follows.

\begin{proposition} \label{cxha:prop:M4muhm}
Let\/ $\mu:A^M \to \CC$ be a smooth function from an open subset of
Minkowski space which 
represents a shear-free ray congruence, i.e., satisfies {\rm(\ref{cxha:EM})}. Then
$\mu$ is a harmonic morphism.
\qed \end{proposition}

\begin{example} \label{cxha:ex:radial-hamo}
The direction field $\UU$ of the SFR conguence of Example \ref{cxha:ex:radial}
defines a harmonic morphism
$\MM^4 \setminus \{(t,0,0,0):t\in \RR\} \to S^2$.  Note that this
harmonic morphism is submersive and surjective.
\end{example}

A horizontally weakly conformal map
$\phi:M \to N$ between semi-Riemann\-ian manifolds
is called \emph{degenerate at $x \in M$} if \,$\ker\dd\phi_x$ is a
degenerate subspace (Definition {\bf \ref{semi:def:HWCdeg}}).  

\begin{example} \label{cxha:ex:circles-hamo}
The direction field $\UU$ of the SFR conguence of Example \ref{cxha:ex:circles}
defines a harmonic morphism from the cone given by
$$
A^M = \{ (t,x_1,x_2,x_3) \in \MM^4: {x_2}^{\! 2}+{x_3}^{\! 2} > t^2 \}
$$
to $S^2$ which is degenerate everywhere and has image given by the
equator of $S^2$. Note that $\UU = \si (\ii\mu)$ with
$$
\mu(t,x_1,x_2,x_3) = \frac{r}{{x_2}^{\! 2}+{x_3}^{\! 2}} \left\{ \left( x_2 +
\frac{t}{r}x_3 \right) + \ii \left(x_3 - \frac{t}{r}x_2 \right)
\right\} \,.
$$
The map $\mu : A^M\ra \CC$ defines a  harmonic morphism
$\mu:A^M \ra \CC$ equivalent to $\UU$; it is also degenerate, with image the unit circle.
The fibre of $\UU$ (or $\mu$) through any point ${\vec{p}}$ is the affine
plane perpendicular to $\UU_{\vec{p}}$\,; this is spanned by $\UU_{\vec{p}}$\,, $\pa/\pa x_1$
and the vector in the $(x_2,x_3)$-plane perpendicular to $\UU_{\vec{p}}$\,.
\end{example}

There is a \emph{complex} version of the last two propositions which we explain in the next section.

\section{Complex-harmonic morphisms} \label{cxha:sec:cxha}

By a \emph{complex-harmonic function} $\phi:A^C \ra \CC$ on an
open subset $A^C$ of $\CC^4$ we mean a \emph{complex-analytic map}
satisfying the complexified Laplace's equation:
\begin{equation}
\sum_{i=0}^3 \frac{\pa^2\phi}{\pa {x_i}^{\! 2}}  =  0 \qquad 
\big( (x_0,x_1,x_2,x_3) \in A^C \big) \,.   \label{cxha:cxha}
\end{equation} 
By a \emph{complex-harmonic morphism} $\phi:A^C \ra \CC$ we mean a
complex-analytic map satisfying (\ref{cxha:cxha}) and a complexified
version of the horizontal weak conformality condition:
\begin{equation}
\sum_{i=0}^3 \left( \frac{\pa \phi}{\pa x_i} \right)^{\!\! 2}  =  0 \qquad
\left( (x_0,x_1,x_2,x_3) \in A^C \right).
\label{cxha:cxHWC} \end{equation}
The chain rule quickly shows that we can characterize complex-harmonic
morphisms $A^C \ra \CC$ as those complex-analytic
maps which pull back complex-analytic functions to complex-harmonic
functions (cf.\ Theorem {\bf\ref{fund:th:char}}).  Note also that, since these
equations are invariant when $\phi$ is composed with a holomorphic
map, we can replace $\CC$ by any Riemann surface.

\begin{proposition} \label{cxha:prop:cxha}
{\rm (i)} Let\/ $\phi:A^C \ra N^2$ be a complex-harmonic morphism from an
open subset\/ $A^C$ of\/ $\CC^4$ to a Riemann surface.  Then, for
any\/ ${\vec{p}} \in A^C$,
\begin{itemize}
\item[{\rm (a)}] $\phi\vert_{A^C\cap \RR^4_{\vec{p}}}$ is a harmonic morphism
(with respect to the standard Euclidean metric);

\item[{\rm (b)}] $\phi\vert_{A^C\cap \MM^4_{\vec{p}}}$ is a harmonic morphism
(with respect to the standard Minkowski metric).
\end{itemize}

{\rm (ii)} All harmonic morphisms from open subsets of\/ $\RR^4$ to Riemann
surfaces and all real-analytic harmonic morphisms from open subsets
of\/ $\MM^4$ to Riemann surfaces arise in this way.
\end{proposition}

\begin{proof}
(i) Immediate from the equations.

(ii) As noted above, any harmonic morphism from an open subset of
$\RR^4$ to a Riemann surface is real analytic.  By analytic
continuation, this is the restriction of a complex-analytic map on an
open subset of $\CC^4$, and this complex-analytic map is
complex-harmonic.  The $\MM^4$ case is similar except that real analyticity is no longer automatic and must be assumed.
\end{proof}

Then we have a complex-analytic version of Proposition
\ref{cxha:prop:R4muhm}, as follows.

\begin{proposition} \label{cxha:prop:C4muhm}
Let\/ $\mu:A^C \to \CC$ be a complex-analytic function from an open subset
of\/ $\CC^4$.  Suppose that\/ $\mu$ satisfies\/ {\rm(\ref{cxha:EC})}, i.e.,
suppose that\/
$\mu$ represents a complex-analytic foliation by null planes.  Then\/
$\mu$ is a complex-harmonic morphism.
\qed \end{proposition}

\section{Harmonic morphisms and shear-free ray congruences}
\label{cxha:sec:ha-morph-SFR}
We have just shown that a Hermitian structure, or any of the related
distributions \eqref{cxha:Q}, defines a harmonic morphism.  We now give some converses.
We first of all reformulate Proposition {\bf\ref{twistor:prop:Kerr}}.
By `fibre component' we shall mean `connected component of a fibre'.

\begin{theorem} \label{cxha:th:hamoJ} Let\/ $\phi:A^4
\to N^2$ be a harmonic morphism
from an open subset\/ $A^4$ of\/ $\RR^4$ to a Riemann surface,  with\/ $\dd\phi$ nowhere zero.  Then there
exists a Hermitian structure\/ $\JJ$ on\/ $A^4$ which is parallel
along each fibre component of\/ $\phi$. Further, for
any\/ ${\vec{p}} \in A^4$\,, there is a neighbourhood\/ $A_1^4$ of\/ ${\vec{p}}$ in\/ $A^4$
and a holomorphic map\/ $\rho:V \to \Ci$ from an open subset\/ $V$ of\/
$N^2$ such that\/ $\mu = \rho \circ \phi$ represents\/ $\JJ$ on\/ $A_1^4$.
\qed \end{theorem}

We have the following analogue for complex-harmonic morphisms:

\begin{theorem} \label{cxha:th:cxhamoalpha} Let\/ $\phi:A^C \ra N^2$ be
a complex-harmonic morphism, from an open subset\/ $A^C$ of\/ $\CC^4$
to a Riemann surface, with\/ $\dd\phi$ nowhere zero.  Then there exists
a holomorphic foliation\/ $\Ff$ of\/ $A^C$ by\/ $\alpha$-planes or
by\/ $\beta$-planes such that each fibre component of\/ $\phi$ is the union of parallel null planes of\/
$\Ff$. Further, for any\/ ${\vec{p}} \in A^C$, there is a neighbourhood\/ $A_1^C$
of\/ ${\vec{p}}$ in\/ $A^C$, and a holomorphic map\/ $\rho:V \ra \Ci$ from an
open subset\/ $V$ of\/ $N^2$, such that\/ $\mu = \rho \circ \phi$
represents\/ $\Ff$ on\/ $A_1^C$.
\end{theorem}

\begin{proof}
Let ${\vec{p}} \in A^C$.  By Proposition \ref{cxha:prop:cxha},  \ $\phi$
restricts to a harmonic morphism on $A^4 = A^C \cap \RR^4_{\vec{p}}$
which is easily seen to be submersive.  By Theorem \ref{cxha:th:hamoJ},
$\phi\vert_{A^4}$ is holomorphic with respect to some Hermitian
structure $\JJ$ which is constant along each fibre component
of $\phi$.  By replacing $q_2$ with $\wt{q}_2$,
if necessary, we can assume that $\JJ$ is positively oriented.
Represent $\JJ$ by the map $\mu:A^4 \ra \Ci$; then, since $\phi$ is holomorphic and
$\JJ$ has $(0,1)$-tangent space with basis (\ref{cxha:tang}), we have at
all points of $A^4$,
\begin{equation}
\frac{\pa\phi}{\pa \wt{q}_1} - \mu \frac{\pa\phi}{\pa q_2} = 0 \,,
\quad
\frac{\pa\phi}{\pa \wt{q}_2} + \mu \frac{\pa\phi}{\pa q_1} = 0 
\,.\label{cxha:muphi}\end{equation}
Since $\JJ$ is integrable, we have also equation (\ref{cxha:EC}):
\begin{equation}
\frac{\pa\mu}{\pa \wt{q}_1} - \mu \frac{\pa\mu}{\pa q_2} = 0 \,,
\quad
\frac{\pa\mu}{\pa \wt{q}_2} + \mu \frac{\pa\mu}{\pa q_1} = 0  
\,. \label{cxha:ER2} \end{equation}
Extend $\mu$ to
$A^C$ by insisting that (\ref{cxha:muphi}) holds; note that $\mu$ is well defined
by one of these equations---
since the hypothesis that $\dd\phi$ is nowhere zero ensures that
not all the partial derivatives of $\phi$ can vanish
simultaneously---and is a complex-analytic function.
By analytic continuation, (\ref{cxha:ER2}) holds at all points of $A^C$
so that $\mu$ defines a holomorphic foliation $\Ff$ by
$\alpha$-planes.  By (\ref{cxha:muphi}), $\phi$ is constant along any
$\alpha$-plane of $\Ff$, so that each fibre of $\phi$ is the union
of $\alpha$-planes of $\Ff$.  Further, $\JJ$ and so $\mu\vert_{A^4}$ is
constant along the fibre components of
$\phi\vert_{A^4}$; by analytic continuation, $\mu:A^C \ra \Ci$ is
constant along the fibre components of
$\phi:A^C \ra \Ci$, so that the $\alpha$-planes of $\Ff$ making up
a fibre component of $\phi$ are all parallel.

Lastly, since $\mu$ is constant on the leaves of the foliation
given by the fibres of $\phi$, it factors through local leaf
spaces as $\mu = \rho \circ \phi$.  Since $\phi$ and $\mu$ are
both holomorphic (with respect to $\ii$ and $\JJ$), $\rho$ must be
holomorphic.
\end{proof}

We have the following analogue for Minkowski space proved in a similar way.

\begin{theorem} \label{cxha:th:hamoSFRanal}
Let\/ $\phi:A^M \ra N^2$ be a real-analytic harmonic
morphism from an open subset\/ $A^M$ of Minkowski\/ $4$-space\/ $\MM^4$ to a Riemann
surface, with $\dd\phi$ nowhere zero.  Then there is a shear-free ray
congruence\/ $\ell$ on\/ $A^M$ such that each fibre component
of\/ $\phi$ is the union of parallel null lines of\/ $\ell$.
Further, for any\/ ${\vec{p}} \in A^M$, there is a neighbourhood\/ $A_1^M$ of
${\vec{p}}$ in\/ $A^M$, and a holomorphic map\/ $\rho:V \ra \Ci$ from an open
subset\/ $V$ of\/ $N^2$, such that\/ $\mu = \rho \circ \phi$ represents\/
$\ell$ on\/ $A_1^M$.
\qed \end{theorem}

Note that in the real and complex cases the condition `$\dd\phi$ nowhere zero' is equivalent to `submersive'. This is not so
in the Minkowski case where $\phi$ may be degenerate (Definition
{\bf\ref{semi:def:HWCdeg})}.  In this case we can be more precise as follows.

\begin{corollary} \label{cxha:cor:deg}
Let\/ $\phi:A^M \to N^2$ be a real-analytic harmonic morphism from an
open subset\/ $A^M$ of\/ $\MM^4$ to a Riemann surface.
Suppose that\/ $\phi$ is degenerate at\/ ${\vec{p}}$ (so that\/ $\dd\phi_{\vec{p}}$ has rank\/ $1$).
Then there is a unique null
direction\/ $W_{\vec{p}} \in T_{\vec{p}}\MM^4$ such that\/ $W_{\vec{p}} \subset \ker
\dd\phi_{\vec{p}}$. Furthermore, $\ker \dd\phi_{\vec{p}} = W_{\vec{p}}^{\perp}$.  If,
further, at each point\/ ${\vec{q}}$ in the fibre component
through\/ ${\vec{p}}$, \ $\phi$ is degenerate,
then that fibre component is an open subset of
the affine null 3-space tangent to\/ $W_{\vec{p}}^{\perp}$.
\end{corollary}

\begin{proof}
Since $\dd\phi_{\vec{p}}$ has rank $1$, $\ker\dd\phi_{\vec{p}}$ is three-dimensional.
By Lemma {\bf\ref{semi:lem:HWCdeg}},
$(\ker\dd\phi_{\vec{p}})^{\perp} \subset \ker\dd\phi_{\vec{p}}$,   so
$(\ker \dd\phi_{\vec{p}})^{\perp}$ is one-dimensional and null.  We set
$W_{\vec{p}} = (\ker\dd\phi_{\vec{p}})^{\perp}$,  so $\ker\dd\phi_{\vec{p}}= W^{\perp}_{\vec{p}}$.

To prove uniqueness of $W_{\vec{p}}$, suppose that $W'_{\vec{p}} \subset \ker \dd\phi_{\vec{p}}$
is another null direction.  Then $W_{\vec{p}}$ and $W'_{\vec{p}}$ together span a $2$-dimensional
subspace of null directions in $W^{\perp}_{\vec{p}}$, which is impossible.

This means that the distribution ${\vec{p}} \mapsto W_{\vec{p}}$ must be tangent
to the shear-free ray congruence of Theorem \ref{cxha:th:hamoSFRanal},
and so each $W_{\vec{p}}$ is parallel for all ${\vec{p}}$ in a fibre component.
The last assertion follows from the fact that the connected
component of the fibre is three-dimensional and has every tangent
space parallel to $W_{\vec{p}}^{\perp}$.
\end{proof}

We can also give the following version of Theorem \ref{cxha:th:hamoSFRanal}
for harmonic morphisms which are only smooth,
though this needs the stronger hypothesis of submersivity and the
conclusion is slightly weaker.

\begin{theorem} \label{cxha:th:hamoSFRsmooth}
Let\/ $\phi:A^M \ra N^2$ be a smooth submersive harmonic map
from an open subset\/ $A^M$ of\/ $\MM^4$ to a Riemann
surface.  Then there there is an open dense set of points\/ $A_0^M$
and a shear-free ray
congruence\/ $\ell$ on\/ $A_0^M$ such that each fibre component
of\/ $\phi$ is the union of parallel null lines of\/ $\ell$.
Further, for any\/ ${\vec{p}} \in A_0^M$, there is a neighbourhood\/ $A_1^M$ of
${\vec{p}}$ in\/ $A_0^M$ and a holomorphic map\/ $\rho:V \ra \Ci$ from an open
subset\/ $V$ of\/ $N^2$ such that\/ $\mu = \rho \circ \phi$ represents\/
$\ell$.
\end{theorem}

\begin{proof}
Let ${\vec{p}} \in U$.  Let $v,w,z,\ov{z}$ be any null coordinates, i.e.,
coordinates in which the metric
is $\dd v \,\dd w + \dd z \,\dd\ov{z}$.
Then the harmonic morphism equations \eqref{cxha:hamo} read
\begin{eqnarray}
\phi_v \phi_w + \phi_z \phi_{\ov{z}} & = & 0 \,, \label{cxha:nullHWC} \\
\phi_{vw} + \phi_{z\ov{z}} & = & 0 \,. \label{cxha:nullha}
\end{eqnarray}

Now choose the null coordinates such that
$\pa/\pa z, \pa/\pa \ov{z}$ span the horizontal space at ${\vec{p}}$ and
$\pa/\pa v \,, \pa/\pa w$ span the vertical space,
thus $\phi_v({\vec{p}}) = \phi_w({\vec{p}}) = 0$. Then from
(\ref{cxha:nullHWC}), $\phi_z({\vec{p}})\, \phi_{\ov{z}}({\vec{p}}) = 0$.
Without loss of generality we can assume that
\begin{equation} \label{cxha:neq}
\phi_z({\vec{p}}) \neq 0\,,
\end{equation}
so that $\phi_{\ov{z}}({\vec{p}}) = 0$.
Now differentiation of (\ref{cxha:nullHWC}) shows that
all first derivatives of $\phi_{\ov{z}}$ are zero at ${\vec{p}}$.
On differentiating (\ref{cxha:nullHWC}) with respect to $z$ then $\ov{z}$,
evaluating at ${\vec{p}}$ and using (\ref{cxha:neq}) we obtain
$\phi_{z\bar{z}\bar{z}}({\vec{p}}) = 0$ then (\ref{cxha:nullha}) gives
$\phi_{\ov{z} vw}({\vec{p}}) = 0$.
On differentiating (\ref{cxha:nullHWC}) with respect to $v$ then $w$ and
evaluating at ${\vec{p}}$ we obtain $\phi_{vv}({\vec{p}})\,\phi_{ww}({\vec{p}}) = 0$, so that
$$
\phi_{vv}({\vec{p}}) = 0 \quad \mbox{or} \quad \phi_{ww}({\vec{p}}) = 0 \,.
$$
Suppose that $\phi_{vv}({\vec{p}}) = 0$.  Then, at ${\vec{p}}$, the function
$\mu = \ii\phi_v/\phi_z$
satisfies
\begin{eqnarray*}  
\frac{\pa\mu}{\pa v} + \ii\mu \frac{\pa\mu}{\pa z} & = & 0 \,,
	\quad \mbox{and} \\
\frac{\pa\mu}{\pa \ov{z}} -\ii\mu \frac{\pa\mu}{\pa w} &  = & 0 \,.
\end{eqnarray*}
comparing these with \eqref{cxha:EM}
shows that $\mu$ represents a shear-free ray congruence $W$ at ${\vec{p}}$.
If, instead, $\phi_{ww}({\vec{p}}) = 0$; interchange $v$ and $w$ in the above;
 again, we get a shear-free ray congruence $W$ at ${\vec{p}}$.

Now let $A_1 = \{{\vec{p}} \in A^M : \phi_{vv} = 0\}$ and $A_2 =
\{{\vec{p}} \in A^M : \phi_{ww} = 0\}$.  Then $A_1 \cup A_2 = A^M$ and we obtain
shear-free ray congruences on the open sets $\mbox{int}(A_1)$ and
$\mbox{int}(A_2)$;
by the Baire category theorem (Sims 1976, Section 6.4), their union is dense in $A^M$.
\end{proof}

\begin{remark} \rm
 On the open set $\mbox{\rm int}(A_1) \cap
\mbox{\rm int}(A_2)$ (which may be empty) both $V$
and $W$ are shear-free ray congruences; the fibres of $\phi$ are spanned by
$V$ and $W$ and are totally geodesic.
\end{remark}

\begin{example} \label{cxha:ex:HopfMin}
It can easily be checked that, for any holomorphic
function $f$, the composite map $\phi = f\bigl((x_2+\ii x_3)/(x_1+t)\bigr)$
is a complex-valued harmonic morphism on $\MM^4 \setminus \{(t,x_1,x_2,x_3): x_1+t = 0 \}$.
In fact, $\phi = f(\ii\mu)$ where $\mu$ is given by
(\ref{cxha:HopfMin}).  If $\dd f$ is nowhere zero, then $\phi$ is submersive,
and the associated shear-free ray
congruence $\ell$ of Theorem \ref{cxha:th:hamoSFRanal} (or
\ref{cxha:th:hamoSFRsmooth}) is that in Example \ref{cxha:ex:Hopf}.

Similar examples can be given for Theorems \ref{cxha:th:hamoJ} and
\ref{cxha:th:cxhamoalpha}.
\end{example}

\section{Compactifications and twistorial formulation}
\label{cxha:sec:twistor}

\begin{definition} {\rm (LeBrun 1983)} \label{cxha:def:holconfstr}
A\/ {\rm holomorphic conformal structure} on a complex manifold\/ $M$ is
a holomorphic line subbundle\/ $L$ of the symmetric square\/ $\odot^2 T^*_{1,0}M$ such that,
for every\/ $p \in M$, any non-zero element of\/ $L_p$ defines a (complex-sym\-metric
bilinear) non-degenerate inner product on\/ $T^{1,0}_pM$.
\end{definition}

Any local nowhere-zero holomorphic section of $L$ is a holomorphic
metric.  If $g_1^C, \ g_2^C$ are any two such holomorphic sections
defined on an open subset $A$ of $M$ then they are holomorphically conformally
equivalent, i.e.,  $g_1^C = \nu\, g_2^C$ for some holomorphic map $\nu:A \to \CC \setminus \{0\}$.  
Conversely, a holomorphic conformal
structure can be specified by giving a collection of holomorphic
metrics whose domains cover $M$ such that any two are conformally related
on the intersection of their domains. 

Let $\Cf = G_2(\CC^4)$ denote the Grassmannian of
all two-dimensional complex subspaces of $\CC^4$ equipped with its standard
complex structure.  Then $\CC^4$ can be embedded holomorphically in
$\Cf$ by the mapping given in standard null coordinates by
\begin{equation}
c: \vec{x} = (q_1, \wt{q}_1, q_2, \wt{q}_2) \mapsto \mbox{ column space of }
A = \left(\begin{array}{rr} 1\ & 0\ \\ 0\ & 1\ \\ q_1  &
		-\wt{q}_2 \\ q_2 & \wt{q}_1 \end{array} \right) \,.
\label{cxha:emb} \end{equation}
We shall call $c$ the \emph{standard coordinate chart} for
$\Cf$; we identify $\CC^4$ with its image under $c$.  Points of $\Cf \setminus \CC^4$ will be called
\emph{points at infinity}.

Note that the matrix $A$ is of the form
$A = \left(\begin{array}{cc} I \\ Q \end{array} \right)$
where $I$ denotes the $2\times2$-identity matrix and
\begin{equation} \label{cxha:Qcomplex}
Q = \left( \begin{array}{rr} 	q_1 & -\wt{q}_2 \\
			q_2 & \wt{q}_1 \end{array} \right).
\end{equation}
Write $\vec{x} = (x_0,x_1,x_2,x_3) = (q_1,q_2) \in \RR^4 = \CC^2$, then
\begin{equation} \label{cxha:Qreal}
Q = \left( \begin{array}{rr} 	q_1 & -\ov{q}_2 \\
			q_2 & \ov{q}_1 \end{array} \right)
  = \left( \begin{array}{rr} x_0 + \ii x_1 & -x_3 + \ii x_4 \\
		x_3 + \ii x_4 & x_0 - \ii x_3 \end{array} \right)   
\end{equation}
which is a standard representation of the quaternion $q_1 + q_2 \jj
=x_0 + \ii x_1 + \jj x_2 + \kk x_3$; we shall call such a matrix $Q$
\emph{quaternionic}. 
Thus $Q$, and so $A$, represents a point of $\RR^4$ if and only if $Q$
is quaternionic.
If, on the other hand, $\vec{x} = (t,x_1,x_2,x_3) \in \MM^4$, then
\begin{equation} \label{cxha:QMin}
Q = \left( \begin{array}{rr} -\ii t + \ii x_1 & -x_3 + \ii x_4 \\
		x_3 + \ii x_4 & -\ii t - \ii x_3 \end{array} \right) ,
\end{equation}
so that $Q$, and so $A$, represents a point of $\MM^4$ if and only if $Q$
is skew-Hermitian.

We give $\Cf$ a holomorphic conformal structure, as follows.
Given any complex chart $c:U \to
\Cf$, choose holomorphic maps $X,Y:U \to \CC^4$ such that 
$$
c({\vec{q}}) =
\mbox{span}\{X({\vec{q}}),Y({\vec{q}})\} \qquad ({\vec{q}} \in U).
$$  
Choose a non-zero element
$\omega^*$ of $\La^4\CC^4$.  For each ${\vec{p}} \in U$ and $w \in T_{\vec{p}}U$ define
$g(w,w)$ to be the unique complex number such that
\begin{equation} \label{cxha:holconf-Cf}
X \wedge Y \wedge \dd X(w) \wedge \dd Y(w) = g(w,w)\,\omega^*.
\end{equation}
On defining $g(v,w)$ for $v,w \in T_pU$ by  polarization: 
$$
g(v,w) =
\half \big\{g(v+w,v+w) - g(v,v) - g(w,w) \big\}
$$ 
we obtain
a holomorphic metric on $U$, and so on
$c(U)$.  It is well defined up to multiplication by nowhere-zero
holomorphic functions, and so defines a holomorphic conformal structure.

For the standard chart (\ref{cxha:emb}) with $X, Y$ given by the columns
of the matrix $A$, this gives the standard holomorphic metric
$\dd q_1 \,\dd\wt{q}_1 + \dd q_2 \,\dd\wt{q}_2$; it follows that the
standard chart is conformal as well as holomorphic.

Now, at any point $\vec{p} \in \Cf$, given any orthogonal basis
$\{ e_0,e_1,e_2,e_3 \}$ of $T_{\vec{p}}\RR^4_{\vec{p}}$ of vectors of the same length,
the two-dimensional subspace 
$$
\Pi_{\vec{p}} = \mbox{span}\{e_0+\ii e_1, e_2 + \ii
e_3\}
$$ 
of $T_{\vec{p}}\Cf$ is null; we call $\Pi_{\vec{p}}$ an \emph{$\alpha$-plane}
(respectively, \emph{$\beta$-plane}) according as the basis $\{ e_0,e_1,e_2,e_3\}$
is positively (respectively, negatively) oriented.  This agrees with our
previous definition (Section \ref{cxha:sec:nullplanes}) if ${\vec{p}} \in
\CC^4$.  As before, any null two-dimensional subspace of $T_{\vec{p}}\Cf$ is an
$\alpha$-plane or a $\beta$-plane.

We now determine the \emph{null surfaces} in $\Cf$, i.e. the surfaces whose tangent
spaces are null.
Fix $w \in \CP^3$, thus $w$ is a one-dimensional subspace of $\CC^4$; then set
\begin{equation} \label{cxha:wtilde}
\wt{w} = \{p \in G_2(\CC^4) : \mbox{the plane } p \mbox{ contains } w \} \,,
\end{equation}
this defines a surface in $\Cf \equiv G_2(\CC^4)$.  Now the formula
(\ref{cxha:holconf-Cf}) shows that all its tangent
spaces are null; it is easily seen that they are all $\alpha$-planes.
All such surfaces are given this way. If it is not at infinity, we easily see
that $\wt{w}$ is the image under $c$ of an affine $\alpha$-plane
of $\CC^4$  so we shall continue to call the surfaces
\emph{affine $\alpha$-planes}.
(Affine $\beta$-planes have a similar definition and can be described as
the sets of all two-dimensional subspaces of $\CC^4$ contained in a given
three-dimensional subspace.)

To understand all this better, embed $G_2(\CC^4)$ in $\CP^5$ by the
Pl\"ucker embedding $pl$ which sends the plane spanned by the vectors $X$ and $Y$ to the
point $[X \wedge Y] \in P(\La^2\CC^4) \cong \CP^5$.  Here 
$[ \ ] = \pi( \ )$ where $\pi:\La^2\CC^4 \setminus \{\vec{0}\} \to P(\La^2\CC^4)$ is
the standard projection.  The image of $pl$ is the complex quadric
$$
\Qq^C = \{ [\omega] : \omega \in \La^2\CC^4 : \omega \wedge
\omega = 0 \} \,,
$$
 the condition $\omega \wedge \omega = 0$ expressing the
decomposibility of $\omega$.  Identify $\La^2\CC^4$ with $\CC^6$ in
the standard way; then the components of $pl$ are given by
the $2 \times 2$ minors of the matrix $A$. Thus the
composition $j = pl \circ c:\CC^4 \to \CP^5$ is given by
\begin{eqnarray} \label{cxha:j}
j(q_1, \wt{q}_1, q_2, \wt{q}_2) & = & 
[z_{12},\,z_{13},\,z_{14},\,z_{23},\,z_{24},\,z_{34}] \nonumber \\
  &  = & 
[1,\,-\wt{q}_2,\,\wt{q}_1,\,-q_1,\,-q_2,\, q_1\wt{q}_1 + q_2\wt{q}_2]
\end{eqnarray}
and $\Qq^C$ has equation $z_{12}z_{34} - z_{13}z_{24} + z_{14}z_{23} = 0$.
Change the coordinates linearly to
\begin{equation} \label{cxha:xi}
\xi_0 = z_{12}+z_{34}, \ \xi_1 = z_{12}-z_{34}, \ \xi_2 = z_{14}-z_{23}, \ 
\xi_3 = \ii(z_{14}+z_{23}), \ \xi_4 = -(z_{13}+z_{24}), \ 
\xi_5 = -\ii(z_{13}-z_{24}).
\end{equation}
In these coordinates, $\Qq^C$ has equation
\begin{equation} \label{cxha:QR}
{\xi_0}^{\! 2} - {\xi_1}^{\! 2} - {\xi_2}^{\! 2} - {\xi_3}^{\! 2} - {\xi_4}^{\! 2} - {\xi_5}^{\! 2} = 0\,,
\end{equation}
and the restriction of $j$ to $\RR^4$ reads
$$
\vec{q}= (x_0,x_1,x_2,x_3) \mapsto [1+\abs{\vec{q}}^2, 1-\abs{\vec{q}}^2, 2x_0, 2x_1, 2x_2,2x_3] =
[1, \si^{-1}(\vec{q})]
$$
where $\si^{-1}:\RR^4 \to S^4 \subset \RR^5$ is the inverse of
stereographic projection.
This shows that (i) the closure $\Rf$ of the image of $\RR^4$ consists of  the points of
$\Qq^C$ with real coordinates $\xi_i$, this is the quadric $\Qq^R$ in $\RP^5$ with
equation (\ref{cxha:QR}); (ii) $\Qq^R$ is conformally equivalent to $S^4$;
(iii) the set of points at infinity of $\RR^4$, i.e., the set
$\Rf \setminus \RR^4$, consists of the single point $[1,-1,0,0,0,0]$.

In the coordinates
\begin{equation} \label{cxha:wtxi}
\wt{\xi}_i = \xi_i \ (i \neq 2), \quad
\wt{\xi}_2 =  \ii \xi_2
\end{equation}
the equation of $\Qq^C$ becomes
\begin{equation} \label{cxha:QM}
g(\wt{\vec{\xi}},\wt{\vec{\xi}}) \equiv
\wt{\xi}_0^{\,\,2} - \wt{\xi}_1^{\,\,2} + \wt{\xi}_2^{\,\,2} - \wt{\xi}_3^{\,\,2}
- \wt{\xi}_4^{\,\,2} - \wt{\xi}_5^{\,\,2} = 0
\end{equation}
and the restriction of $j$ to $\MM^4$ reads
\begin{equation} \label{cxha:j|M4}
\vec{q}= (t,x_1,x_2,x_3) \mapsto [1+\abs{\vec{q}}_1^2, 1-\abs{\vec{q}}_1^2, 2t, 2x_1, 2x_2, 2x_3]\,.
\end{equation}
Hence the closure $\Mf$ of the image of $\MM^4$ is the points of $\Qq^C$ with
the coordinates $\wt{\xi}_i$ real valued; this is the quadric $\Qq^M$
in $\RP^5$ with equation (\ref{cxha:QM}).  The set of points at infinity of
$\MM^4$, i.e., the set $\Mf \setminus \MM^4$, consists of the null cone
$\{[1,-1,\vec{q}]: \vec{q} \in \MM^4, \ \abs{\vec{q}}_1^2 = 0 \}$.

\begin{remark} \label{cxha:rem:stereoM4}  \rm
The intersection of $\Mf$ with the chart $\wt{\xi}_0 = 1$ is the pseudosphere given by
$S^4_1 = \{(\wt{\xi}_1,\wt{\xi}_2,\wt{\xi}_3,\wt{\xi}_4,\wt{\xi}_5) \in \RR^5: \wt{\xi}_1^{\,\,2} - \wt{\xi}_2^{\,\,2}
+ \wt{\xi}_3^{\,\,2} + \wt{\xi}_4^{\,\,2} + \wt{\xi}_5^{\,\,2} = 1 \}$.  The map $j:\MM^4 \to \Qq^M$
restricts to the map $\MM^4 \setminus \{\vec{q} \in \MM^4: \ \abs{\vec{q}}_1^2 = -1 \}$
given by
$$
\vec{q} \mapsto \frac{1}{1+\abs{\vec{q}}^2_1} \big( 1-\abs{\vec{q}}^2_1, \vec{q} \big) \,.
$$
This is a higher-dimensional version of the stereographic projection
discussed in Example {\bf\ref{semi:ex:char-coords}}.
\end{remark}

Note that the mapping from $S^1\times S^3$ given by 
\begin{equation} \label{cxha:S1S3QM}
 \big( (\wt{\xi}_0, \wt{\xi}_2),
(\wt{\xi}_1,\wt{\xi}_3,\wt{\xi}_4,\wt{\xi}_5) \big) \mapsto
[\wt{\xi}_0,\wt{\xi}_1,\wt{\xi}_2,\wt{\xi}_3,\wt{\xi}_4,
\wt{\xi}_5]
\end{equation}
is a double covering; it is also
conformal if we give $S^1 \times S^3$ the Lorentian product metric
$-g^{S^1}+g^{S^3}$. 

Think of $S^1 \times S^3$ as $\big\{ (z_1,z_2,z_3) \in
\CC \times \CC^2 : \abs{z_1}^2 =
1, \ \abs{z_2}^2 + \abs{z_3}^2 = 1 \big\}$, and define an equivalence
relation $\sim$ on $S^1 \times S^3$
by $(z_1, z_2,z_3) \sim \pm(z_1, z_2, z_3)$. 
The map (\ref{cxha:S1S3QM}) factors to a conformal diffeomorphism
\begin{equation} \label{cxha:S1S3simQM}
S^1 \times S^3 \big/ \!\! \sim \ \to \Qq^M.
\end{equation}
Note that the mapping
$$
S^1 \times S^3 \ni (z_1,z_2,z_3)
   \mapsto ({z_1}^{\! 2}, z_1 z_2, z_1 z_3) \in S^1 \times S^3
$$
factors to a diffeomorphism of $S^1 \times S^3 \big/ \!\! \sim$
to $S^1 \times S^3$.

\begin{remark} \label{cxha:rem:pa/pa-t}  \rm
We give $S^1 \times S^3$ the time- and space-orientations induced from
the standard orientations on $S^1$ and $S^3$.  These factor to $S^1
\times S^3\big/ \!\! \sim$.  
The map $j$ given by \eqref{cxha:j|M4} sends $\vec{x} \in \RR^3$ to
the equivalence class of the point
$\big( (1,0), \si^{-1}(\vec{x}) \big) \in S^1 \times S^3\big/\!\! \sim$\,.  The differential of $j$ maps
the vector $\pa/\pa t$ at $\vec{x}$, which is normal to $\RR^3$ in $\MM^4$, to $(u,0) \in T_{(1,0)} S^1 \times
T_{\si^{-1}(\vec{x})} S^3 $ where $u$ denotes the unit positive tangent
vector.  By a slight abuse of notation, we shall continue to denote
$(u,0)$ by $\pa/\pa t$ even at points at infinity.  Then any
null direction in $\Mf \cong S^1 \times S^3\big/ \!\! \sim$
is spanned by a vector $\pa/\pa t + \UU = (u,\UU) \in TS^1 \times TS^3$, as
in $\MM^4$.
\end{remark}

We give $\Qq^C$ a holomorphic conformal structure, as follows.  Let $g^Q$
be a symmetric bilinear form on $\CC^6$ such that
$\Qq^C$ has equation $g^Q(\vec{\xi},\vec{\xi}) = 0$, i.e., (\ref{cxha:QM}) in the
coordinates $\wt{\xi}_i$.
Let $X:A \to T^{1,0}\Qq^C$ be a holomorphic vector field defined on an
open subset of $\Qq^C$, and let $\wh{X}:A \to
\CC^6$ be a holomorphic lift of it, i.e.,
$\dd\pi(\wh{X}) = X$ where $\pi: \CC^6 \setminus \{\vec{0}\} \to \CP^5$ is the natural projection.
Set $g(X,X) = g^Q(\wh{X}, \wh{X})$.  Note that $\wh{X}$ is defined up to a
transformation $\wh{X}'_p = \la(p)\wh{X}_p + \mu(p) \wh{p}$ \
$(p \in A)$
where $\la:A \to \CC \setminus \{0\}$ is holomorphic and
$\pi(\wh{p})= p$. Since
$g(\wh{p},\wh{p}) = g(\wh{p}, \wh{X}) = 0$, the first being the
equation of $\Qq^C$ and the second the tangency condition, we have
$g^Q(\wh{X}', \wh{X}') = \la^2 g^Q(\wh{X}, \wh{X})$.  Hence $g$ is
well defined up to multiplication by a nowhere-zero holomorphic function;
it thus  defines a holomorphic conformal structure on $\Qq^C$.  Clearly
this corresponds under $pl$ to the holomorphic
conformal structure defined above on $G_2(\CC^4)$.  Thus $j$ is
a holomorphic conformal inclusion.

More generally, for any ${\vec{p}} \in \CC^4$,
let $\Rf_{\vec{p}}$, $\Rt_{\vec{p}}$, $\Mf_{\vec{p}}$ denote the closures of the slices
$j(\RR^4_{\vec{p}})$, $j(\RR^3_{\vec{p}})$, $j(\MM^4_{\vec{p}})$ in $\Cf$, respectively;
we write $\Rf$, $\Rt$, $\Mf$ if $\vec{p} = \vec{0}$.
Then the holomorphic conformal structure on $\Cf$ induces real conformal
structures on these closures; clearly $c$ restricts to conformal
inclusions $\RR^4_{\vec{p}} \to \Rf_{\vec{p}}$, \ $\RR^3_{\vec{p}} \to \Rt_{\vec{p}}$ and 
$\MM^4_{\vec{p}} \to \Mf_{\vec{p}}$\,. 

Let ${F}_{1,2}$ be the complex flag manifold
$$
\{ (w,p) \in \CP^3 \times G_2(\CC^4) :
\mbox{ the line } w \mbox{ lies in the plane } p \} \,;
$$
this is called the \emph{correspondence space}.  The restrictions to ${F}_{1,2}$ of the
natural projections from $\CP^3 \times G_2(\CC^4)$  define a \emph{double holomorphic
fibration}:
\begin{equation}
\begin{array}{ccccc} 
       &              & (w,p) \in {F}_{1,2} &              &   \\
       & \mu \swarrow &               & \searrow \nu &   \\
w \in \CP^3 &        &                &  & p \in G_2(\CC^4) \equiv \Cf
\end{array}
\label{cxha:doublefibr} \end{equation}
For any $w \in \CP^3$, set
\begin{equation} \label{cxha:wt}
\wt{w} = \nu \circ \mu^{-1}(w) \,;
\end{equation}
as in (\ref{cxha:wtilde}), $\wt{w}$  
is an $\alpha$-plane in $\Cf$, which we call the $\alpha$-plane
\emph{determined} or \emph{represented by} $w$.
Conversely, for any point $p \in \Cf$ we write
\begin{equation} \label{cxha:wh}
\wh{p} = \mu \circ \nu^{-1}(p) \,;
\end{equation}
then $\wh{p}$ is a projective line (i.e., a
one-dimensional projective
subspace) of $\CP^3$ which represents the $\CP^1$'s-worth of
$\alpha$-planes through $p$.

Explicitly, in the standard null coordinates $c:\CC^4 \hookrightarrow \Cf$
given by (\ref{cxha:emb}), we have $(w,p) \in {F}_{1,2}$ if and only the
following \emph{incidence relations} are satisfied:
\begin{equation}
w_0 q_1 - w_1 \wt{q}_2 = w_2 \,, \quad w_0 q_2 + w_1 \wt{q}_1 =
w_3 \,;
\label{cxha:I} \end{equation}
indeed,  these express the condition that $w$ be a linear combination
of the columns of the matrix (\ref{cxha:emb}) and so lie in the plane
represented by the point $p \in \Cf$.  Now note that,
for any $[w_0,w_1,w_2,w_3]$ with $[w_0,w_1] \neq [0,0]$,
(\ref{cxha:I}) defines an affine $\alpha$-plane in $\CC^4$: indeed,
its tangent space at any point is the set of
vectors annihilated by $\{w_0 \dd q_1 - w_1 \dd\wt{q}_2\,,\, w_0 \dd q_2 +
w_1 \dd\wt{q}_1 \}$ and so is given by (\ref{cxha:NT0}).  Points of
$\CP^3$ on the projective line $\CP^1_0 = \{ [w_0,w_1,w_2,w_3]
:[w_0, w_1] = [0,0] \}$ correspond to affine $\alpha$-planes \emph{at
infinity}, i.e., in $\Cf \setminus \CC^4$.  Thus $\CP^3$
parametrizes all $\alpha$-planes in $\Cf$ and (\ref{cxha:I})
\emph{expresses the condition that the point  in $\CC^4$ with
null coordinates $(q_1, \wt{q}_1,q_2, \wt{q}_2)$
lies on the $\alpha$-plane $\wt{w}$ determined by
$w = [w_0,w_1,w_2,w_3] \in \CP^3 \setminus \CP^1_0$}, i.e.,
(\ref{cxha:I}) is the equation of the $\alpha$-plane $\wt{w}$.

For $w = [w_0,w_1,w_2,w_3] \in \CP^3 \setminus \CC P^1_0$, the
$\alpha$-plane $\wt{w}$ given by (\ref{cxha:I}) intersects
$\RR^4 = \RR^4_0$ at $(q_1,q_2)$ where
\begin{equation}
\left. \begin{array}{lll}
                      w_0 q_1 - w_1 \ov{q}_2 & = & w_2 \\
\ov{w}_1 q_1 + \ov{w}_0 \ov{q}_2 & = &
\ov{w}_3 \end{array}
\right\} \,.
\label{cxha:IR} \end{equation}
This system has the unique solution
\begin{equation} 
\left( \begin{array}{l} q_1 \\ q_2 \end{array} \right) =
\frac{1}{\abs{w_0}^2 + \abs{w_1}^2}
\left( \begin{array}{l} \ov{w}_0 w_2 + w_1 \ov{w}_3 \\
                \ov{w}_0 w_3 - w_1 \ov{w}_2 \end{array}
\right)\,;
\label{cxha:int} \end{equation}
and so defines a map $\CP^3 \setminus \CP^3_0 \to \RR^4$
which is the projection map of the twistor bundle
of $\RR^4$ given by ({\bf\ref{twistor:CP2R4}}).  On mapping points of
$\CP^1_0$ to the point at infinity $\Rf \setminus \RR^4$,  
this map extends to a smooth map
\begin{equation} \label{cxha:CP3S4}
\pi = \pi_0:\CP^3 \ra \Rf \cong S^4 
\end{equation}
which is the projection map of the twistor bundle of $S^4$ as given in
({\bf\ref{twistor:CP3HP1}}).
Similary, for any ${\vec{p}} \in \CC^4$ we have a smooth map
$\pi_{\vec{p}}:\CP^3 \to \Rf_{\vec{p}} \cong S^4$ which sends $w \in \CP^3$ to the
unique point of intersection of the $\alpha$-plane $\wt{w}$
with $\Rf_{\vec{p}}$\,; on writing
${\vec{p}} = (a_1, \wt{a}_1, a_2, \wt{a}_2)$ in standard null coordinates,
a simple calculation shows that, in the chart (\ref{cxha:R4p})
for $\Rf_{\vec{p}}$\,,  this point of intersection is
given by (\ref{cxha:int}) with $w_2, \  w_3$ replaced by $w_2',
\ w_3'$ where
\begin{equation} \label{cxha:transl}
w_2' = w_2 - w_0 a_1 + w_1 \wt{a}_2, \quad
w_3' = w_3 - w_0 a_2 - w_1 \wt{a}_1
\end{equation}
(see Example \ref{cxha:ex:SL4C}
for another explanation of these formulae). 

The affine $\alpha$-plane (\ref{cxha:I})
intersects the Minkowski slice $\MM^4_{\vec{p}}$ if and only if it
intersects $\RR^3_{\vec{p}}$;  the intersection is then an (affine) null
line. For ${\vec{p}}={\vec{0}}$ this holds if and only if the point (\ref{cxha:int})
lies in $\RR^3$, i.e., $\Re q_1 = 0$, and this holds if and only if
$[w_0,w_1,w_2,w_3]$ lies on the real quadric
\begin{equation}
N^5 = \pi^{-1}(\RR^3) =
\bigl\{ [w_0,w_1,w_2,w_3] : w_0 \ov{w}_2 + \ov{w}_0
w_2 + w_1 \ov{w}_3 + \ov{w}_1 w_3 = 0 \bigr\} \subset
\CP^3.
\label{cxha:N5} \end{equation}
Points of $N^5 \setminus \CP^1_0$ thus represent affine null lines
of $\MM^4$; points of $\CP^1_0$ lie in $N^5$ and
represent the null lines at infinity.
Note that the map (\ref{cxha:CP3S4}) restricts to the projection map of a bundle
\begin{equation} \label{cxha:N5S3}
\pi:N^5 \to \Rt \cong S^3
\end{equation}
which gives the point of intersection of the null line with $\Rt$.
For general ${\vec{p}}$ we replace $N^5$ by $N^5_{\vec{p}} = \pi_{\vec{p}}^{-1}(\Rt_{\vec{p}})$.

Note that, on interchanging $q_2$ and $\wt{q}_2$ in (\ref{cxha:I}) we obtain the
standard parametrization of \emph{$\beta$-planes} by $\CP^3$.

The incidence relations (\ref{cxha:I}) define the following fundamental map which
gives the point of $\CP^3$ representing the $\alpha$-plane
through $(q_1, \wt{q}_1, q_2, \wt{q}_2)$ with direction
vector $\si^{-1}(\ii w_1/w_0)$: \begin{equation}
\begin{array}{rrl} \CC^4 \times \CP^1 & \lra   & \CP^3 \setminus
{\CC P^1_0} \ \subset \ \CP^3    \\
\iota:\big( (q_1, \wt{q}_1, q_2, \wt{q}_2), [w_0,w_1] \big)
& \mapsto & [w_0, w_1, w_0 q_1 - w_1 \wt{q}_2, w_0 q_2 + w_1
\wt{q}_1] \;. \end{array}
\label{cxha:fund} \end{equation}
The restriction of this to $\RR^4 \times \CP^1$ gives a
trivialisation of the twistor bundle (\ref{cxha:CP3S4}) over
$\RR^4$, cf.\ Section {\bf \ref{twistor:sec:twistorR4}}.

In summary, \emph{{\rm (i)} the map\/ $w \mapsto \wt{w}$ defined by\/
{\rm (\ref{cxha:wt})} determines a
bijection from\/ $\CP^3$ to the set of all affine\/ $\alpha$-planes in\/ $\Cf$;
if\/ $w \in \CP^3 \setminus \CP^1_0$, \ $\wt{w}$ is given on\/
$\CC^4$ by\/ {\rm (\ref{cxha:I})}; {\rm (ii)} a point\/ $(w,p) \in {F}_{1,2}$
represents the\/ $\alpha$-plane\/ $\Pi_p$ at\/ $p$ tangent to\/ $\wt{w}$}.

\section{Group actions} \label{cxha:sec:group-actions}
In order to decide what distributions are conformally equivalent,
we need to discuss group actions on the diagram (\ref{cxha:doublefibr}).
The group $\SL(4,\CC)$ acts in a canonical way on $\CC^4$ and
induces transitive actions on the three spaces in (\ref{cxha:doublefibr})
such that $\mu$ and $\nu$
are equivariant.  {}From our description of the holomorphic conformal
structure on $\Cf = G_2(\CC^4)$ it is clear that $\SL(4,\CC)$ acts
on that space by holomorphic conformal transformations.

The subgroup
$\bigl\{B \in \GL(6,\CC): \det B = 1,\, g^{\cC}\bigl(B(\vec{\xi}),B(\vec{\xi})\bigr) =
g^{\cC}(\vec{\xi},\vec{\xi})\bigr\}$ is isomorphic to the
complex orthogonal group $\SO(6,\CC)$.  This acts on $\CP^5$ preserving
the complex quadric $\Qq^C$ and restricts to an action on
$\Qq^C$ by holomorphic conformal tranformations.

As before, identify $\La^2\CC^4$ with $\CC^6$ in the standard way. Then,
given $P \in \SL(4,\CC)$, set $R(v \wedge w) = Pv \wedge Pw$; this defines
a  homomorphism
\begin{equation} \label{cxha:SL4CSO6C}
\SL(4,\CC) \to \SO(6,\CC)
\end{equation}
which is easily seen to be a double covering.

Write
\begin{equation} \label{cxha:PABCD}
P = \left[ \begin{array}{cc} A & B \\ C & D \end{array} \right]
\end{equation}
where $A,B,C,D$ are $2 \times 2$ complex matrices.  We discuss which
subgroups give transformations of $\Rf$ and $\Mf$.

For $\Rf$, let $\SL(2,\HH)$ be the subgroup of matrices $P$ of
$\SL(4,\CC)$ which commute
with $x \mapsto \jj x$; i.e., with $A,B,C,D$ quaternionic.  Then any $P$
sends the fibres of the twistor map (\ref{cxha:CP3S4}) into other
fibres and so
restricts to a transformation of $\Rf$.  Set $F_{1,2}^R = \{(w,p)
\in {F}_{1,2} : w \in \CP^3, p \in \Rf \}$; then $\mu$ restricts to a
diffeomorphism on $F_{1,2}^R$, and so the double fibration collapses to the twistor map
(\ref{cxha:CP3S4}).

For $\Mf$, let $h$ be the (pseudo-)Hermitian form on $\CC^4$ given by
$$
h(v,w) = v_0 \ov{w}_2 + v_1 \ov{w}_3
	+ v_2 \ov{w}_0 + v_3 \ov{w}_1
$$
so that $N^5$ is given by $h(w,w) = 0$.  Let
$$
\SU(4,h) = \{ P \in \SL(4,\CC) : h(Pv,Pw) = h(v,w) \}\,.
$$
Then $P \in \SU(4,h)$ if and only if
$$
\left[ \begin{array}{cc} A^* & C^* \\ B^* & D^* \end{array} \right]
	\left[ \begin{array}{cc} 0\ & I \\ I\ & 0 \end{array} \right]
	\left[ \begin{array}{cc} A & B \\ C & D \end{array} \right]
	= \left[ \begin{array}{cc} 0\ & I \\ I\ & 0 \end{array} \right] \,,
$$
i.e.,
\begin{equation} \label{cxha:SU4h-condn}
A^*C = -C^*A, \ A^*D+C^*B = I, \ B^*D = - D^*B \,.
\end{equation}
An easy linear
change of coordinates shows that $h$ has signature $(2,2)$ so that
$\SU(4,h)$ is
isomorphic to the group $\SU(2,2)$.  The group $\SU(4,h)$ partitions $\CP^3$ into
three disjoint orbits, namely,\linebreak[4]
$\CP^3_+ = \{z = [z_0,z_1,z_2,z_3] \in \CP^3 : h(z,z) > 0\}$, 
$\CP^3_- = \{z = [z_0,z_1,z_2,z_3] \in \CP^3 : h(z,z) < 0\}$ and
$N^5$.   The diagram (\ref{cxha:doublefibr}) restricts to
\begin{equation}
\begin{array}{ccccc} 
       &              & (w,p) \in \FF_{1,2} &              &   \\
       & \mu \swarrow &                   & \searrow \nu &   \\
w \in N^5 &        &                &                 & p \in \Mf
\end{array}
\label{cxha:doublefibr-M} \end{equation}
where $\FF_{1,2} = \{(w,p) \in {F}_{1,2} : w \in N^5, p \in \Mf \}$.

Since the quadric hypersurface $\Qq^M$ is given by the quadratic form
\eqref{cxha:QM} of signature $(2,4)$, the group
$\{ B \in \GL(6, \CC) : g(B\wt{\vec{\xi}}, B\wt{\vec{\xi}})
= g(\wt{\vec{\xi}},\wt{\vec{\xi}})\}$ is isomorphic to $\Orthog(2,4)$.
It clearly acts by conformal transformations on $\Qq^M$; in fact we have
a double cover
\begin{equation} \label{cxha:O24C13} 
\Orthog(2,4) \to \C(1,3)
\end{equation}
of the full conformal group of $\Qq^M$.

The homomorphism (\ref{cxha:SL4CSO6C}) restricts to a double cover
\begin{equation} \label{cxha:homSU22}
\SU(2,2) \to \Orthog_0(2,4) \
\end{equation}
of the identity component of $\Orthog(2,4)$.
The composition of this with (\ref{cxha:O24C13}) gives a $4:1$ covering
\begin{equation} \label{cxha:SU22C13}
\SU(2,2) \to \C_0(1,3) = \C^{\uparrow}_+(1,3)
\end{equation}
of the identity component of the conformal group of $\Qq^M$ which
consists of those conformal transformations which preserve
time- and space-orientations.

\begin{remark} \label{cxha:rem:groups}  \rm 
{\rm (i)} The double cover {\rm (\ref{cxha:SL4CSO6C})} exhibits
$\SL(4,\CC)$ as the \emph{universal covering} of\linebreak $\SO(6,\CC)$,
i.e., $\SL(4,\CC) \cong \Spin(6,\CC)$.

{\rm (ii)}  Since $\MM^4$ is conformally embedded in $\Mf$, any conformal
transformation $R \in C(1,3)$ of $\Mf$ induces \emph{local} conformal
transformations of $\MM^4$, or a \emph{global} conformal transformation if $R$ maps the points at infinity to
points at infinity.
\end{remark}

\begin{example} \label{cxha:ex:SL4C}
Given $P \in \SL(4, \CC)$ as in (\ref{cxha:PABCD}), in the standard chart
(\ref{cxha:emb}), the action on $\CC^4$ is
given by
\begin{equation} \label{cxha:Mobius}
Q \mapsto (C+DQ)(A+BQ)^{-1}
\quad \mbox{where} \quad
Q = \left( \begin{array}{rr} 	q_1 & -\wt{q}_2 \\
				q_2 & \wt{q}_1 \end{array} \right) \,.
\end{equation}
If $Q \in \SL(2,\HH)$ then this restricts to a conformal transformation
of the space $\Rf$; if $Q \in \SU(4,h)$ it restricts to a conformal
transformation of $\Mf$.
We now discuss some special cases.

(i) If $B = C = 0$, the transformation \eqref{cxha:Mobius} is $Q \mapsto DQA^{-1}$.  If $D$ and $A$ are in
$\SU(2) = \Symp(1)$, i.e., represent unit quaternions, then $P \in \SL(2,\HH)$
and (\ref{cxha:Mobius}) restricts to a
rotation of $\RR^4$.  If, on the other hand, $A = (D^*)^{-1}$ with $D \in
\SL(2, \CC)$, then $P$ lies in $\SU(4,h) \cong \SU(2,2)$ and
(\ref{cxha:Mobius}) restricts to the \emph{Lorentz transformation}
$Q \mapsto DQD^*$ of $\MM^4$;  in fact this is a `restricted' Lorentz
transformation, i.e., it lies in the identity component
$\Orthog_0(1,3) = \Orthog^{\uparrow}_+(1,3)$
of $\Orthog(1,3)$ consisting of those Lorentz transformations which
preserve time- and space-orientations, and so we obtain a double covering
of $\Orthog_0(1,3)$ by $\SL(2,\CC)$, important in spinor theory.

(ii) If $A = \la^{-1/2}I, \ B = \la^{1/2}I$ are multiples of the identity with $\la$ real and
positive, then $P \in \SL(2,\HH) \cap \SU(2,2)$, and on $\CC^4$, $\RR^4$
or $\MM^4$, the transformation
(\ref{cxha:Mobius}) is the \emph{dilation} $Q \mapsto \la Q$.

(iii) If $A = D = I$ and $B=0$, then (\ref{cxha:Mobius}) is the
\emph{translation}
$Q \mapsto Q+C$.  If $C$ is of the form (\ref{cxha:Qreal}), then
$P \in \SL(2,\HH)$ and the translation restricts to a translation
through $\vec{x}$ in
$\RR^4$; if $C$ is of the form (\ref{cxha:QMin}), i.e., is skew-Hermitian,
then $P \in \SU(4,h)$ and it restricts to a translation through $\vec{x}$
in $\MM^4$.

(iv) If $A = D = 0$ and $B = C = I$, then (\ref{cxha:Mobius}) is the
\emph{inversion} $Q \mapsto Q^{-1}$, i.e.
$$
\vec{x} = (x_0,x_1,x_2,x_3) \mapsto \frac{1}{g^{\cC}(\vec{x},\vec{x})}(x_0,-x_1,-x_2,-x_3) \,.
$$
Since $P \in \SL(2,\HH) \cap
\SU(2,2)$, this restricts to inversions in $\RR^4$ and $\MM^4$, the
latter given by
\begin{equation} \label{cxha:inversion}
\vec{x} = (t,x_1,x_2,x_3) \mapsto \frac{1}{\abs{\vec{x}}^2_1}(t,-x_1,-x_2,-x_3) \,.
\end{equation}
\end{example}

\begin{remark} \label{cxha:rem:inversion}  \rm 
If {\rm (\ref{cxha:SU4h-condn})} is replaced by
\begin{equation} \label{cxha:hto-h}
A^*C = -C^*A, \ A^*D+C^*B = -I, \ B^*D = - D^*B \,,
\end{equation}
then $P$ sends $h$ to $-h$ and so interchanges $\CP^3_{\pm}$. 
For example, if $A = C = 0$, $ B=-I$ and $C=I$, then $P$ satisfies
{\rm (\ref{cxha:hto-h})}.  The induced
action {\rm (\ref{cxha:Mobius})} on $\CC^4$ is reflection through the
origin: $\vec{x} \mapsto -\vec{x}$, this restricts to an orientation preserving
conformal transformation on $\MM^4$; however this is 
time- and space-orientation reversing.

Similarly, $A = C = 0, \ B=I, \ C=-I$ satisfies {\rm (\ref{cxha:hto-h})}
and induces the inversion $Q \mapsto -Q^{-1}$ which is, again,
time- and space-orientation reversing.

\end{remark}

\begin{example} \label{cxha:O6C}
Take the coordinates for $\RR^6 = \La^2 \RR^4$ given by
$$
(\vec{\eta}, \wh{\vec{\xi}}) = (\eta_0, \eta_1,
\wt{\xi}_2, \wt{\xi}_3, \wt{\xi}_4, \wt{\xi}_5)
$$
where $\eta_0 =
\wt{\xi}_0 + \wt{\xi}_1$, \ $\eta_1 = \wt{\xi}_0 - \wt{\xi}_1$ and
$(\wt{\xi}_0, \ldots, \wt{\xi}_5)$ are the coordinates given in
(\ref{cxha:wtxi}).  Thus $\Qq^{\cC}$ has equation $\eta_0 \eta_1 -
\abs{\wh{\vec{\xi}}}^2_1 = 0$ and $\Qq^M$ is the set of points of $\Qq^{\cC}$ with
real coordinates.  Using these coordinates, 
given $R \in \SO(6,\CC)$, write it in block form as $R =
	\left( \begin{array}{cc} 	E & F \\
					G & H \end{array} \right)$
where $E$ is $2 \times 2$.
Since the embedding of $\CC^4$ in $\Qq^{\cC}$ is given by
$$
{\vec{x}} \mapsto (\eta_1, \eta_2, {\vec{x}})
	= 2\bigl(1,\,g^{\cC}({\vec{x}},{\vec{x}}),\,{\vec{x}}\bigr),
$$
the action of $R$ on $\CC^4 \subset \Qq^{\cC}$ is easily calculated to be
\begin{equation} \label{cxha:actionR}
{\vec{x}} \mapsto \frac{1}{e_{11} + e_{12} \,g^{\cC}({\vec{x}},{\vec{x}}) + f_{1\cdot} \cdot \wh{\vec{\xi}}}
	\big( g_{\cdot 1} + g^{\cC}({\vec{x}},{\vec{x}}) g_{\cdot 2} + H\vec{x} \big)
\end{equation}
where $E = (e_{ij})$, $f_{1\cdot}$ (respectively,  $g_{\cdot j}$) denotes the
first row of $F$ (respectively,  $j$'th column of $G$) and $f_{1\cdot} \cdot \wh{\vec{\xi}} = 
\sum_{j=1}^4f_{1j}\wh{\xi}_{j+1}$.

If all the entries of $R$ are real so that $R \in \SO(2,4)$, then this
action restricts to a conformal transformation of $\Mf$.
We now discuss some particular cases.

(i)  Set $E = I, F=G=0$; then $R \in \Orthog(2,4)$ if $H \in \Orthog(1,3)$.
In this case the action (\ref{cxha:actionR}) is the
\emph{Lorentz transformation} $x \to Hx$, this is a \emph{restricted}
Lorentz transformation if $H$ is in the identity component of
$\Orthog(1,3)$, in which case $R \in \Orthog_0(2,4)$.

(ii)  Set $E = \left( \begin{array}{ll}	\la^{-1} & 0 \\
					0        & \la \end{array}\right)$
where $\la$ is a non-zero real number, and $F = G = 0$, $H = I$.	
Then $K \in \Orthog_0(2,4)$ and the action is the \emph{dilation}
${\vec{x}} \mapsto \la {\vec{x}}$.

(iii)  For ${\vec{a}} \in \CC^4$ regarded as a column vector, set
$E = \left( \begin{array}{cc}	1       & 0 \\
				\abs{{\vec{a}}}^2_1 & 1 \end{array}\right)$,
$F = \left( \begin{array}{cc}	{\vec{0}} \\ 2{\vec{a}}^{\tT} \end{array}\right)$,
$G = \left( \begin{array}{cc}   {\vec{a}} & {\vec{0}} \end{array}\right)$
and $H = {\vec{0}}$.  Then $R \in \Orthog_0(2,4)$ and the action is the
\emph{translaton} ${\vec{x}} \mapsto {\vec{x}} + 2{\vec{a}}$.

(iv)  Set $E = \left( \begin{array}{cc}	0 & 1 \\
					1 & 0 \end{array} \right)$,
and $F = G = {\vec{0}}$, then $R \in \Orthog(2,4)$ if and only if $H \in \Orthog(1,3)$,
and the action on
$\MM^4$ is given by ${\vec{x}} \mapsto H{\vec{x}}/ \norm{\vec{x}}^2_1$. For example, if $H
= \mbox{\rm diag}\{1,-1,-1,-1\}$, then $R$ lies in $\Orthog_0(2,4)$ and
the transformation is the \emph{inversion} (\ref{cxha:inversion}).

On the other hand if $H = -I$, then $R$ is in $\Orthog(2,4)$ but not in the
identity component; the corresponding transformation ${\vec{x}} \mapsto
-{\vec{x}}/\norm{\vec{x}}^2_1$ is in $\C(1,3)$ but not in $\C^{\uparrow}_+(1,3)$.

Note that all the examples with matrices $R \in
\Orthog_0(2,4)$ correspond to matrices
$P$ in $\SU(2,2)$ under the double cover (\ref{cxha:SU22C13}), the form of
such $P$ for (i)--(iv) above being as in the corresponding number in
Example \ref{cxha:ex:SL4C}.

If we use the coordinates \eqref{cxha:xi} instead,
we can give similar examples of matrices
$R \in \SO(1,5)$ which give conformal transformations of $\Rf$. 
\end{example}

\section{A Kerr theorem} \label{cxha:sec:Kerr}

Let $J$ be a Hermitian structure on an open subset $A^4$ of $\RR^4$.
Then, as in Section \ref{cxha:sec:Herm}, $J$ is represented by a map
$J:A^4 \to \CP^1$; combining this with the fundamental map
(\ref{cxha:fund}) gives a map
\begin{equation} \label{cxha:compo}
A^4 \stackrel{(I,J)}{\lra} A^4 \times \CP^1 \stackrel{\iota}{\to}
	\CP^3 \setminus \CP^1_0 \subset \CP^3,
\end{equation}
where $I$ denotes the identity map.
Comparing definitions we see that this is just the section $\si_J$
of the twistor bundle (\ref{cxha:CP3S4}) defined by $J$, as in Section 
{\bf\ref{twistor:sec:twistorspace}}.  As proved there, $\si_J$ is a
holomorphic map $(A^4,J) \to (\CP^3, \mbox{standard})$; this can also be
seen by noting that the two maps in (\ref{cxha:compo}) are both
holomorphic.  Thus, as in Section {\bf\ref{twistor:sec:twistorspace}},
$\si_J(A^4)$ is a complex hypersurface of $\CP^3$,
In the case of $\RR^4$, or its compactification $\Rf$,
Proposition {\bf\ref{twistor:prop:JsigmaJ}}(iii) gives the following Kerr-type
result.

\begin{proposition} \label{cxha:prop:Kerr}
Given a complex hypersurface\/ $S$ of\/ $\CP^3$, any
smooth map\/ $w:A^4 \to \CP^3$ from an open subset of\/ $\Rf$
with image in\/ $S$ defines a Hermitian structure\/ $J$ on\/ $A^4$.

All Hermitian structures\/ $J$ on open subsets\/ $A^4$ of\/ $\Rf$ are given
this way by setting\/ $S = \si_J(A^4)$.
\qed \end{proposition}

On using Theorem \ref{cxha:th:unifgerms} we obtain similar results for the
other three sorts of distribution in (\ref{cxha:wh}), but now we must use
the mapping \ $\wh{}$ \ defined by (\ref{cxha:doublefibr}).  For completeness,
we include the last theorem with this notation, as follows.

\begin{corollary} \label{cxha:cor:Kerr1}
Given a complex hypersurface\/ $S$ of\/ $\CP^3$, any smooth map\/ $w:A \to
\CP^3$ from an open subset of\/ {\rm (i)} $\Cf$ {\rm (}respectively, {\rm (ii)} $\Rf$,
{\rm (iii)} $\Mf$, {\rm (iv)} $\Rt${\rm )} with
\begin{equation} \label{cxha:inc}
w({\vec{p}}) \in \wh{{\vec{p}}} \qquad ({\vec{p}} \in A)
\end{equation}
 defines a
{\rm (i)} holomorphic foliation\/ $\Pi$ by\/ $\alpha$-planes {\rm (}respectively,
{\rm (ii)} a positive Hermitian structure\/ $\JJ$,
{\rm (iii)} a real-analytic shear-free ray congruence\/ $\ell$,
{\rm (iv)} a real-analytic conformal foliation\/ $\Cc$ by curves{\rm )}
on\/ $A$.

All such distributions are given this way.
\end{corollary}

We call $S$ the \emph{twistor surface of\/ $\Pi, \JJ, \ell$ or $\Cc$}.

Note that in case (ii) (respectively, (iv)) the `incidence condition'
(\ref{cxha:inc}) reduces
to the condition that $w$ be a section of the twistor bundle
(\ref{cxha:CP3S4}) (respectively, its restriction (\ref{cxha:N5S3})).

We can give a more explicit version for uncompactified spaces,
as follows.

\begin{corollary} \label{cxha:cor:Kerr2}
Let\/ $\psi$ be a homogeneneous complex-analytic function of four complex
variables. 
Let\/ $[w_0, w_1] = \mu(q_1, \wt{q}_1, q_2, \wt{q}_2)$ be
a smooth solution to the equation
\begin{equation}
\psi (w_0,\, w_1,\, w_0 q_1 - w_1 \wt{q}_2,\, w_0 q_2 + w_1
\wt{q}_1) = 0 \;.
\label{cxha:K} \end{equation}
on an open subset\/ $A$ of\/ {\rm (i)} $\CC^4$, {\rm (}respectively,
{\rm (ii)} $\RR^4$, {\rm (iii)} $\MM^4$, {\rm (iv)} $\RR^3${\rm )}.
Then\/ $\mu = w_1/w_0$ represents\/
{\rm (i)} a complex-analytic foliation by null planes {\rm (}respectively,
{\rm (ii)} a positive Hermitian structure,
{\rm (iii)} a real-analytic shear-free ray congruence,
{\rm (iv)} a real-analytic conformal foliation by curves{\rm )}
on\/ $A$.

Further, the twistor surface of the distribution {\rm (i)--(iv)} is given by
\begin{equation}
\psi(w_0, w_1, w_2, w_3) = 0 \,.
\label{cxha:S} \end{equation}

All such distributions are given this way locally.
\qed \end{corollary}

\begin{example} \label{cxha:ex:HopfKerr}
Set $\psi(w_0, w_1, w_2, w_3) = w_3$.  Then (\ref{cxha:K}) reads
$$
q_2 + \mu\wt{q}_1 = 0 \,;
$$
this has solution $\mu:\CC^4 \setminus \{\vec{0}\} \to \Ci$ given by
(\ref{cxha:Hopfcx}), this defines the distributions $\Pi, \JJ, \ell, \Cc_0$
of Example \ref{cxha:ex:Hopf}.
\end{example}

For further examples, see Section \ref{cxha:sec:examples}.

\section{CR interpretation} \label{cxha:sec:CR}

A \emph{CR structure} of class $C^r$ (Jacobowicz 1990)
on an odd-dimensional manifold $M =
M^{2k+1}$ is a choice of $C^r$ rank $k$ complex subbundle $\Hh$ of the
complexified tangent bundle $T^{\cC}M = TM \otimes \CC$ with $\Hh \cap
\overline{\Hh} = \{\vec{0}\}$ and
\begin{equation}
\bigl[C^{\infty}(\Hh), C^{\infty}(\Hh)\bigr] \subset C^{\infty}(\Hh)\,.
\label{cxha:CRint}
\end{equation}
Given $\Hh$ we define the \emph{Levi subbundle} $H$ of $TM$ by $H \otimes
\CC = \Hh \oplus \ov{\Hh}$ and a complex-linear endomorphism
$J:H \otimes
\CC \ra H \otimes \CC$ by multiplication by $-\ii$ (respectively $+\ii$)
on $\Hh$ (respectively $\ov{\Hh}$).  Then $J$ restricts to an endomorphism
of $H$ with $J^2 = -I$, and (\ref{cxha:CRint}) is equivalent to
the condition:
\begin{equation} \label{cxha:CRint2}
\left. \begin{array}{l}
\mbox{(i) if } X,Y \mbox{ are sections of } H \mbox{ then so is }
[J X,Y]+[X,J Y]\,; \\
\hspace{-0.3em}\mbox{(ii) } J\bigl([J X,Y]+[X,J Y]\bigr) = [J X,J Y] - [X,Y] \,. \end{array} \right\}
\end{equation}
Conversely, such a pair $(H,J)$ determines $\Hh$ by $\Hh = \{ X + \ii J X : X
\in H \}$;
thus a CR structure can be defined equivalently as a subbundle $H$ of
$TM$ of rank $k$ together with an endomorphism $J:H \to H$ with $J^2 = -I$
which satisfies (\ref{cxha:CRint2}).

Say that a map $\phi:M \to N$ is \emph{CR with respect to CR
structures\/ $(H,J)$ on $M$ and\/ $(H',J')$ on\/ $N$} if $\dd\phi$ maps $H$
to $H'$ and $\dd\phi \circ J = J' \circ \dd\phi$ on $H$, equivalently
$\dd\phi:T^{\cC}M \to T^{\cC}N$ maps $\Hh$ to $\Hh'$.   If $\phi$ is also
a diffeomorphism, it is called a \emph{CR diffeomorphism}, clearly its
inverse is automatically CR.

Replacing $H$ by $TN$ and $\Hh$ by $T^{0,1}N$, we obtain the definition of a
\emph{CR map from a CR manifold to a complex manifold}.

\begin{example} \label{cxha:ex:fCR}
Let $A^3$ be an open subset of $\RR^3$, or of any Riemannian $3$-ma\-ni\-fold, and 
let $\UU$ be the unit positive tangent to an oriented foliation $\Cc$ by
curves of $A^3$. Give $A^3$ the CR structure defined by
$(H,J) = (\UU^{\perp},\JJ^{\perp})$.

Suppose that the leaves of the foliation $\Cc$ are given by the level sets of a 
submersion $f = f_1 + \ii f_2:A^3 \to \CC$. By replacing $f$ by its
conjugate $f_1 - \ii f_2$, if necessary, we
can assume that $\{\grad f_1, \grad f_2, U \}$ is positively oriented.
 Then \emph{$f$ is horizontally conformal if and
only if it is CR with respect to the structure $(H,J)$, and this holds if and only if\/
$\Cc$ is a conformal foliation}.
\end{example}

Any real hypersurface $M^{2k+1}$ of a complex manifold $(\wt{M},J)$
has a canonical CR structure called the \emph{hypersurface CR
structure} given at each point $p \in M$ by
(i) $\Hh_p = T^{0,1}_p\wt{M} \cap T^{\cC}_p M$. 
Equivalently, (ii) $H_p = T_pM \cap JT_pM$.  Note that this necessarily has
dimension $2k$ at all points $p$.

More generally, for any submanifold $M^{2k+1}$ of a complex manifold $(\wt{M},J)$,
we may define $\Hh_p$ and $H_p$ by (i) and (ii); if the dimension of
$H_p$ is $2k$ at all points $p \in M^{2k+1}$, it defines a CR structure
on $M^{2k+1}$ and $M^{2k+1}$ is called a \emph{CR submanifold}
of $(\wt{M},J)$.  The inclusion map is then a CR map.

We give the unit tangent bundle $T^1 \RR^3 = \RR^3 \times S^2$ a
CR structure $(H,J)$, as follows.  At each $({\vec{p}},\UU) \in \RR^3
\times S^2$, use the canonical isomorphism $\RR^3 \cong T_{\vec{p}}\RR^3$
to regard $\UU^{\perp}$ as a subspace $\UU_{\vec{p}}^{\perp}$ of
$T_{\vec{p}}\RR^3$\,; then set $H_{({\vec{p}},\UU)} = \UU_{\vec{p}}^{\perp} \oplus T_{\UU}
S^2$ and define $J$ to be rotation through $+\pi/2$ on $\UU_{\vec{p}}^{\perp}$
together with the standard complex structure $J^{S^2}$ on $T_{\UU}
S^2$.  Define a complex structure $\Jj$ on $\RR^4 \times S^2$ by
\begin{equation} \label{cxha:Jfund}
\Jj: \RR^4 \times S^2 \lra
\mbox{End}\,T(\RR^4 \times S^2)\,, \qquad ({\vec{p}},\UU) \mapsto (\JJ(\UU)_{\vec{p}},
J^{S^2})
\end{equation}
with $\JJ(\UU)_{\vec{p}}$ the unique positive almost Hermitian structure which is determined by
$\JJ(\UU)_{\vec{p}}(\pa/\pa x_0) = \UU_{\vec{p}}$\,, cf.\ (\ref{cxha:unif}). A calculation,
for example using Proposition
{\bf\ref{twistor:prop:JsigmaJ}}(ii), shows that $\Jj$ is integrable.
Then the CR structure $(H,J)$ on $T^1\RR^3$ is the hypersurface CR structure given by regarding
$\RR^3 \times S^2$ as a real hypersurface of the manifold ($\RR^4
\times S^2, \Jj)$.

More generally, for any oriented Riemannian $3$-manifold $M^3$,
we can give the unit tangent bundle $T^1 M^3$ a CR structure
$(H,J)$, as follows.  At each
point $(p,\UU) \in T^1 M^3$ ($p \in M^3, \ \UU \in T^1_p M^3$),
the Levi-Civita connection on $M^3$ defines a splitting
$T_{(p,\UU)}(T^1 M^3) = T_p M^3 \oplus T_{\UU}(T^1_p M^3)$. Since
$T^1_p M^3$ is canonically oriented and is isometric to a $2$-sphere, it has a canonical
K{\"a}hler structure $J^{S^2}$.  We set $H_{(p,\UU)} =
\UU_p^{\perp} \oplus T_{\UU} (T^1_p M^3)$ and define $J$ to be rotation through
$+\pi/2$ on $\UU_p^{\perp}$ together with $J^{S^2}$ on $T_{\UU} (T^1_p
M^3)$.  It can be checked that (\ref{cxha:CRint2}) is satisfied.

For $M^3 = S^3$ this can be described more explicitly.   The
differential of the canonical embedding $S^3 \hookrightarrow
\RR^4$ defines an embedding $i:T^1 S^3 \hookrightarrow T\RR^4 =
\RR^4 \times \RR^4$.  At a point $({\vec{p}},\UU) \in T^1 S^3 \ ({\vec{p}} \in
S^3, \ \UU \in T^1_{\vec{p}} S^3)$, we have
\begin{equation*}
\dd i(T_{({\vec{p}},\UU)} T^1 S^3) = \{ (X,u) \in {\vec{p}}^{\perp} \times \UU^{\perp} : \inn{
X, \UU} + \inn{{\vec{p}}, u} = 0 \} 
	  \subset \RR^4 \times \RR^4 \cong T_{({\vec{p}},\UU)}(\RR^4 \times \RR^4) \;;
\end{equation*}
then we choose $H_{({\vec{p}},\UU)} = (\UU \oplus {\vec{p}})^{\perp} \times (\UU \oplus
{\vec{p}})^{\perp}$ and $J =$ rotation through $+\pi/2$ on each oriented plane $(\UU
\oplus {\vec{p}})^{\perp}$.  That this satisfies (\ref{cxha:CRint2})
can be seen either by direct calculation or by
noting that any stereographic projection $S^3 \setminus
\{\mbox{point} \} \ra \RR^3$ is conformal and induces a CR
diffeomorphism between $T^1(S^3 \setminus \{\mbox{point} \}$) and
$T^1 \RR^3$.

Let $A^3$ be an open subset of $\RR^3$ or $S^3$, and 
let $\UU$ be the unit positive tangent to an oriented foliation by
curves of $A^3$.  We can regard $\UU$ as a smooth map from $A^3$ to $T^1\RR^3$ or $T^1S^3$.
Give $A^3$ the CR structure defined by $(H,J) = (\UU^{\perp},\JJ^{\perp})$.
Then we can interpret (\ref{cxha:conffoln2}) as a CR condition, as follows.

\begin{proposition} \label{cxha:prop:conffoln-CR}
An oriented foliation on\/ $A^3$ is conformal if and only if its positive
unit tangent vector field\/ $\UU: A^3 \ra T^1 \RR^3$ or\/ $T^1 S^3$ is CR.
\qed \end{proposition}

Define a map $k:T^1 \Rt \ra N^5$ by sending a unit tangent
vector $\UU$ at a point ${\vec{p}}$ of $\Rt \cong S^3$ to the point of $\CP^3$
representing the affine null geodesic through ${\vec{p}}$ tangent to
$\pa/\pa t + \UU$.  (This is well defined by Remark
\ref{cxha:rem:pa/pa-t}.)
Since every null geodesic of $\Mf$ meets $\Rt$
in precisely one point, $k$ is a diffeomorphism.  Its restriction to $\RR^3$
is given by the restriction
\begin{equation} \label{cxha:T1R3N5}
T^1\RR^3 \cong \RR^3 \times S^2 \ra N^5
\end{equation}
of the fundamental map (\ref{cxha:fund}). 
Since the latter map is holomorphic with respect to the complex
structure $\Jj$ of (\ref{cxha:Jfund}), \emph{the map
{\rm (\ref{cxha:T1R3N5})} is a CR diffeomorphism}.

Now, given a congruence of null lines on an open subset $A^M$ we have
a map $w:A^M \to N^5$ which sends ${\vec{p}} \in A^M$ to the point of $N^5$
representing (see Section \ref{cxha:sec:twistor}) the null line of the
congruence through ${\vec{p}}$. 
{}From Theorem \ref{cxha:th:SFRconf} we obtain a characterisation of \emph{shear-free}, as follows.

\begin{proposition} \label{cxha:prop:SFR-CR}
Let\/ $\ell$ be a\/ $C^{\infty}$ congruence of null lines on an open
subset\/ $A^M$ of\/ $\MM^4$.  Then\/ $\ell$ is
shear-free if and only if the map\/ $w:A^M \ra N^5$ representing it
is CR when restricted to\/ $A^3$.
\qed \end{proposition}

\begin{corollary} \label{cxha:cor:SFR-CR}
Let\/ $N^3$ be a smooth CR submanifold of\/ $\CP^3$ contained in\/ $N^5$.
Then any smooth section\/
$w:A^3 \to N^5$ of\/ {\rm (\ref{cxha:N5S3})} defined on an open set of\/ $\RR^3$
with image in\/ $N^3$ defines a conformal foliation by curves\/ $\Cc$ of\/ $A^3$.
\end{corollary}

\begin{proof}
It is easy to check that such a smooth section $w$ of $N^5 \cong T^1 A^3$
is CR if and only if its image is a CR submanifold of $\CP^3$.
\end{proof}

Thus, given a conformal foliation $\Cc$ by curves, or the associated
SFR congruence $\ell$, its image under $w$ is a CR
submanifold $N^3$ of $N^5$.  If $\ell$ is real analytic then $N^3$ is
the intersection of a complex hypersurface $S$ with $N^5$; the hypersurface $S$ is the
twistor surface of $\Cc$ and $\ell$, as in Corollary \ref{cxha:cor:Kerr1}.

If $\Cc$ (or, equivalently, $\ell$) is only $C^{\infty}$, then it cannot necessarily
be extended to an open subset of $\CC^4$, but it can be extended to one side, as follows.

\begin{proposition} \label{cxha:prop:Cinfty1}
Let\/ $\Cc$ be a\/ $C^{\infty}$ conformal foliation by curves of an open subset $A^3$ of
$\RR^3$, with horizontal distribution everywhere non-integrable,
and let\/ $\UU$ denote its unit positive tangent.
Then the vector field\/ $\UU$ can be
extended to a positive Hermitian structure\/ $J$ on one side of\/ $\RR^3$;
precisely, set\/ $\RR^4_+=\{ (x_0,x_1,x_2,x_3)
\in \RR^4 : x_0 > 0 \}$ and\/ $\RR^4_-=\{ (x_0,x_1,x_2,x_3)
\in \RR^4 : x_0 < 0 \}$, then there exists a positive\/
$C^{\infty}$ Hermitian structure\/ $J$ on an open subset\/ $A^4$ of\/
$\RR^4_+ \cup \RR^3$ or of\/ $\RR^4_- \cup \RR^3$ with\/ $A^4 \cap \RR^3
= A^3$ and\/ $\JJ(\pa/\pa x_0) = \UU$ on $A^3$.
\end{proposition}

\begin{proof}
The \emph{Levi form} of a CR manifold $(M,\Hh)$ is the mapping
$\Hh \times \ov{\Hh} \to T^{\cC}M/(\Hh \oplus \ov{\Hh})$ induced by
$(Z,\ov{W}) \mapsto -\half\ii\,[Z,\ov{W}]$.  If $\Cc$
has horizontal distribution everywhere non-integrable, then clearly the Levi form
of the associated CR stucture on $A^3$ is non-zero, as is that of the isomorphic 
CR manifold $N^3$.  Then,  by a theorem of Harvey and
Lawson (1975, Theorem 10.2), $N^3$ is the boundary of a complex hypersurface
$S$.  This means that $J$ can be defined on one side of $\RR^3$.
\end{proof}

\begin{remark}  \rm  The non-integrability of the horizontal spaces means that the
twist ($g([e_1,e_2],\ww) = g([e_1,e_2],\UU)$ (see Remark \ref{cxha:rem:twist}) is
non-zero, and $J$ can be extended to the side $x_0 > 0$ (respectively, $< 0$)
according as $g([e_1,e_2],\UU) > 0$ (respectively, $<0$) for a
positive orthonormal frame $\{e_1,e_2,\UU\}$.
\end{remark}

\section{The boundary of a hyperbolic harmonic morphism}
\label{cxha:sec:hypha}

We shall now discuss how horizontally conformal submersions are the boundary values
at infinity of harmonic morphisms from hyperbolic $4$-space.
\begin{proposition} \label{cxha:prop:restr-R3}
Let\/ $\phi:A^4 \ra \CC$ be a submersive map from an open subset
of\/ $\RR^4$ which is holomorphic with respect to a positive
Hermitian structure\/ $J$ on\/ $A^4$. Let\/ ${\vec{p}} \in A^4$ and set\/ $\UU =
\JJ(\pa/\pa x_0)$ on
the open subset\/ $A^3 = A^4 \cap \RR^3_{\vec{p}}$ of\/ $\RR^3_{\vec{p}}$.  Let\/ $\Cc$
denote the foliation of\/ $A^3$ given by the integral curves of\/ $\UU$.
Set $f= \phi\vert_{A^3}$, so that $f$ is real analytic.
Then
\begin{equation}
{\rm (i)} \quad \frac{\pa \phi}{\pa x_0} = 0 \mbox{ on } A^3
\label{cxha:orthog} \end{equation}
if and only if\/ {\rm (ii)} $f= \phi\vert_{A^3}$ is constant on the leaves of\/ $\Cc$.

Further, in this case, {\rm (iii)} $f$ is a horizontally conformal submersion.
\end{proposition}

\begin{proof}
Let ${\vec{q}} \in A^3$.  By the holomorphicity of $\phi$, since $\UU_{\vec{q}} =
\JJ_{\vec{q}} (\pa/\pa x_0)$, equation \eqref{cxha:orthog} is equivalent to the vanishing of
the directional derivative $\UU_{\vec{q}}(f)$, which is equivalent to constancy of
$f$ on the leaves of $\Cc$.  In this case, let $\{ e_2,\, e_3 = \JJ e_2 \}$ be a
basis for $\UU_{\vec{q}}^{\perp} \cap \RR_{\vec{p}}^3\,$; then holomorphicity of $\phi$
implies that $e_3(f) = \ii\, e_2(f)$ so that $f$ is horizontally
conformal.  Submersivity of $f$ easily follows from that of
$\phi$.
\end{proof}

\begin{remark}  \rm
{\rm (i)} In particular, $\Cc$ is a conformal foliation; this also follows from
Corollary {\rm \ref{cxha:cor:JtoU}}.

{\rm (ii)} Any two of the conditions (i), (ii), (iii) in the proposition imply the other.
\end{remark}

We now show how to find such functions $\phi$ using harmonic
morphisms. Equip $\breve{\RR}^4 \equiv \RR^4 \setminus \RR^3$
with the \emph{hyperbolic metric} $g^H = \bigl( \sum_{i=0}^{3}
\dd {x_i}^{\! 2} \bigr) \big/ {x_0}^{\! 2}$ so that each component $\RR_+^4$,
$\RR_-^4$ is isometric to hyperbolic $4$-space $\RH^4$. Let
$\breve{A}^4$ be an open subset of $\breve{\RR}^4$, then we call
a smooth map $\phi:\breve{A}^4 \ra \CC$  a \emph{hyperbolic
harmonic map} if it is a harmonic map with respect to the
hyperbolic metric $g^H$.  By calculating the tension field from
({\bf\ref{hamap:tensionlc}}) we see that $\phi$ \emph{is hyperbolic harmonic
if and only if}
\begin{equation}
x_0 \sum_{i=0}^3 \frac{\pa^2\phi}{\pa {x_i}^{\! 2}} -
2 \frac{\pa\phi}{\pa x_0} = 0
\label{cxha:hyp-ha} \end{equation}
\emph{at all points of} $\breve{A}^4$.

Similarly, $\pi:\breve{A}^4 \ra \CC$ will be called a
\emph{hyperbolic harmonic morphism} if it is a harmonic morphism with
respect to the metric $g^H$; by the fundamental result of Fuglede and Ishihara (see Theorem {\bf \ref{fund:th:char}}),
such maps are characterized as
satisfying (\ref{cxha:hyp-ha}) and the horizontal weak conformality
condition:

\begin{equation} \label{cxha:HWC}
\sum_{i=0}^3 \left(\frac{\pa\phi}{\pa x_i}\right)^{\!\!2} = 0 \,.
\end{equation}

Now, from Theorem {\bf\ref{twistor:th:Einstein4d}} we have the following.

\begin{corollary} \label{cxha:cor:Einsteinhyp}
{\rm (i)} Any submersive hyperbolic
harmonic morphism $\phi:\breve{A^4} \to \CC$ is holomorphic with
respect to some Hermitian structure\/ $J$ on\/ $\breve{A^4}$ and has
superminimal fibres with respect to\/ $J$, i.e., $\ker \dd\phi \subseteq
\ker \na^H \! J$ on\/ $\breve{A}^4$  where\/ $\na^H$ is the
Levi-Civita connection of the hyperbolic metric on\/ $\breve{\RR}^4$.

{\rm (ii)} Conversely, let\/ $J$ be a Hermitian structure on an open
subset\/ $\breve{A}^4$ of\/ $\breve{\RR^4}$, and\/ $\phi:\breve{A}^4
\to \CC$ a non-constant map which is holomorphic with respect to
$J$.  Then\/ $\phi$ is a hyperbolic harmonic map, and hence a hyperbolic
harmonic morphism, if and only if, at points
where\/ $\dd\phi$ is non-zero, its fibres are superminimal with respect to
$J$.
\qed \end{corollary}

To formulate this analytically, let $\Theta$ be the
holomorphic contact form on complex projective space $\CP^3$ given in homogeneous
coordinates by (i.e., its
lift to $\CC^4 \setminus \{\vec{0}\}$ is given by) 
\begin{equation}
w_1 \,\dd w_2 - w_2 \,\dd w_1 - w_0 \,\dd w_3 + w_3 \,\dd w_0\,.
\label{cxha:HCF} \end{equation}
Then $\ker\Theta$ gives the horizontal distribution of the
restriction   $\pi:\CP^3
\setminus N^5 \to (\breve{\RR}^4,g^H)$ of the twistor map (\ref{cxha:CP3S4}). Set $\Phi = \phi \circ
\pi$. Then $\phi$ has superminimal fibres if and only if
\begin{equation}
\ker \dd\Phi \subseteq \ker \Theta\,;
\label{cxha:hypsupermin} \end{equation}
this is thus the condition that $\phi$ be a hyperbolic harmonic morphism.
Equivalently, let $w:\breve{A}^4 \to \CP^3$ be the section of
$\pi$ corresponding to $J$, i.e., with image the twistor surface
$S$ of $J$ (see  Theorem \ref{cxha:prop:Kerr}); then we may pull back
$\Theta$ to a $1$-form $\theta = w^*\Theta$ on $\breve{A}^4$. 
Condition (\ref{cxha:hypsupermin}) now reads $\ker \dd\phi \subseteq \ker
\theta$,  this condition being equivalent to superminimality of the fibres of $\phi$.

Note that, if $A^4$ is an open subset of $\RR^4$, $\RR_+^4
\cup \RR^3$ or $\RR_-^4 \cup \RR^3$, with $A^4 \cup \RR^3$
non-empty, any ($C^2$, say) map which is a hyperbolic harmonic
map on $\breve{A}^4 = A^4 \setminus \RR^3$ satisfies
(\ref{cxha:hyp-ha}), and so (\ref{cxha:orthog}), at the boundary
$A^3 = A^4 \cup \RR^3$. A key property for us is the following converse.

\begin{proposition} \label{cxha:prop:orthog-hypha}
Let\/ $A^4$ be a connected open subset of\/ $\RR^4$, $\RR_+^4 \cup
\RR^3$ or\/ $\RR_-^4 \cup \RR^3$ such that\/ $A^3 = A^4 \cap \RR^3$ is
non-empty, and let\/ $\phi:A^4 \to \CC$ be a non-constant\/ $C^1$ map
which is holomorphic with respect to a Hermitian structure\/ $J$ on\/
$\breve{A}^4 = A^4 \setminus \RR^3$ and submersive at almost all
points of\/ $A^3$. Then\/ $\phi$ satisfies {\rm (\ref{cxha:orthog})}
on\/ $A^3$ if and only if\/ $\phi\vert_{\breve{A}^4}$ is a hyperbolic harmonic
morphism.
\end{proposition}

\begin{proof}
It suffices to work at points where $\phi$ is submersive. At such
points, note that (\ref{cxha:orthog}) holds if and only if $\ker \dd\phi =
\mbox{span} \{\pa/\pa x_0 , J \pa/\pa x_0 \}$. Now let $S$ be the
twistor surface of $J$, and let $\Phi:S \to \CC$ be defined by
$\Phi = \phi \circ \pi$.

We show that
the pull-back $\theta = w^* \Theta$ to $A^4$ satisfies
\begin{equation} \label{cxha:span-ker}
\mbox{\rm span} \{\pa/\pa x_0 , J \pa/\pa x_0 \} \subseteq \ker \theta
\end{equation}
at all points of $A^3$.  To do this,
since the (complexified) normal to $\RR^3$ is given by the
annihilator of $\mbox{span} \{\dd q_1 - \dd\wt{q}_1,\,
\dd q_2,\, \dd\wt{q}_2 \}$, it suffices to show that, on $\RR^3 =
\{ q_1 + \wt{q}_1 = 0 \}$, the form $\theta$ is a linear
combination of those three forms.  Now, on taking
differentials in (\ref{cxha:I}) we obtain
\begin{eqnarray*}
\dd w_2 & = & q_1 \,\dd w_0 + w_0 \,\dd q_1 - \wt{q}_2 \,\dd w_1 - w_1 \,\dd\wt{q}_2 \;,\\
\dd w_3 & = & q_2 \,\dd w_0 + w_0 \,\dd q_2 + \wt{q}_1 \,\dd w_1 + w_1 \,\dd\wt{q}_1 \;.
\end{eqnarray*}
Substituting these into (\ref{cxha:HCF}) and rearranging gives the expression
\begin{align*}
\theta =& \ (w_3 - w_0 q_2 + w_1 q_1)\,\dd w_0 - (w_2 + w_0 \wt{q}_1
+w_1 \wt{q}_2)\,\dd w_1\\ &+ w_0 w_1 (\dd q_1 - \dd\wt{q}_1) -
{w_0}^{\! 2} \dd q_2 - {w_1}^{\! 2} \dd\wt{q}_2 \;.
\end{align*}
But, by (\ref{cxha:I}), the coefficients of $\dd w_0$ and $\dd w_1$ vanish
when $q_1 + \wt{q}_1 = 0$, and (\ref{cxha:span-ker}) follows. 

Thus (\ref{cxha:orthog}) is equivalent to the condition 
\,$\ker \dd\phi \subseteq \ker \theta$\, of
superminimality of the fibres of $\phi$ at points of $A^3$
 or, equivalently \,$\ker \dd\Phi
\subseteq \ker \Theta$\, on the real hypersurface $N^3 = w(A^3)$ of
$S$.  But this is a complex-analytic condition, so by analytic
continuation, if $\phi$  has superminimal fibres at points of
$A^3$ then it has superminimal fibres on the whole of $A^4$.  By
Corollary \ref{cxha:cor:Einsteinhyp}, $\phi$ is a harmonic morphism.
\end{proof}

Combining Propositions \ref{cxha:prop:restr-R3} and
\ref{cxha:prop:orthog-hypha}, we obtain our main result, as follows.

\begin{theorem} \label{cxha:th:HWC-hyp}  
Let\/ $f:A^3 \to \CC$ be a real-analytic horizontally conformal
submersion on an open subset of\/ $\RR^3$. Then there is an open
subset\/ $A^4$ of\/ $\RR^4$ with\/ $A^4 \cap \RR^3 = A^3$, and a
real-analytic submersion\/ $\phi:A^4 \ra \CC$  with\/ $\phi\vert_{A^3} = f$,
such that\/ $\phi\vert_{{A^4} \setminus \RR^3}$ is a hyperbolic harmonic
morphism. In fact the restriction $\phi \mapsto f = \phi\vert_{A^3}$ defines a bijective
correspondence between germs at\/ $A^3$ of real-analytic
submersions\/ $\phi:A^4 \ra \CC$ on open neighbourhoods of\/ $A^3$ in
$\RR^4$ which are hyperbolic harmonic on\/ $A^4\setminus \RR^3$ and
real-analytic horizontally conformal submersions\/ $f:A^3 \to \CC$.
\end{theorem}

\begin{proof}
Let $\Cc$ be the conformal foliation on $A^3$ given by the level sets of
$f$, and let $\UU$ be the unit tangent vector field of $\Cc$ such
that $(\UU,\grad f_1, \grad f_2)$ is positively oriented.  Let $\JJ$
be the unique positive almost Hermitian structure on $A^3$ with
$\UU = \JJ(\pa/\pa x_0)$ and set $\phi = f$ on $A^3$.
Then, as in Theorem \ref{cxha:th:SFRconf}, the null lines
tangent to the vectors $\pa/\pa t + \UU$ define a shear-free ray
congruence $\ell$ on some open neighbourhood $A^M$ of $A^3$ in $\MM^4$;
we extend $\JJ$ and $\phi$ to $A^M$ by making them constant along the
leaves of $\ell$.  Write $\UU = \si^{-1}(\ii\mu)$, then $\mu$ satisfies
equations (\ref{cxha:EM}), and $\phi$ satisfies a similar pair of equations.
We extend these quantities to an open neighbourhoood of $A^M$ in $\CC^4$ by
analytic continuation, i.e., by insisting that they be
complex analytic, and finally we restrict to $\RR^4$.
We have thus extended $\JJ$ and $\phi$ to an open neighourhood $A^4$ of
$A^3$ in $\RR^4$; then $\mu$ satisfies (\ref{cxha:ER})
and $\phi$ satisfies similar equations, hence $\JJ$ is a Hermitian structure
on $A^4$ and $\phi$ is holomorphic with respect to $\JJ$.  Since $\UU = J(\pa/\pa x_0)$, we have
$\pa \phi/\pa x_0 = -\ii \, \UU(\phi) = 0$.  It follows from Proposition
\ref{cxha:prop:orthog-hypha} that $\phi$ is a hyperbolic harmonic
morphism.
\end{proof}

\begin{remark} \label{cxha:rem:Cinfty2}  \rm
If $f$ is only $C^{\infty}$, then, as in Proposition
{\rm \ref{cxha:prop:Cinfty1}}, if the distribution given by
$\,\mbox{\rm span} \{ \grad f_1, \grad f_2 \}$
is nowhere integrable, we can extend $f$
to one side of\/ $\RR^3$; precisely, there is an open subset $A^4$
of\/ $\RR^4_+ \cup \RR^3$ or $\RR^4_- \cup \RR^3$ with $A^4 \cap
\RR^3 = A^3$ and a $C^{\infty}$ map $\phi:A^4 \ra \CC$  with
$\phi\vert_{A^3} = f$ such that $\phi\vert_{{A^4} \setminus \RR^3}$ is a
hyperbolic harmonic morphism.
\end{remark}

\begin{corollary} \label{cxha:cor:minsurf}
{\rm (i)} Let\/ $\Cc$ be a real-analytic conformal foliation by curves of
an open subset of\/ $\RR^3$. Then there is a real-analytic
foliation of an open
subset\/ $A^4$ of\/ $\RR^4$ by surfaces which are minimal in\/ $A^4
\setminus \RR^3$ with respect to the hyperbolic metric and which
intersect\/ $\RR^3$ in leaves of\/ $\Cc$.

{\rm (ii)} Let\/ $c$ be an embedded real-analytic curve in\/ $\RR^3$. Then
there is an embedded real-analytic surface\/ $S$ in an open subset\/
$A^4$ of\/ $\RR^4$ which is minimal in\/ $A^4 \setminus \RR^3$ with
respect to the hyperbolic metric and which intersects\/ $\RR^3$ in\/ $c$.
\end{corollary}

\begin{proof}
(i) Represent the leaves of\/ $\Cc$ as the level curves of a
real-analytic horizontally conformal submersion $f:A^3 \to \CC$
on an open subset of $\RR^3$ and construct a hyperbolic harmonic
morphism $\phi$ as in the theorem; then its fibres give the desired
foliation.

(ii) Embed $c$ in a real-analytic conformal foliation by curves of
an open subset of $\RR^3$, as follows.  Construct the normal planes
to $c$ and integrate the vector field given by the normals to
these. This gives a foliation on an open neighbourhood of $c$ in
$\RR^3$ which has totally geodesic integrable horizontal distribution;
by Proposition {\bf\ref{pre:prop:confumbilic}}, it is Riemannian. (To
get a conformal foliation which is not Riemannian, replace the
planes by spheres, possibly of varying radii.) Now apply (i).
\end{proof}

\begin{example} \label{cxha:ex:circles-hyp}
Let $f$ be the horizontally conformal submersion of Example
\ref{cxha:ex:circles}. Recall that its level sets are given by
the leaves of the conformal foliation $\Cc$ with tangent vector
field $\UU$ given by (\ref{cxha:Ucircles}). The extension
of this to a shear-free ray congruence $\ell$ is described by
(\ref{cxha:NG}).  As in the proof of Theorem \ref{cxha:th:HWC-hyp},
extend $f$ to a function $\phi$ on an open subset of $\MM^4$ by
insisting that it be constant along the leaves of $\ell$.  Using
(\ref{cxha:circlesSFRsoln}) we see that this function is
\begin{eqnarray*}
\phi(t,x_1,x_2,x_3) &=& f\left(x_1,\,\,
\frac{r}{{x_2}^{\! 2}+{x_3}^{\! 2}}(rx_2+tx_3),\,\,
\frac{r}{{x_2}^{\! 2}+{x_3}^{\! 2}}(rx_3-tx_2) \right) \\
&=&-\ii x_1 + r \quad \text{where} \quad r = \sqrt{{x_2}^{\! 2}+{x_3}^{\! 2}-t^2} \:;
\end{eqnarray*}
this is smooth on the cone ${x_2}^{\! 2}+{x_3}^{\! 2} > t^2$.  It extends by
analytic continuation to the function
\begin{equation}
\phi(x_0,x_1,x_2,x_3) = -\ii x_1 + \sqrt{{x_0}^{\! 2}+{x_2}^{\! 2}+{x_3}^{\! 2}}
\label{cxha:circleshamo} \end{equation}
which is a complex-analytic function on suitable domains of
$\CC^4$.  Its restriction to $\RR^4$ is
smooth on $\RR^4 \setminus \{x_1\mbox{-axis} \}$ and defines the
hyperbolic harmonic morphism $\phi$ with boundary values at
infinity given by $f$.  See Example \ref{cxha:ex:circles-twistor} for further developments.
\end{example}

\section{Finding horizontally conformal functions} \label{cxha:sec:HCfns}

We now show how to use Theorem \ref{cxha:th:HWC-hyp} to find explicitly
the horizontally conformal submersion and conformal foliation by
curves on an open subset of $\RR^3$ which corresponds to a given
complex hypersurface $S$ of $\CP^3$.  In fact, by introducing a
parameter $\vec{a} \in \CC^4$, we can obtain a $5$-parameter family of
horizontally conformal submersions whose level sets give the
$5$-parameter family of conformal foliations of open subsets of
$\RR^3$ discussed in Remark \ref{cxha:rem:projn}.  For this, we need to
translate the hyperbolic metric to different slices.  First,
recall from Section \ref{cxha:sec:twistor} the map $\pi_{\vec{a}}:\CP^3 \to
\Rf_{\vec{a}}$ defined by $w \mapsto$ the intersection of the
$\alpha$-plane $\wt{w}$ given by (\ref{cxha:I}) with $\Rf_{\vec{a}}$. Then,
with the first component of ${\vec{a}}$ denoted by $a_0$, we have the following.

\begin{lemma} \label{cxha:lem:Theta_a}
Let\/ ${\vec{a}} \in \CC^4$. Equip\/ $\breve{\RR}^4_{\vec{a}} = \RR^4_{\vec{a}} \setminus
\RR^3_{\vec{a}}$ with the hyperbolic metric $$g_{\vec{a}}^H = \bigl( \tsum_{i=0}^3
\dd {x_i}^{\! 2} \bigr) \big/ \bigl(x_0- \Re a_0\bigr)^2,$$ and let\/ $N^5_{\vec{a}} =
\pi_{\vec{a}}^{-1}(\RR^3_{\vec{a}})$.  Then the kernel of the holomorphic contact
form
\begin{equation} \label{cxha:Theta_a}
\Theta_{\vec{a}}  = -2a_0(w_1 \,\dd w_0 - w_0 \,\dd w_1) +
w_1 \,\dd w_2 - w_2 \,\dd w_1 - w_0 \,\dd w_3 + w_3 \,\dd w_0 
\end{equation}
restricted to the manifold\/ $\CP^3 \setminus N^5_{\vec{a}}$ gives the horizontal
distribution of the map\/ $\pi_{\vec{a}}:\CP^3 \setminus N^5_{\vec{a}} \ra
(\breve{\RR}_{\vec{a}}^4, g^H_{\vec{a}})$.
\end{lemma}

\begin{proof}
Write ${\vec{a}}$ in null coordinates as
$(a_1,\wt{a_1}, a_2, \wt{a_2})$.   {}From Example \ref{cxha:ex:SL4C}
we see that the translation $T_{\vec{a}}:{\vec{x}} \mapsto {\vec{x}}+{\vec{a}}$ in $\CC^4$
corresponds to the map $\wt{T}_{\vec{a}}: \CP^3 \to \CP^3, \
[w_0,w_1,w_2,w_3] \mapsto [w_0,w_1,w_2 + a_1 {w}_0 -
\wt{a}_2 {w}_1, {w}_3 + a_2 {w}_0 + \wt{a}_1 {w}_1]$;
i.e., $\pi_{\vec{a}} \circ \wt{T}_{\vec{a}} = T_{\vec{a}} \circ \pi$.
In fact, this is a rephrasing of (\ref{cxha:transl}).
Then $\Theta_{\vec{a}} = (\wt{T}_{\vec{a}}^{-1})^*\Theta = \wt{T}_{-{\vec{a}}}^*
\Theta$; on calculating this, (\ref{cxha:Theta_a}) follows.
\end{proof}

Now let $S \subset \CP^3$ be a given complex hypersurface, and let ${\vec{a}}
\in \CC^4$. Let $U$ be an open set in $S$ such that $\pi_{\vec{a}}$ maps
$U$ diffeomorphically onto an open set $A^4$ of $\RR_{\vec{a}}^4$. Then
$U$ defines a Hermitian structure $\JJ$ on $A^4$ represented by the
section $w:A^4 \to U$ of $\pi_{\vec{a}}$. If $S$ is given by
\begin{equation}
\psi(w_0,w_1,w_2,w_3) = 0 \quad (\psi\mbox{ homogeneous complex-analytic})
\label{cxha:psi} \end{equation}
$w$ is given by solving
(\ref{cxha:K}).  Given a submersive holomorphic map $\zeta:U \to \Ci$, set
$\phi_{\vec{a}} = \zeta \circ w:A^4 \to \Ci$. Then $\phi_{\vec{a}}$ is
holomorphic with respect to $\JJ$ and so, by Proposition
{\bf\ref{twistor:prop:superint}}, it is a harmonic morphism with respect to the
hyperbolic metric on $\breve{\RR}_{\vec{a}}^4$ if and only if the level
surfaces of $\zeta$ are horizontal, i.e. tangent to
$\ker\Theta_{\vec{a}}$.

Set $\Si^S_{\vec{a}} = \{w \in S : \ker\Theta_{\vec{a}} = TS \}$.  Then on $S
\setminus \Si^S_{\vec{a}}$, \ $\ker\Theta_{\vec{a}} \cap TS$ is a
one-dimensional holomorphic distribution, so that its integral
(complex) curves foliate $S \setminus \Si^S_{\vec{a}}$. Note that
$\pi_{\vec{a}}(\Si^S_{\vec{a}}) \cap A^4$ is the set of K\"ahler points of $\JJ$,
i.e., the set $\{ {\vec{p}} \in A^4 : \nabla^H_w \JJ = 0 \ \forall w \in T_{\vec{p}} A^4 \}$;
here $\nabla^H$ denotes the Levi-Civita connection of the hyperbolic
metric on $\breve{\RR}^4_{\vec{a}}$.

To find $\phi_{\vec{a}}$ we proceed as follows.   Firstly, for an open subset $V$ of $U$,  let $c:V \ra
\CC^2$, \ $[w_0,w_1,w_2,w_3] \mapsto (\zeta,\eta)$ be complex
coordinates for $S$. Then we can
solve (\ref{cxha:I}) locally to find the composite map $c \circ w$, \ ${\vec{p}}
\mapsto (\zeta({\vec{p}}),\eta({\vec{p}}))$.

Next let $\wt{\zeta} = \wt{\zeta}(\zeta,\eta)$ be a
holomorphic function with $\pa\wt{\zeta}/ \pa\zeta \neq 0$
and set $\wt{\eta} = \eta$.  Then $(\zeta, \eta) \mapsto
(\wt{\zeta}, \wt{\eta})$ is locally a complex
analytic diffeomorphism.  The map $\phi_{\vec{a}}({\vec{p}}) =
\wt{\zeta}(\zeta({\vec{p}}),\eta({\vec{p}}))$ restricts to a hyperbolic
harmonic morphism on $\breve{\RR}^4_{\vec{a}}$ if and only if the level
sets of $\wt{\zeta}$ are superminimal, the condition for
this is
\begin{equation}
\Theta_{\vec{a}} \!\!\left( \frac{\pa}{\pa \wt{\eta}} \right) = 0 \,.
\label{cxha:h_a} \end{equation}
By the chain rule,
$$
\frac{\pa}{\pa \wt{\eta}} = \frac{\pa}{\pa \eta} -
\frac{\pa{\wt\zeta}}{\pa\eta}\frac{\pa}{\pa{\wt\zeta}} =  \frac{\pa}{\pa \eta} -
\biggl( \frac{\pa\wt{\zeta}}{\pa \eta} \bigg/ \frac{\pa\wt{\zeta}}{\pa \zeta}\biggr)\frac{\pa}{\pa \zeta}
$$
so that (\ref{cxha:h_a}) reads
\begin{equation}
\Theta_{\vec{a}}\!\! \left( \frac{\pa}{\pa \zeta} \right)
\frac{\pa\wt{\zeta}}{\pa \eta}
- \Theta_{\vec{a}} \!\! \left( \frac{\pa}{\pa \eta} \right)
\frac{\pa\wt{\zeta}}{\pa \zeta} =
0 \;. \label{cxha:E} \end{equation}
This equation can be solved to get a holomorphic function
$\wt{\zeta} = \wt{\zeta}_{\vec{a}} (\zeta,\eta)$ such that
$\phi_{\vec{a}}:{\vec{p}} \mapsto \wt{\zeta}_{\vec{a}} (\zeta({\vec{p}}),\eta({\vec{p}}))$
restricts to a submersive hyperbolic harmonic morphism on $(\breve{\RR}^4_{\vec{a}},
g^H_{\vec{a}})$.  Note that the solution is not unique; however, as the level
sets of $\wt{\zeta}_{\vec{a}}$ are uniquely determined on $V$, the map
$\phi_{\vec{a}}$ is unique up to postcomposition with a
biholomorphic function.  Restriction of $\phi_{\vec{a}}$ to $\RR^3_{\vec{a}}$ then gives a
horizontally conformal submersion $f$ on an open subset of
$\RR^3_{\vec{a}}$ and hence the conformal foliation by curves
corresponding to $S$ and ${\vec{a}}$. Note that the method actually finds
a holomorphic function $\phi_{\vec{a}}$ on an open subset of $\CC^4$.

In summary, given a complex hypersurface $S$ of $\CP^3$, there is
defined locally a real Hermitian structure $\JJ$ and all the
related distributions (\ref{cxha:Q}), in particular, $\UU = \JJ(\pa /\pa
x_0)$. Then given any ${\vec{a}} \in \CC^4$, the integral curves of $\UU$
define a real-analytic conformal foliation $\Cc_{\vec{a}}$ on an open subset of
the slice $\RR^3_{\vec{a}}$ and all such foliations are given this way.
Then we have shown above how to find a holomorphic function
$\phi_{\vec{a}}:A^C \to \CC$ on an open subset of $\CC^4$ such that the
level curves of $\phi_{\vec{a}}\vert_{\RR^3_{\vec{a}}}$ give the leaves of the conformal
foliation $\Cc_{\vec{a}}$.

\section{Examples} \label{cxha:sec:examples}

We start by considering twistor surfaces $S$ in $\CP^3$ which are linear.
The first example may be compared with Example {\bf 7.11.3}.

\begin{example} \label{cxha:ex:linear}
Let $S$ be the $\CP^2$ given by the zero set of the (homogeneous) linear function
\begin{equation} \label{cxha:linear-S}
\psi(w_0,w_1,w_2,w_3) = b_0 w_0 + b_1 w_1 + b_2 w_2 +b_3 w_3
\end{equation}
where $\vec{b} = (b_0,b_1,b_2,b_3) \neq (0,0,0,0)$.  
We consider the case
$$
[b_0,b_1,b_2,b_3] = [0,s,0,1]
$$
with $s$ a real parameter
so that (\ref{cxha:linear-S}) reads
\begin{equation} \label{cxha:Robinson}
s w_1 + w_3 = 0 \,.
\end{equation}
Note that $s = 0$ if and only if
$[b_0,b_1,b_2,b_3] \in N^5$; this will be a special case.
Equation (\ref{cxha:K}) reads
$$
s \mu + (q_2 + \mu \wt{q}_1) = 0\,,
$$
so that the direction field $\UU$ of the distributions (\ref{cxha:Q})
with twistor surface $S$ is given by $\UU = \si^{-1}(\ii\mu)$ where
$$
\mu = -q_2/(\wt{q}_1 + s)\,.
$$
Parametrize $S\setminus \{w_0 = 0 \}$ by
$$
(\zeta, \eta) \mapsto [1,w_1,w_2,w_3] =[1,\zeta,\eta, -s\zeta] \,.
$$
In terms of these parameters, we have from equation (\ref{cxha:Theta_a}),
\begin{eqnarray*}
\theta_{\vec{a}}  & = & 2a_0 \,\dd w_1 + w_1 \,\dd w_2 - w_2 \,\dd w_1 - \dd w_3 \\
   & = & 2a_0 \,\dd\zeta + \zeta \,\dd\eta - \eta \,\dd\zeta + s\,\dd\zeta
\end{eqnarray*}
so that equation (\ref{cxha:E}) reads
$$
 (-2a_0 + \eta - s)\frac{\pa\wt{\zeta}}{\pa\eta} +
\zeta \frac{\pa \wt{\zeta}}{\pa\zeta} = 0 \,;
$$
this has a solution
$$
\wt{\zeta} = - (-2a_0 + \eta - s)/\zeta \;.
$$
The incidence relations (\ref{cxha:I}) read
\begin{equation} \left. \begin{array}{rcl}
q_1 - \zeta \wt{q}_2 & = & \ \eta \\
q_2 + \zeta \wt{q}_1 & = & -s\zeta \end{array} \right\}
\end{equation}
with solution
$$
\zeta   =  -\frac{q_2}{\wt{q}_1 + s} \,, \quad
\eta  =  \frac{q_1 \wt{q}_1 + q_2 \wt{q}_2 +
q_1 s}{\wt{q}_1 + s}
$$
so that
\begin{equation} \label{cxha:linear}
\phi_{\vec{a}} = \wt{\zeta} =  \frac{ q_1 \wt{q}_1 + q_2 \wt{q}_2
-2a_0(\wt{q}_1 + s ) + (q_1 -\wt{q}_1)s -s^2}{q_2} \;.
\end{equation}
Set $a_0 = -\ii c$; this gives the horizontally conformal function
on the slice $\RR^3_{\vec{a}}$:
\begin{equation} \label{cxha:Robinson-HC}
\phi_c = \frac{(x_1+c)^2 + {x_2}^{\! 2} + {x_3}^{\! 2} - s^2
     + 2\ii(x_1 + c)s}{x_2+\ii x_3} \;.
\end{equation}
Write $\rho = \abs{s}$; if $s > 0$, the map $\phi_c$ is the composition
$$
\RR^3  \stackrel{T^1_c}{\lra}  \RR^3
\stackrel{D_{1/\rho}}{\lra} \RR^3 
\stackrel{\si^{-1}}{\lra} S^3
\stackrel{\ov{H}}{\lra}  S^2  \stackrel{\si}{\lra}  \Ci
$$
where $T^1_c$ is the translation $(x_1,x_2,x_3) \mapsto
(x_1+c,x_2,x_3)$, \ $D_{1/\rho}$ is the dilation $x \mapsto x/\rho$, \
$\si^{-1}$ is the inverse of
stereographic projection from $(-1,0,0,0)$,
and $\ov{H}$ is the
`conjugate' Hopf fibration given by $\si\circ \ov{H}(q_1, q_2) = \ov{q}_1/q_2$.  If $s < 0$ we replace
$\ov{H}$ by the Hopf fibration,  given by  $\si\circ H (q_1,q_2) = q_1/q_2$.  In either case
the fibres of $\phi$ give the classic conformal foliation $\Cc$
of $\RR^3$ by
the circles of Villarceau, see Fig.\ {\bf \ref{pre:fig:Hopf}}.
The shear-free ray congruence
$\ell$ whose projection is $\Cc$  (see Theorem \ref{cxha:th:SFRconf})
is called a \emph{Robinson congruence}. 
For $s = 0$, the foliation degenerates to the foliation by the
bunch of circles tangent to the $x_1$-axis at $(-c,0,0)$ described in
Example \ref{cxha:ex:bunch}, and $\ell$ is called a \emph{special
Robinson congruence}.  The hyperbolic harmonic morphism given by
\eqref{cxha:linear} with $a_0 = -\ii c$
and with boundary values given by \eqref{cxha:Robinson-HC} (with $s=0$ in
both cases) has fibres consisting of Euclidean hemispheres based on these circles.
Now note that $A \in \SU(4,h)$ transforms the coefficients of the
linear function (\ref{cxha:linear-S}) by
$$
(b_0,b_1,b_2,b_3) \mapsto (b_0,b_1,b_2,b_3) \, A \,.
$$
Since $\SU(4,h)$ acts transitively on the three subsets $\CP^3_+,
\CP^3_-$ and $N^5$ of $\CP^3$,
if $[b_0,b_1,b_2,b_3] \in \CP^3_+$ (respectively, $\CP^3_-$), the linear
function is $\SU(4,h)$-equivalent to (\ref{cxha:Robinson})
with $s > 0$ (respectively, $s < 0)$; if
$[b_0,b_1,b_2,b_3] \in N^5$, it is equivalent to (\ref{cxha:Robinson})
with $s=0$.  

Note that the special Robinson congruence is conformally equivalent via the
inversion (Example \ref{cxha:ex:SL4C}(iv)) to (\ref{cxha:linear-S}) with $\vec{b} =
(0,1,0,0)$;  this gives $\mu = 0$ and corresponds to a shear-free ray
congruence consisting of parallel rays.  The corresponding Hermitian structure $\JJ$ is
K\"ahler; in fact, any globally defined Hermitian structure on $\RR^4$ is K\"ahler, see
Lemma {\bf\ref{twistor:lem:HermR4}}.  

It is easy to see that any linear hypersurface (\ref{cxha:linear-S}), 
whether $\vec{b}$ is in $N^5$ or not, is $\SL(2,\RH)$-equivalent
to (\ref{cxha:linear-S}) with $\vec{b} = (0,1,0,0)$;  so
that any two Hermitian structures on a subset of $S^4$ given by a
linear equation are conformally equivalent.  
\end{example}

\medskip

We now give three examples where the twistor surface is a quadric
so that we can solve all equations explicitly.

\begin{example} \label{cxha:ex:radial-twistor}
Let the twistor surface be the quadric 
$$
S =
\{ [w_0,w_1,w_2,w_3] \in \CP^3 : w_0 w_3 - w_1 w_2 = 0 \}\,.
$$
 Then
the direction field $\UU$ of the corresponding distributions (\ref{cxha:Q})
is given by 
$\UU = \si^{-1} (\ii \mu)$ where
$$
q_2 + \mu \wt{q}_1 - \mu (q_1 - \mu \wt{q}_2) = 0 \,;
$$
this has solutions
\begin{equation}
\mu  =  \frac{q_1 - \wt{q}_1 \mp
\sqrt{(q_1-\wt{q}_1)^2
- 4q_2 \wt{q}_2}} {2\wt{q}_2}
\label{cxha:mu1} \end{equation}
so that
$$
\ii \mu   = 
\frac{- x_1 \pm \abs{x}}{x_2 - \ii x_3}
 =  \frac{x_2 + \ii x_3}{x_1 \pm \abs{x}}
$$
where $\abs{x} = \sqrt{{x_1}^{\! 2}+{x_2}^{\! 2}+{x_3}^{\! 2}} \,.$
Thus, on any slice  $\RR^3_t$, \ $\UU$ is given by
$$
\UU(t,x_1,x_2,x_3) = \pm (x_1,x_2,x_3)
 / \abs{x}
$$
which gives the tangent vector field of the foliation by radial lines of
Example \ref{cxha:ex:radial}.

Carrying out the calculations of Section \ref{cxha:sec:HCfns} gives
$\phi_{\vec{a}} = \mu$ for all $\vec{a}$; we omit the details which are
similar to the next example. Then, for each ${\vec{a}}$, the map $\si^{-1}
\circ \phi_{\vec{a}}$ restricted to $\breve{\RR}^4_{\vec{a}}$ gives the harmonic
morphism $\RH^4 = (\breve{\RR}^4_{\vec{a}}  \setminus \{x_0 \mbox{-axis} \},
g_{\vec{a}}^H ) \ra S^2$ given by orthogonal projection to $\RR^3_{\vec{a}}$
followed by radial projection.
\end{example}

\begin{figure}[htb]
\centerline{\includegraphics*[bb=140 58 580 480,scale=0.5]{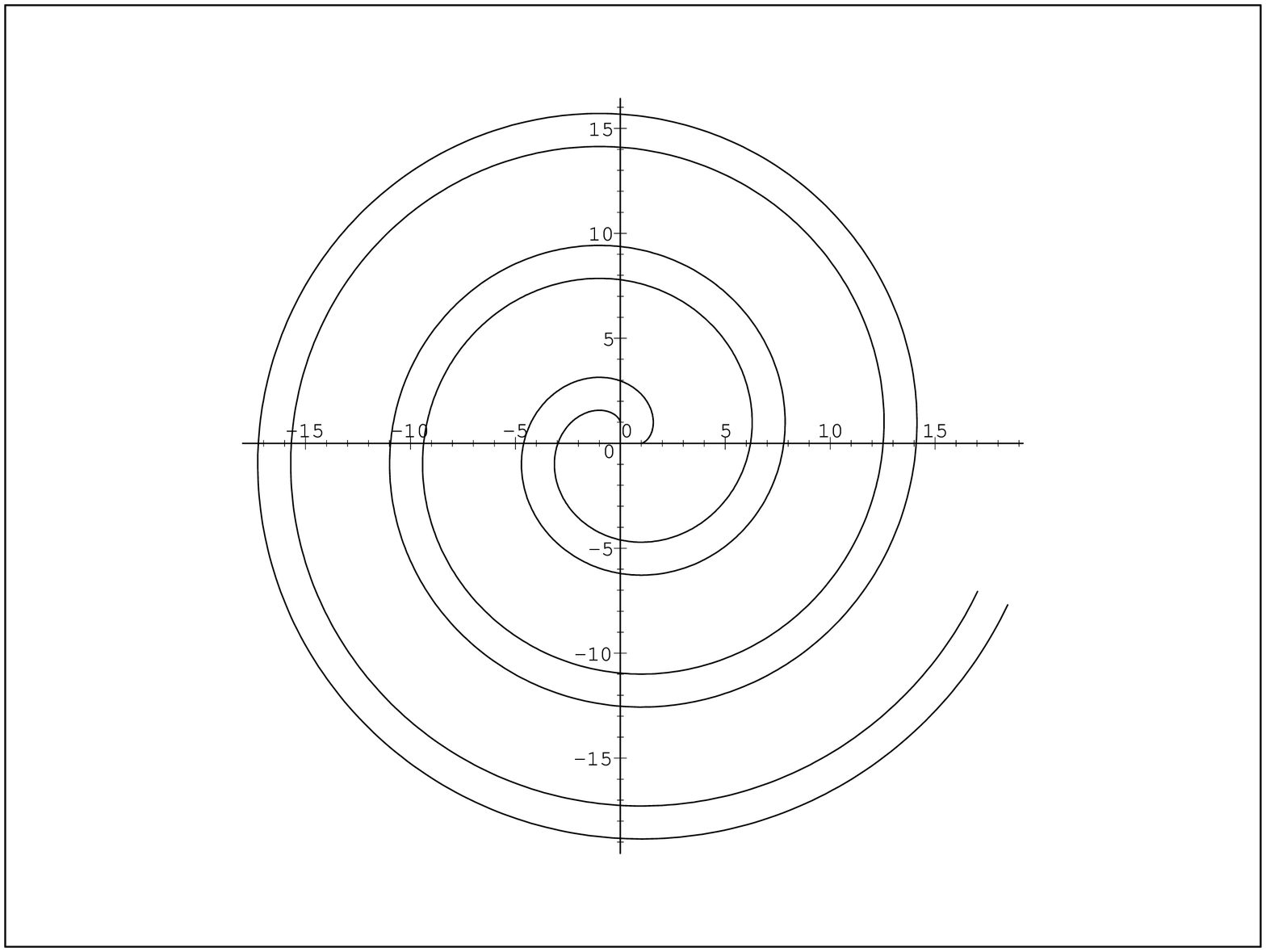}}
\caption{Two leaves of the Riemannian foliation $\Cc_t$
with $t=1$.}
\label{cxha:fig:evolutes}
\end{figure}

\begin{example} \label{cxha:ex:circles-twistor}
Let the twistor surface be the quadric
$$
S = \{[w_0,w_1,w_2,w_3] \in \CP^3 : w_0 w_3 + w_1 w_2 = 0 \}\,.
$$
 Then
the direction field $\UU$ of the corresponding distributions is given
by the vector field $\UU = \si^{-1}(\ii \mu)$ where
$$
q_2 + \mu \wt{q}_1 + \mu (q_1 - \mu \wt{q}_2) = 0 \,;
$$
this has solutions
$$
\mu  =  \frac{q_1 + \wt{q}_1 \pm
\sqrt{(q_1+\wt{q}_1)^2 + 4q_2 \wt{q}_2}}
{2\wt{q}_2} =  \frac{x_0 \pm s}{x_2 - \ii x_3}
$$
where we write $s = \sqrt{{x_0}^{\! 2} + {x_2}^{\! 2} + {x_3}^{\! 2}}$.  
Note that
$$
\mu\vert_{\RR^3_{\vec{0}}} = \pm \frac{\sqrt{{x_2}^{\! 2} +{x_3}^{\! 2}}}{x_2-\ii x_3}
= \pm \frac{x_2+\ii x_3}{\sqrt{{x_2}^{\! 2} +{x_3}^{\! 2}}}
$$
so that on $\RR^3_{\vec{0}}$\,,
$$
\UU(x_1,x_2,x_3) = \pm \frac{1}{\sqrt{{x_2}^{\! 2} +{x_3}^{\! 2}}}
(0,-x_3,x_2)\,;
$$
this is the tangent vector field of the conformal foliation
$\Cc$  by circles around the $x_1$-axis of Example
\ref{cxha:ex:circles}. Further note that, writing
$r = \sqrt{{x_2}^{\! 2} +{x_3}^{\! 2} - t^2}$, the restriction
$\mu\vert_{\MM^4}$ is given by
\begin{eqnarray*}
\mu(t,x_1,x_2,x_3) & = & \frac{-\ii t \pm r}{x_2-\ii x_3} \\
& = & \frac{r}{{x_2}^{\! 2} + {x_3}^{\! 2}}
\left\{ \left( \pm x_2 + \frac{t}{r}x_3 \right) + \ii
\left(\pm x_3 -\frac{t}{r}x_2 \right) \right\} \,,
\end{eqnarray*}
and so, on the open set ${x_2}^{\! 2} +{x_3}^{\! 2} > t^2$, we have
$$
\UU(t,x_1,x_2,x_3) = \frac{r}{\sqrt{{x_2}^{\! 2} + {x_3}^{\! 2}}} \left(
0,\, \mp x_3 + \frac{t}{r}x_2,\, \pm x_2 + \frac{t}{r} x_3 \right)
\,;
$$
then $\ww = \pa/\pa t + \UU$ gives the tangent field of the shear-free
congruence $\ell$ extending $\Cc$ discussed in Example
\ref{cxha:ex:circles}.

Parametrize $S$ away from $w_0 = 0$ by $(\zeta,\eta) \mapsto
[1,w_1,w_2,w_3] = [1,\eta,-\zeta, \zeta \eta]$.  Then the incidence
relations (\ref{cxha:I}) read
$$
\left. \begin{array}{ccr}
q_1 - \eta \wt{q}_2 & = & - \zeta \\
q_2 + \eta \wt{q}_1 & = & \zeta \eta 
\end{array} \right\}\,;
$$
these have solutions
$$
\zeta = -\ii x_1 \pm s \quad \mbox{and} \quad
\eta = \frac{x_0 \pm s}{x_2 - \ii x_3} \,. 
$$  
On $\MM^4$,
\begin{equation}
\zeta = -\ii x_1 \pm r \quad \mbox{and} \quad \eta =
\frac{-\ii t \pm r}{x_2-\ii x_3}
\label{cxha:zeta-eta-1} \;.\end{equation}
These solutions are real analytic except on the cone ${x_2}^{\! 2} +
{x_3}^{\! 2} = t^2$.  We obtain smooth solutions if we avoid this
cone.

We remark that
the $3$-manifold $S \cap N^5$ has equation:
$\bigl(\zeta+\ov{\zeta}\bigr)\bigl(\abs{\eta}^2-1\bigr) = 0$ and so consists of the
union of two submanifolds:  $S_1:\Re \zeta = 0$ and $S_2: \abs{\eta}
= 1$.  Smooth branches of $(\zeta,\eta)$ will correspond to CR
maps $A \to \CP^3$ with image lying in $S_1 \setminus S_2$
or $S_2 \setminus S_1$; \ $S_1 \cap S_2$ corresponds to the
\emph{branching set} (cf.\ Chapter {\bf 9})
of $(\zeta,\eta)$ which is the subset of $A$ on
which the square root $r$ vanishes.  Explicitly, note from
{\rm (\ref{cxha:zeta-eta-1})} that, if $t > \abs{q_2}$ then $\Re \zeta = 0$,
this corresponds to points of $S_1$; if $t < \abs{q_2}$ then
$\abs{\eta} = 1$, this corresponds to points of $S_2$.

We have
\begin{eqnarray*}
\theta_{\vec{a}} & = & 2a_0 \,\dd w_1 + w_1 \,\dd w_2 - w_2 \,\dd w_1 - \,\dd w_3 \\
& & 2a_0 \,\dd\eta - \eta \,\dd\zeta + \zeta \,\dd\eta - \zeta \,\dd\eta {-} \eta \,\dd\zeta
\end{eqnarray*}
so that (\ref{cxha:E}) reads
$$
\eta \frac{\pa \wt{\zeta}}{\pa \eta} + a_0
\frac{\pa\wt{\zeta}}{\pa \zeta} = 0 \,;
$$
this has a solution
\begin{equation}
\wt{\zeta}_{\vec{a}} = \zeta - a_0 \ln \eta \;.
\label{cxha:sol1} \end{equation}
Set $\phi_{\vec{a}} = \wt{\zeta}_{\vec{a}} \big(\zeta({\vec{p}}),\eta({\vec{p}}) \big)$;
then this is the complex-valued map on a dense subset of $\CC^4$ given by
\begin{equation} 
\phi_{\vec{a}}(x_0,x_1,x_2,x_3) = -\ii x_1 \pm s
- a_0 \ln \frac{x_0 \pm s}{x_2-\ii x_3} \,.
\label{cxha:phia2} \end{equation}
For any ${\vec{a}} \in \CC^4$, this restricts to a complex-valued hyperbolic harmonic
morphism on a dense subset of $(\breve{\RR}_{\vec{a}}^4, g_{\vec{a}}^H)$. In
particular, when ${\vec{a}} = {\vec{0}}$, this simplifies to the hyperbolic
harmonic morphism on $\breve{\RR}^4\setminus \{x_1 \mbox{-axis}\}$
given by \eqref{cxha:circleshamo}.
As in Example \ref{cxha:ex:circles-hyp}, this further restricts on $\RR^3$ to
\begin{equation}
\phi(x_1,x_2,x_3) = -\ii x_1 \pm \sqrt{{x_2}^{\! 2} + {x_3}^{\! 2}}\,,
\label{cxha:L} \end{equation}
the level surfaces of which give the conformal foliation by
circles round the $x_1$-axis of  Example \ref{cxha:ex:circles}. The
harmonic morphism \eqref{cxha:circleshamo} has fibres given by the
Euclidean spheres having these circles as equators, these
spheres are totally geodesic in $(\breve{R}^4, g^H)$.

For definiteness, take plus signs in the above.
Set $a_0 = -\ii t$ in (\ref{cxha:phia2}).  Then the restriction of $\phi_{\vec{a}}$ to the
open set $\{(x_1,x_2,x_3) : {x_1}^{\! 2} + {x_2}^{\! 2} > t^2 \} \subset \RR^3_t
= \RR^3_{\vec{a}}$  is the horizontally conformal map
\begin{eqnarray*}
\phi_t  =  \phi_{\vec{a}} & = & -\ii x_1 + r + \ii t \ln
\frac{r-\ii t}{x_2- \ii x_3} \\
       & = & -\ii x_1 + r - t \arg \frac{r- \ii t}{x_2- \ii x_3} \;;
\end{eqnarray*}
the level curves of this are the leaves of the conformal (in fact, Riemannian)
foliation $\Cc_t$ of
Example \ref{cxha:ex:circles}; they lie in the planes $x_1 =$ constant \ and are
the involutes of the unit circle, see Fig.\ \ref{cxha:fig:evolutes}.
\end{example}

\begin{remark} \label{cxha:rem:cxsame}  \rm
The matrix $P = \mbox{\rm diag}(1,\ii,\ii,1)$---or rather
$P = \mbox{\rm diag}(\theta,\ii\theta,\ii\theta,\theta)$ where $\theta^4
= -1$ so that
$P \in \SL(4,\CC)$---transforms Example
{\rm \ref{cxha:ex:radial-twistor}} into Example
{\rm \ref{cxha:ex:circles-twistor}}.
Thus they give conformally equivalent foliations by $\alpha$-planes on
subsets of\/ $\CC^4$; however, $P$ does not lie in $\SU(4,h)$ or
$\SL(2,\HH)$; in fact we can easily see that the other pairs of distributions
associated to this pair of foliations as in {\rm (\ref{cxha:Q})} are not conformally equivalent.
\end{remark}

\begin{figure}[htb]
\centerline{\includegraphics*[bb=120 170 500 580,scale=0.5]{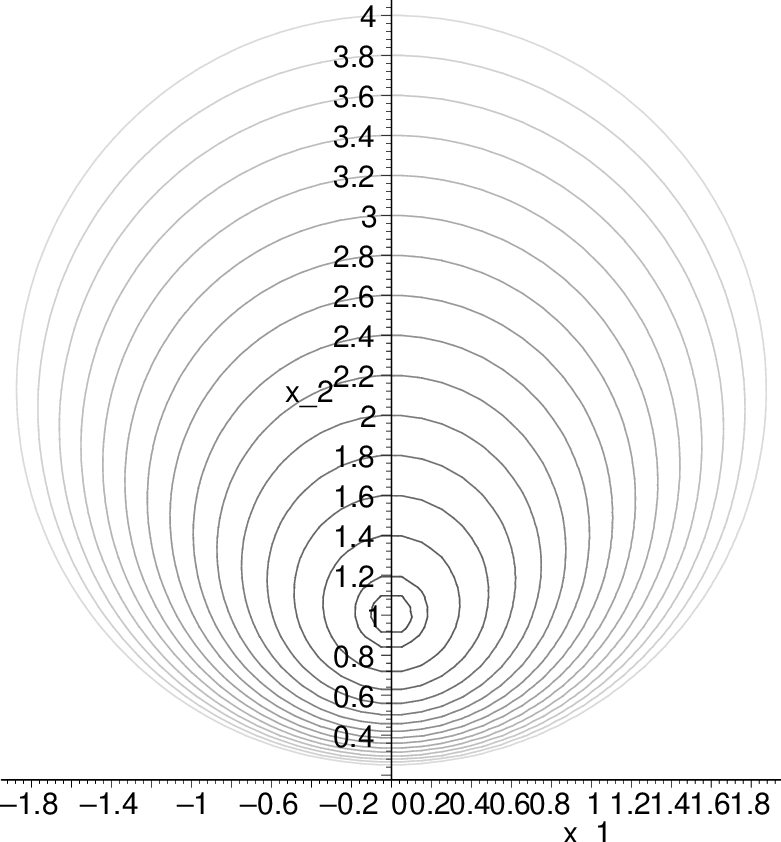}}
\caption{Some leaves of the foliation $\Cc_0$ in the plane $\arg(x_2+\ii x_3) = 0$.}
\label{cxha:fig:contours}
\end{figure}

\begin{example}
Let the twistor surface be the quadric 
$$
S = \{ [w_0,w_1,w_3,w_4] \in
\CP^3 : w_0 w_1 + w_2 w_3 = 0 \}\,.
$$
  Then $\mu$ satisfies
$$
\mu + (q_1 - \mu \wt{q}_2)(q_2 + \mu \wt{q}_1) = 0 \,;
$$
this has solutions
$$
\mu = \frac{(1 + q_1 \wt{q}_1 - q_2 \wt{q}_2) \pm
\sqrt{(1 + q_1 \wt{q}_1 - q_2 \wt{q}_2)^2 + 4 q_1 \wt{q}_1 q_2 \wt{q_2}}}
{2\wt{q}_1\wt{q}_2} \;.
$$
Parametrize $S$ away from $w_0 = 0$ by
$$
(\zeta,\eta) \mapsto [1,w_1,w_2,w_3] = [1,\, \zeta\eta,\, -\eta,\, \zeta] \,.
$$
Then the incidence relations (\ref{cxha:I}) read
$$
\left.
\begin{array}{ccc}
q_1 - \zeta\eta\, \wt{q}_2 & = & -\eta  \\
q_2 + \zeta\eta\, \wt{q}_1 & = & \zeta 
\end{array}
\right\}\,;
$$
solving these gives
\begin{equation}
\zeta = \frac{1+q_1 \wt{q}_1 + q_2 \wt{q}_2 \pm s}{2\wt{q}_2}
\label{cxha:Q1} \end{equation}
where $s = \sqrt{(1 + q_1 \wt{q}_1 + q_2 \wt{q}_2)^2 - 4 q_2 \wt{q_2}}$.
This has branching set when the square root $s$ is zero; on
$\RR^4 = \RR^4_{\vec{0}}$ this occurs on
$C = \{(q_1,q_2) \in \CC^2 = \RR^4: q_1 q_2 =0, \ \abs{q_1}^2 + \abs{q_2}^2 = 1\}$;
this set is the disjoint union of two circles.
On any open set $A \subseteq \RR^3_{\vec{0}} \setminus
C$, on fixing the $\pm$ sign in (\ref{cxha:Q1}), we obtain a smooth
solution. Next note that
\begin{eqnarray*}
\theta_{\vec{a}} & = & 2a_0 \,\dd w_1 + w_1 \,\dd w_2 - w_2 \,\dd w_1 - \dd w_3 \\
  & = & 2a_0(\eta \,\dd\zeta + \zeta \,\dd\eta) - \zeta\eta \,\dd\eta +
  \eta (\eta \,\dd\zeta + \zeta \,\dd\eta) - \dd\zeta \\
  & = & (\eta^2 + 2a_0 \eta - 1)\,\dd\zeta + 2a_0 \zeta \,\dd\eta\,,
\end{eqnarray*}  
so that (\ref{cxha:E}) reads
$$
(\eta^2 + 2a_0\,\eta - 1) \frac{\pa\wt{\zeta}}{\pa\eta}
      -2a_0 \,\zeta \frac{\pa\wt{\zeta}}{\pa\zeta} = 0 \,;
$$
this has a solution
$$
\wt{\zeta} = \zeta \left\{ \frac{\eta + a_0 +
\sqrt{a_0^{{}2} +1}} {\eta + a_0 - \sqrt{a_0^{{}2} + 1}} \right\}
^{-a_0/\sqrt{a_0^{{}2} + 1}} \;.
$$
When $a_0 = 0$ this gives $\wt{\zeta} = \zeta$, so that, for
any smooth branch $\zeta:A \ra \CC$ of (\ref{cxha:Q1}), the level
curves of $\phi = \zeta\vert_A$ give a conformal foliation. 
Explicitly,
$$
\phi(x_1,x_2,x_3) = \frac{1+{x_1}^{\! 2}+{x_2}^{\! 2}+{x_3}^{\! 2} \pm s}{2(x_2 + \ii
x_3)}
$$
with $s = \sqrt{(1 + {x_1}^{\! 2}+ {x_2}^{\! 2}+ {x_3}^{\! 2})^2
	- 4({x_2}^{\! 2}+{x_3}^{\! 2})}$\;; 
this defines a rotationally symmetric foliation $\Cc_0$.   In fact, on each plane $\arg(x_2+\ii x_3) = \const$
it restricts to a foliation by coaxal circles as shown in Fig.\ \ref{cxha:fig:contours} (the plus and minus
signs give the same foliation).
\end{example}

\section{Notes and comments} \label{cxha:sec:notes}  \small

\smallskip

\noindent \emph{Section} \ref{cxha:sec:SFR}

\smallskip
\noindent
The \emph{shear tensor} of a foliation
$\ell$ of a Lorentzian manifold by null lines is defined by
$S(X,Y) =$ the trace-free part of
\begin{equation} \label{cxha:shear}
(X,Y) \mapsto \half g(\nabla_X Y + \nabla_Y X, \ww) \qquad
\big(X,Y \in C^{\infty}(\Si) \big).
\end{equation}
Geometrically it measures how fast infinitesimal
circles in the screen space distort into ellipses under Lie transport
along the congruence.  It is easy to see that $\ell$ is shear-free
if and only if its shear tensor vanishes.

The \emph{expansion} tensor is defined as the trace of
(\ref{cxha:shear}).  Together with the twist tensor discussed in
Remark \ref{cxha:rem:twist}, these are the three 
fundamental tensors associated to a congruence of null
lines.  For more information see Benn (1994), Beem, Ehrlich and Easley (1996).

\bigskip

\noindent \emph{Section} \ref{cxha:sec:cxha}

\smallskip
\noindent
Consider the set of complex-analytic maps $\CC^m \supset A^C \to \CC^n$ which satisfy
$$
\sum_{i=1}^m \frac{\pa^2\phi}{\pa {x_i}^{\! 2}}  =  0
	\quad \text{and} \quad
\sum_{i=1}^m \frac{\pa \phi^{\al}}{\pa x_i}
	\frac{\pa \phi^{\be}}{\pa x_i} = \la\, \delta_{\al\be} \qquad 
\big( (x_1,\ldots,x_m) \in A^C,  \quad  \al,\be \in
\{1,\ldots,n\} \big) \,.
$$
By a complex-analytic version of Theorem {\bf\ref{fund:th:char}}, 
such maps can be characterized as those complex-analytic maps which pull
back complex-harmonic functions to complex-harmonic functions,
hence we shall call them \emph{complex-harmonic morphisms}.  If $\la
\equiv 0$, they pull back \emph{arbitrary} complex-analytic functions to
complex-harmonic functions; 
we shall call them \emph{degenerate complex-harmonic morphisms}.
When $m=4, n=1$ this is the case studied in Section \ref{cxha:sec:cxha};
for arbitrary $m$ and $n=1$, they are the analogue of the degenerate harmonic morphisms studied
in Sections {\bf\ref{semi:sec:HWC}} and  {\bf\ref{semi:sec:hamorph}}.
See Pambira (2002, 2003) for more on complex-harmonic morphisms, as well as
harmonic morphisms between \emph{degenerate} semi-Riemannian manifolds.

\bigskip

\noindent \emph{Section} \ref{cxha:sec:ha-morph-SFR}

\smallskip
\noindent
The correspondence between Hermitian structures and harmonic morphisms from a $4$-dimensional Einstein manifold was
established by the second author (Wood 1992); see also \ Chapter {\bf\ref{ch:twistor}}.

\bigskip

\noindent \emph{Section} \ref{cxha:sec:twistor}

\smallskip

\noindent 1.  An invariant description of the compactification of $\RR^n_p$ for any
dimension and signature is given by Cahen, Gutt and Trautman (1993); see
Akivis and Goldberg (1996), Ward and Wells (1990),
Fillmore (1977), Mason and Woodhouse (1996) for other approaches.

\smallskip

\noindent 2.  The double holomorphic fibration \eqref{cxha:doublefibr} is the basis of the 
well-known \emph{twistor correspondence} of relativity theory, which has led to \emph{twistor  
theory} as pioneered by Penrose and others, see Penrose (1967, 1975), Penrose and Rindler (1987, 1988) and 
Huggett, Mason, Tod, Tsou and Woodhouse (1998) for a more recent compendium of articles.  The \emph{Twistor Newletter}
has provided a forum for discussion on the subject and contributions have been collected into three volumes:
Mason and Hughston (1990), Mason, Hughston and Kobak (1995), Mason, Hughston, Kobak  and Pulverer (2001).

The space $\CP^3$ is thought of as (projectivized) \emph{twistor space} and points of this space become 
the fundamental objects
from a physical point of view.  Points of the real quadric hypersurface $N^5\subset \CP^3$ defined by 
\eqref{cxha:N5} are called \emph{null twistors}; these correspond to null geodesics in Minkowski space (see
Section \ref{cxha:sec:twistor}).  At one level, the correspondence is merely formal, but the distinction 
becomes fundamental when quantization takes place.  As Penrose describes
the situation `... in the usual view, one imagines a space-time in which the points remain intact but where 
(the metric) $g_{ab}$ becomes quantised ... according to the twistor view ... the concept of a 
twistor remains intact, while that of a space-time point becomes fuzzy.'

One motivation for this complex-analytic view of the world is the empirically observed 
discreteness in nature, e.g., charge, spin etc., which has its analogue in the discrete values arising 
from contour integration.  This is the basis of the so-called \emph{twistor diagrams}, which give 
a combinatorial description of particle interactions (Penrose 1975, Sparling 1975).  

In general, defining twistors on a curved space-time presents formidable difficulties. Two possible extensions to this case are \emph{hypersurface twistors} (Penrose and Rindler 1988, part of Section 7.4)
and \emph{asymptotic twistors} (Penrose and Rindler 1988, part of Section 9.8).    See also Mason, Hughston, Kobak and Pulverer (2001), and the next Note.

\smallskip

\noindent 3.  The extension of the ideas in this paper to more general manifolds
will be complicated by the presence of curvature which means that
(i) there are shear-free ray congruences only if the curvature is
sufficiently special (Goldberg--Sachs Theorem), (ii)
Lemma \ref{cxha:lem:Sachs} no longer
holds---instead we have the Sachs equations (Penrose and Rindler 1988,
Chapter 7),  (iii) the CR structure on the unit
tangent bundle of a $3$-manifold may not be `realizable' as a
hypersurface CR structure (LeBrun 1984).
For discussions of these and related matters, see also
Trautman (1985), Robinson and Trautman (1986),
Lewandowski, Nurowski and Tafel (1990), Nurowski (1996).

\bigskip

\noindent \emph{Section} \ref{cxha:sec:Kerr}

\smallskip

\noindent  The original Kerr Theorem (Kerr 1963), expresses the correspondence
between a complex-analytic surface in the twistor space $\CP^3$ and a shear-free null geodesic congruence 
in Minkowski space.  This correspondence permits one to express solutions to the 
zero-rest-mass field equation in terms of contour integrals of 
holomorphic functions defined on domains of $\CP^3$ (\emph{twistor functions}), see Penrose (1969). 
By analogy, Baird (1993) (see also Baird 2002) expresses harmonic morphisms defined on the three-dimensional space 
forms in terms of contour integrals of functions defined on the corresponding mini-twistor space
(cf.\ Chapter {\bf\ref{ch:global}}).

\bigskip

\noindent \emph{Section} \ref{cxha:sec:CR}

\smallskip

\noindent 1.  For general information on CR manifolds and embeddability and extension
problems, see Jacobowitz (1990) and Boggess (1991).

\smallskip

\noindent 2.  For descriptions of the circles of Villarceau, see Wilker (1986)
and compare with a different treatment in Baird (1998).

\smallskip

\noindent  3.  For other descriptions of the CR structure on the unit tangent bundle
of a Riemannian $3$-manifold see Sato and Yamaguchi (1989,
Proposition 5.1) and LeBrun (1984).
Note that this is \emph{not} the same as
the CR structure on the unit tangent bundle of a Riemannian manifold of
arbitrary dimension discussed, for
example, in Blair (1976, Chapter 7) or Tanno (1992).

\bigskip

\noindent \emph{Section} \ref{cxha:sec:hypha}
\smallskip

\noindent
1.  A description of harmonic morphisms on hyperbolic space $\RH^4$ in terms of
holomorphic data---essentially the twistorial formulation described
in Sections \ref{cxha:sec:ha-morph-SFR} and \ref{cxha:sec:twistor}---is given by Baird (1992). 
\smallskip

\noindent
2. Thinking of $H^m$ as the unit disc with the Poincar\'e metric
(see Example {\bf \ref{pre:ex:std-spaces})},
Li and Tam (1991) show that,  given any $C\sp 3$ map from the boundary at infinity
$S^{m-1}$ of $\RH^m$ to the boundary at infinity $S\sp {n-1}$ of $\RH\sp n$, whose energy density
is nowhere zero, there exists a harmonic map from $H\sp m$ to $H\sp n$ which realizes the given boundary
map; see Li and Tam (1993) for regularity and uniquenes, Li, Tam  and Wang (1995) for an extension to Hadamard surfaces,
and Donnelly (1994, 1999) for the case of
maps between rank one symmetric spaces.
\smallskip

\noindent
2.  LeBrun (1984) shows that any real-analytic $3$-manifold $M^3$ is the boundary at infinity of some anti-self-dual Einstein
$4$-manifold $P^4$.  In view of the description of harmonic morphisms with values in a surface on such a 
$4$-manifold,  given in Sections {\bf\ref{twistor:sec:hamo-super}} and {\bf\ref{twistor:sec:converse}} in terms 
of integrable Hermitian structures,  it is tempting to conjecture that 
a horizontally weakly conformal map $\phi : M^3\ra N^2$ onto a Riemann surface is the boundary
values at infinity of a 
harmonic morphism $\Phi : P^4\ra N^2$.

\section{References}

\small

\noindent Articles marked with an asterisk (*) are, 
in part at least, about harmonic
morphisms.  This part of the list of references is based on the
\emph{Bibliography of Harmonic Morphisms}
regularly updated by Sigmundur Gudmundsson, see
{\tt http://www.maths.lth.se/matematiklu/personal/sigma/harmonic/bibliography.html}.

\smallskip

The numbers of the sections in which articles are referred to are shown at the end of each entry;
those in italics refer to the `Notes and comments' for that section.

\bigskip

\rm

\parindent=-2em \leftskip=2em \par

Akivis, M. A. and Goldberg, V. V. (1996).
{\it Conformal differential geometry and its generalizations.}
Pure and Applied Mathematics.
A Wiley-Interscience Publication.
Wiley, New York.
MR~98a:53023.  Zbl.863.53002.
\secno {\it S.9}.

Baird, P. (1992{\it b\/})*. Riemannian twistors and Hermitian structures
on low-dimen\-sional space forms. {\it J. Math. Phys.}, {\bf  33}, 3340--55.
MR~93h:53069,  Zbl.763.32017.
\secno {\it S.13}

Baird, P. (1993)*. Static fields on three-dimensional space
forms in terms of contour integrals of twistor functions. {\it Phys.  Lett.
A}, {\bf  179},  279--83. 
MR~94j:58042.  
\secno {\it S.11}.

Baird, P. (1998)*.  Conformal foliations by circles and complex
isoparametric functions on Euclidean $3$-space.
{\it Math. Proc. Cambridge Philos. Soc.}, {\bf 123}, 273--300.
MR~99h:53031,  Zbl.899.58012.
\secno {\it S.12}.

Baird, P.  (2002)*.  An introduction to twistor theory.  Notes from the conference
{\it Integrable systems and quantum field theory} held at Peyresq, France, June 2002, see\\
{\tt http://www.univ-mlv.fr/physique/phys-math/baird-twistors.pdf}
\secno {\it S.11}.

Baird, P. and Wood, J.~C. (1998)*.  Harmonic morphisms,
conformal foliations and shear-free ray congruences.  {\it Bull. Belg.
Math. Soc.}, {\bf 5}, 549--64. 
MR~2000d:53098,  Zbl.990.13511.
\secno S.0.

Baird, P. and Wood, J. C. (2003)
{\it Harmonic morphisms between Riemannian manifolds.}
London Mathematical Society Monographs, No. 29, 
Oxford University Press;
see {\tt http://www.amsta.leeds.ac.uk/Pure/\\staff/wood/BWBook/BWBook.html}.
\secno S.0--3,~S.6--9,~S.11--15,~{\it S.7},~{\it S.11},~{\it S.13}. 

Beem, J.~K., Ehrlich, P.~E. and Easley, K.~L.  (1996).
{\it Global Lorentzian geometry}.
Monographs and Textbooks in Pure and Applied Mathematics, vol.~202, 2nd edn.
Marcel Dekker, New York.
MR~97f:53100,  Zbl.846.53001.
\secno {\it S.3}.

Benn, I. M.  (1994).
A unified description of null and nonnull shear-free congruences.
{\it J. Math. Phys.} {\bf 35}, 1796--1802.
MR~95a:53103,  Zbl.808.53018.
\secno {\it S.3}.

Blair, D. E.  (1976).
{\it Contact manifolds in Riemannian geometry.}
Lecture Notes in Mathematics, vol. 509. Springer, Berlin.
MR~57\#7444,   Zbl.319.53026.
\secno {\it S.12}.

Boggess, A.  (1991).
{\it CR manifolds and the tangential Cauchy-Riemann complex}.
Studies in Advanced Mathematics. CRC Press, Boca Raton, FL.
MR~94e:32035,  Zbl.760.32001.
\secno {\it S.12}.

Cahen, M., Gutt, S. and Trautman, A. (1993).
Spin structures on real projective quadrics.
{\it J. Geom. Phys.}, {\bf 10}, 127--54. 
MR~94c:57041,  Zbl.776.57011.
\secno {\it S.9}.

Donnelly, H. (1994).
Dirichlet problem at infinity for harmonic maps: rank one symmetric spaces.
{\it Trans. Amer. Math. Soc.} {\bf 344}, 713--35.
MR 95c:58045,  Zbl.812.58020.
\secno {\it S.13}.

Donnelly, H. (1999).
Harmonic maps with noncontact boundary values.
{\it Proc. Amer. Math. Soc.} {\bf 127} (1999), 1231--41.
MR 99f:58048,  Zbl.399.15016.
\secno {\it S.13}.

Fillmore, J. P. (1977).
The fifteen-parameter conformal group.
{\it Internat. J. Theoret. Phys.} {\bf 16}, 937--63.
MR~58\#22417,  Zbl.399.15016.
\secno {\it S.9}.

Fuglede, B. (1978)*. Harmonic morphisms between Riemannian
manifolds. {\it Ann. Inst. Fourier (Grenoble)}, {\bf  28} (2), 107--44.  
MR~80h:58023,  Zbl.369.53044.
\secno{S.6}.

Fuglede, B. (1996)*. Harmonic morphisms between semi-Riemannian
manifolds. {\it Ann. Acad. Sci. Fennicae},  {\bf 21}, 31--50.
MR~97i:58035,  Zbl.84753013.
\secno{S.6}.

Harvey, F. R., Lawson, H. B. (1975).
On boundaries of complex analytic varieties. I.
{\it Ann. of Math.} (2) {\bf 102}, 223--90.
MR~54\#13130,  Zbl.317.32017.
\secno S.12.

Huggett, S. A., Mason, L. J., Tod, K. P., Tsou, S. T. and Woodhouse, N. M. J. (1998).
{\it The geometric universe. Science, geometry, and the work of Roger Penrose.} Papers
   from the Symposium on Geometric Issues in the Foundations of Science held in honor of the 65th
   birthday of Sir Roger Penrose at St. John's College, University of Oxford, Oxford, June 1996.
Oxford University Press, Oxford.
MR~99b:00015,  Zbl.890.00046.  
\secno {\it S.9}.

Ishihara, T. (1979)*.  A mapping of Riemannian manifolds which
preserves harmonic functions.  {\it J. Math. Kyoto Univ.}, {\bf  19},
215--29.  
MR~80k:58045,  Zbl.421.31006.
\secno{S.6}.

Jacobowitz, H.  (1990).
{\it An introduction to CR structures}. 
Mathematical Surveys and Monographs, 32. 
American Mathematical Society, Providence, RI.
MR~93h:32023,  Zbl.712.32001.
\secno S.12~{\it S.12}.

Kerr, R. P, (1963).
Gravitational field of a spinning mass as an example of algebraically special
   metrics. {\it Phys. Rev. Lett.} {\bf 11}, 237--8. MR~27 \#6594, Zbl.0112.21904
\secno {S.11}.

LeBrun, C.  (1983). 
Spaces of complex null geodesics in complex-Riemannian geometry. 
{\it Trans. Amer. Math. Soc.} {\bf 278}, 209--231.
MR~84e:32023.  Zbl.562.53018.
\secno S.4,~S.9.

LeBrun, C. R.  (1984).
Twistor CR manifolds and three-dimensional conformal geometry. 
{\it Trans. Amer. Math. Soc.} {\bf 284}, no. 2, 601--616.
MR~86m:32033.  Zbl.513.53006.
\secno {\it S.9},~{\it S.12},~{\it S.13}.

Lewandowski, J., Nurowski, P. and Tafel, J.  (1990).
Einstein's equations and realizability of CR manifolds.
{\it Classical Quantum Gravity} {\bf 7}, no. 11, L241--L246.
MR~91j:32021,    Zbl.714.53047.
\secno {\it S.9}.

Li, P. and Tam, L.-F. (1991)
The heat equation and harmonic maps of complete manifolds.
{\it Invent. Math.} {\bf 105}, 1--46.
MR~93e:58039,  Zbl.748.58006.
\secno {\it S.13}.

Li, P. and Tam, L.-F. (1993)
Uniqueness and regularity of proper harmonic maps.
{\it Ann. of Math.} (2) {\bf 137}, 167--201.
MR~93m:58027,  Zbl.776.58010.
\secno {\it S.13}.

Li, P., Tam, L.-F. and Wang, J.  (1995).
Harmonic diffeomorphisms between Hadamard manifolds.
{\it Trans. Amer. Math. Soc.}, {\bf 347}, 3645--58.  
MR~95m:58044,  Zbl.855.58020.
\secno {\it S.13}.

Mason, L.~J. and Hughston, L.~P.  (ed.)  (1990).
{\it Further advances in twistor theory Volume I:
The Penrose transform and its applications}.
Pitman Research Notes in Mathematics Series, vol.~231.
Longman Scientific and Technical, Harlow;
Wiley, New York.
MR~92d:32040,  Zbl.693.53022.
\secno {\it S.9}.

Mason, L.~J. and Woodhouse, N. M. J. (1996).
{\it Integrability, self-duality, and twistor theory}.
Oxford University Press, Oxford.
MR 98f:58002, Zbl.856.58002.
\secno {\it S.9}.

Mason, L.~J., Hughston, L.~P. and Kobak, P.~Z. (ed.) (1995).
{\it Further advances in twistor theory Volume II:
Integrable systems, conformal geometry and gravitation}.
Pitman Research Notes in Mathematics Series, vol.~232.
Longman Scientific and Technical, Harlow.
MR~97m:83067,  Zbl.828.53070.
\secno {\it S.9}.

Mason, L.~J., Hughston, L.~P., Kobak P.~Z. and Pulverer, K. (ed.) (2001).
{\it Further advances in twistor theory  Volume III:
Curved twistor spaces}.
Research Notes in Mathematics, vol.~424.
Chapman and Hall/CRC, Boca Raton, FL.
MR~2002d:83063,  Zbl.968.53035.
\secno {\it S.9}.

Nurowski, P.  (1996). 
Optical geometries and related structures.
{\it J. Geom. Phys.} {\bf 18}, no. 4, 335--348.
MR~97g:83083,  Zbl.849.53014.
\secno {\it S.9}

Pambira, A. (2002)*.
{\it Harmonic maps and morphisms in semi-Riemannian and com\-plex-Riemannian geometry.}
PhD thesis. University of Leeds.
\secno {\it S.7}

Pambira, A. (2003)*.
Harmonic morphisms between degenerate semi-Riemannian manifolds,
Preprint, University of Leeds. 
{\tt http://xxx.soton.ac.uk/abs/math.DG/0303275}.
\secno {\it S.7}

Penrose, R.  (1967).  Twistor algebra.
{\it J. Math. Phys.}, {\bf 8}, 345--66.   
MR~35\#7657,  Zbl.163.22602.
\secno {\it S.9}~{\it S.11}.

Penrose, R.  (1969).   Solutions of the zero-rest-mass equations {\it J. Math. Phys.}, {\bf 10}, No. 1, 38 - 39.
\secno {\it S.11}.

Penrose, R. (1975). Twistor theory, its aims and achievements. 
{\it Quantum gravity.} An Oxford Symposium held at the Rutherford Laboratory, Chilton,
   February 15--16, 1974. Ed. C. J. Isham, R. Penrose and D. W. Sciama.
Clarendon Press, Oxford.  MR~56\#17537.
\secno {\it S.9}.

Penrose, R. and  Rindler, W. (1987).
{\it Spinors and space-time. Vol.~1. Two-spinor calculus
and relativistic fields}.
Cambridge Monographs on Mathematical Physics,
2nd edn.  (1st edn., 1984).
Cambridge University Press, Cambridge.
MR~88h:83009,  Zbl.602.53001.
\secno {\it S.9}.
  
Penrose, R. and Rindler, W. (1988).
{\it Spinors and space-time. Vol.~2.  Spinor and twistor
methods in space-time geometry}.
Cambridge Monographs on Mathematical Physics,
2nd edn.  (1st edn., 1986).
Cambridge University Press, Cambridge.
MR~89d:83010,  Zbl.591.53002.
\secno {\it S.9},

Robinson, I. and Trautman, A. (1986).
Cauchy--Riemann structures in optical geometry.
{\it Proceedings of the fourth Marcel Grossmann meeting on general relativity, Part
A, B (Rome, 1985)}, 317--324, North-Holland, Amsterdam.
MR~88d:83026,   
\secno {\it S.9}

Sato, H. and Yamaguchi, K.  (1989).   Lie contact manifolds.
{\it Geometry of manifolds (Matsumoto, 1988)}, 191--238.
Perspect. Math., 8, Academic Press, Boston, MA.
MR~91m:53027,  Zbl.705.53019.
\secno {\it S.12}.

Sims, B.~T. (1976).
{\it Fundamentals of topology}.
Macmillan, New York; Collier
Macmillan, London. \linebreak
MR~54\#3638,  Zbl.324.54002.
\secno S.8.

Sparling, G. (1975).
Homology and twistor theory.
{\it Quantum gravity.} 
An Oxford Symposium held at the Rutherford Laboratory, Chilton, February 15--16, 1974.
Ed. C. J. Isham, R. Penrose and D. W. Sciama. 
Clarendon Press, Oxford, 1975.
MR~56 \#17537.
\secno {\it S.9}.

Tanno, S.  (1992).
The standard CR structure on the unit tangent bundle.
{\it Tohoku Math. J.} (2) {\bf 44}, 535--543.
MR~93k:53033,   Zbl.779.53024.
\secno {\it S.12}.

Trautman, A.  (1985).
Optical structures in relativistic theories.
{\it The mathematical heritage of \'Elie Cartan (Lyon, 1984).}
Ast\'erisque 1985, Numero Hors Serie, 401--420.
MR~88a:83085,  Zbl.599.53055.
\secno {\it S.9}.

Ward, R.~S. and Wells, R.~O., Jr. (1990).
{\it Twistor geometry and field theory}.
Cambridge Monographs on Mathematical Physics.
Cambridge University Press, Cambridge.
MR~91b:32034,  Zbl.729.53068.  
\secno {\it S.9}.

Wilker, J.~B. (1986).  Inversive geometry and the Hopf fibration. {\it
Studia Sci. Math.  Hungar.}, {\bf 21}, 91--101.
MR~88h:57008,  Zbl.627.55009.
\secno {\it S.12}.

Wood, J.~C. (1992)*. Harmonic morphisms and Hermitian structures
on Einstein $4$-manifolds. {\it Internat. J. Math.}, {\bf  3},  415--39.  
MR~94a:58054,  Zbl.763.53051.
\secno {\it S.8}

\end{document}